\documentclass[10pt]{amsart}

\newcommand{\comm}[1]{}

\numberwithin{equation}{section}
\renewcommand{\P}{\mathbb{P}}
\newcommand{\RR}{\mathfrak{R}}

\def\Ga{\Theta}

\def\Om{\mho}

\def\Up{\Upsilon}

\def\De{\Xi}

\def\di{\,|\,}

\def\be{\begin{equation}}
\def\ee{\end{equation}}

\def\ba{{\begin{align}}}
\def\ea{{\end{align}}}

\def\proj{{\operatorname{proj}}}

\def\Leb{{\mathrm{Leb}}}

\def\TV{{\mathcal{TV}}}
\def\VT{{\mathcal{VT}}}
\def\TF{{\mathcal{TF}}}

\def\e{{\Lambda}}
\def\n{{\mathrm{N}}}
\def\m{{\mathrm{m}}}
\def\M{{\mathrm{M}}}

\def\wde{{\widehat \Delta_r}}

\def\ess{{\mathrm{ess}}}
\def\red{{\mathrm{red}}}

\def\mm{{\underline m}}

\def\floor{\mathrm{floor}}
\def\roof{\mathrm{roof}}
\def\lat{\mathrm{lat}}

\newcommand\ssigma{\mathfrak{S}^0}
\newcommand\sssigma{\mathfrak{S}}

\newcommand{\SO}{{\mathrm {SO}}}
\newcommand{\GL}{{\mathrm {GL}}}
\newcommand{\SL}{{\mathrm {SL}}}

\newtheorem*{mainthm}{Main Theorem}

\newtheorem{thm}{Theorem}[section]
\newtheorem{cor}[thm]{Corollary}
\newtheorem{lem}[thm]{Lemma}
\newtheorem{lemma}[thm]{Lemma}

\newtheorem{prop}[thm]{Proposition}
\theoremstyle{remark}
\newtheorem{rem}[thm]{Remark}

\theoremstyle{definition}
\newtheorem{definition}[thm]{Definition}

\newcommand{\dist}{\operatorname{dist}}
\newcommand{\id}{\operatorname{id}}

\renewcommand{\AA}{{\mathcal A}}
\newcommand{\BB}{{\mathcal B}}
\newcommand{\CC}{{\mathcal C}}
\newcommand{\DD}{{\mathcal D}}
\newcommand{\EE}{{\mathcal E}}

\newcommand{\FF}{{\mathcal F}}

\newcommand{\HH}{{\mathcal H}}

\newcommand{\MM}{{\mathcal M}}

\newcommand{\QQ}{{\mathcal Q}}

\newcommand{\C}{{\mathbb C}}

\newcommand{\N}{{\mathbb N}}
\newcommand{\Q}{{\mathbb Q}}
\newcommand{\R}{{\mathbb R}}

\newcommand{\Z}{{\mathbb Z}}

\newcommand{\Cb}{C_0}
\newcommand{\Cd}{C_1}
\newcommand{\Ce}{C_2}
\newcommand{\eb}{\epsilon_0}
\newcommand{\norm}[1]{\left\| #1 \right\|}
\DeclareMathOperator{\dLeb}{dLeb}
\newcommand{\dd}{\, {\rm d}}
\newcommand{\tq}{\, :\, }
\newcommand{\moins}{\backslash}
\DeclareMathOperator{\sys}{sys} \DeclareMathOperator{\Var}{Var}
\newcommand{\Ca}{C_3}
\newcommand{\betaa}{\beta_0}
\newcommand{\Cc}{C_4}
\renewcommand{\hat}{\widehat}

\DeclareMathOperator{\type}{type} \DeclareMathOperator{\Lip}{Lip}

\DeclareMathOperator{\diam}{diam}

\renewcommand{\ln}{\log}

\begin{document}

\title{Exponential mixing for the Teichm\"{u}ller flow}
\date{\today}
\author{Artur Avila, S\'{e}bastien Gou\"{e}zel and Jean-Christophe Yoccoz}

\address{
CNRS UMR 7599, Laboratoire de Probabilit\'{e}s et Mod\`{e}les Al\'{e}atoires,
Universit\'{e} Pierre et Marie Curie, Bo{\^\i}te Postale 188,
75252 Paris Cedex 05, France.
}
\email{artur@ccr.jussieu.fr}

\address{
CNRS UMR 6625, IRMAR, Universit\'e de Rennes 1, Campus de Beaulieu,
35042 Rennes Cedex, France.
}
\email{sebastien.gouezel@univ-rennes1.fr}

\address{
Coll\`{e}ge de France, 3 rue d'Ulm, 75005 Paris, France.
}
\email{jean-c.yoccoz@college-de-france.fr}

\begin{abstract}

We study the dynamics of the Teichm\"{u}ller flow in the moduli space of
Abelian differentials (and more generally, its restriction to any connected
component of a stratum).  We show that the (Masur-Veech)
absolutely continuous invariant probability measure is exponentially
mixing for the class of H\"{o}lder observables.  A geometric consequence is
that the $\SL(2,\R)$ action in the moduli space has a spectral gap.

\end{abstract}

\maketitle

\setcounter{tocdepth}{1}\tableofcontents

\section{Introduction}

Let $\MM_g$ be the moduli space of non-zero Abelian differentials on a
compact Riemann surface of genus $g \geq 1$.  Alternatively, $\MM_g$ can be
seen as the moduli space of translation surfaces of genus $g$: outside the
zero set of an Abelian differential $\omega$
there are preferred local charts where $\omega=dz$, and the coordinate
changes of those charts are translations.  Let $\MM^{(1)}_g \subset \MM_g$
denote the subspace of surfaces with normalized area $\int |\omega|^2=1$.

By postcomposing the preferred
charts with an element of $\GL(2,\R)$ one obtains another translation
structure: this gives a natural $\GL(2,\R)$
action on $\MM_g$.  The $\SL(2,\R)$ action preserves $\MM^{(1)}_g$.
The Teichm\"{u}ller flow on $\MM_g$ is defined as the
diagonal action of $\SL(2,\R)$:
$\TF_t=\left (\begin{matrix}e^t&0\\0&e^{-t}\end{matrix} \right ):\MM_g
\to \MM_g$.

The space $\MM_g$ is naturally
stratified: given an unordered list $\kappa=(\kappa_1,...,\kappa_s)$ of
positive integers with $\sum (\kappa_i-1)=2g-2$, 
we let $\MM_{g,\kappa}$ be the
space of Abelian differentials whose zeroes have order
$\kappa_1-1,\dots,\kappa_s-1$.  The strata are obviously invariant by the
$\GL(2,\R)$ action.

The strata $\MM_{g,\kappa}$ are not necessarily
connected (a classification of connected components is given in \cite {KZ}).
Let $\CC$ be a connected component of some stratum $\MM_{g,\kappa}$, and let
$\CC^{(1)}=\CC \cap \MM^{(1)}_g$.  It has
a natural structure of an analytic variety, and hence a natural Lebesgue
measure class.  By the
fundamental work of Masur \cite {M} and Veech \cite {V1}, there exists a
unique probability measure $\nu_{\CC^{(1)}}$ on $\CC^{(1)}$
which is equivalent to
Lebesgue measure, invariant by the Teichm\"{u}ller flow, and ergodic.
Veech later showed in \cite {V2} that
$\nu_{\CC^{(1)}}$ is actually mixing, meaning that for any observables
$f,g \in L^2(\nu_{\CC^{(1)}})$ one has
\be \label {mix}
\int_{\CC^{(1)}} f \circ \TF_t(x) g(x) \dd\nu_{\CC^{(1)}}(x)
-\int_{\CC^{(1)}}
f \dd\nu_{\CC^{(1)}} \int_{\CC^{(1)}} g \dd \nu_{\CC^{(1)}}\to 0.
\ee
In this paper we are concerned with the speed of mixing of the Teichm\"{u}ller
flow, that is, the rate of the convergence in (\ref {mix}), for a suitable
class of observables.

\begin{mainthm}

The Teichm\"{u}ller flow (restricted to any connected component of any
stratum of the moduli space of Abelian differentials) is exponentially
mixing for H\"{o}lder observables.

\end{mainthm}

The complete formulation of this result, specifying in particular what is
understood by a H\"older observable in this non-compact setting,
is given in Theorem \ref{decay_teichmuller}.

Previously it had been shown by Bufetov \cite {Bu} that the Central Limit
Theorem holds for the Teichm\"{u}ller flow (for suitable classes of
observables).  Though he did not obtain rates of mixing for the
Teichm\"{u}ller flow itself, he did obtain {\it stretched exponential}
estimates for a related discrete time transformation (the Zorich
renormalization algorithm for interval exchange transformations).  In this
paper we will also work with a discrete time transformation, though not
directly with the Zorich renormalization.

This paper has two main parts: first we obtain exponential recurrence
estimates, and then, using ideas first introduced by Dolgopyat \cite {Do}
and developed in \cite {BV},
we obtain exponential mixing.  The proof of exponential recurrence
uses an ``induction on the complexity'' scheme.
Intuitively, the dynamics at ``infinity'' of the Teichm\"{u}ller flow
can be partially described by the dynamics in simpler (lower dimensional)
connected components of strata, and we obtain estimates by induction all the
way from the simplest of the cases.  A simpler version of this scheme was
used to show some combinatorial richness of the Teichm\"{u}ller flow in the
proof of the Zorich-Kontsevich conjecture \cite {AV}.  The recurrence
estimates thus obtained are close to optimal.

It should be noted that our work does not use the
$\SL(2,\R)$ action for the estimates, and can be used to obtain new proofs
of some previously known results which used to depend on the $\SL(2,\R)$
action.  In the other direction, however, it was pointed out to us by
Bufetov that our main theorem has an important new corollary for the
$\SL(2,\R)$ action.  It regards the nature of the corresponding
unitary representation of $\SL(2,\R)$ on the space
$L^2_0(\nu_{\CC^{(1)}})$ of $L^2$ zero-average functions.

\begin{cor} \label {isolated}

The action of $\SL(2,\R)$ on $L^2_0(\nu_{\CC^{(1)}})$ has a spectral
gap.

\end{cor}

The notion of spectral gap and the derivation of the corollary from the
Main Theorem are discussed in Appendix \ref {sg}.

\begin{rem}

Exponential recurrence estimates for the Teichm\"{u}ller flow were
first obtained by Athreya \cite {At}, who used the $\SL(2,\R)$ action to
prove them for some large compact sets (which are, in particular,
$\SO(2,\R)$ invariant).  Our work allows us to obtain exponential recurrence
for certain smaller compact sets, for which the first return map has good
hyperbolic properties.  Bufetov has independently obtained a proof
of exponential recurrence estimates for such small compact sets, using the
method of \cite {Bu}.  Those estimates, while non-optimal, are enough to
obtain exponential mixing using the remainder of our argument.

We should also point out that recurrence estimates are
often useful in statistical arguments in a very practical sense.
For instance, the proof of typical weak mixing in \cite {AF} can be
made more transparent using such estimates.

\end{rem}

\comm{
\subsection{Outline}

In \S \ref {teich}, we give
a more detailed presentation of the Teichm\"{u}ller flow, which
will enable us to give the precise formulation of our main results.

In \S \ref {veech},
we then turn to the Veech flow, which provides a more combinatorially
accessible model for the Teichm\"{u}ller flow.  We give a description of its
basic combinatorics and its relation to interval exchange transformations.

In \S \ref {reduc}, we show that the Teichm\"{u}ller flow can be seen as a
hyperbolic extension of a suspension semiflow over an
(uniformly) expanding Markov system.  After stating a general result on
exponential mixing for such systems, we are able to reduce our main results
to certain recurrence estimates for the Veech flow.

In \S \ref {main estimate}, we provide an almost optimal distortion estimate
for the renormalization algorithm of interval exchange transformations, and
in \S \ref {recurrence estimat} we use the distortion estimate to prove the
recurrence estimates.

In \S \ref {expanding semiflows}, we carry out the analysis of exponential
mixing for suspension over expanding Markov systems, closely following \cite
{BV}, and in \S \ref {hyperbolic flows}, we consider the case of hyperbolic
extensions.

In Appendix \ref {appendixc}, we give a simple proof of a non-optimal
distortion estimate.
}

{\bf Acknowledgments:}  We thank Nalini Anantharaman, Sasha Bufetov, 
Giovanni Forni and Viviane Baladi for several discussions.

\section{Statements of the results}

\newcommand{\expansion}{\kappa}

\subsection{Exponential mixing for excellent hyperbolic semiflows}

\label{par:define_exp}

To prove exponential decay of correlations for the Teichm\"{u}ller flow,
we will show that this flow can be reduced to an abstract flow with
good hyperbolic properties. In this paragraph, we describe some
assumptions under which such a flow is exponentially mixing.

By definition, a Finsler manifold is a smooth manifold endowed with
a norm on each tangent space, which varies continuously with the
base point.

\begin{definition}
A \emph{John domain} $\Delta$ is a finite dimensional connected
Finsler manifold, together with a measure $\Leb$ on $\Delta$, with
the following properties:
\begin{enumerate}
\item For $x,x' \in \Delta$, let  $d(x,x')$ be the infimum
of the length of a $C^1$ path contained in $\Delta$ and
joining $x$ and $x'$. For this distance, $\Delta$ is bounded and
there exist constants $\Cb$ and
$\eb$ such that, for all $\epsilon<\eb$, for all $x\in \Delta$,
there exists $x'\in \Delta$ such that $d(x,x') \leq \Cb \epsilon$
and such that the ball $B(x',\epsilon)$ is compactly contained in
$\Delta$.
\item The measure $\Leb$ is a fully supported finite
measure on $\Delta$, satisfying the following inequality: for all
$C>0$, there exists $A>0$ such that, whenever a ball $B(x,r)$ is
compactly contained in $\Delta$, $\Leb(B(x,Cr))\leq A \Leb(B(x,r))$.
\end{enumerate}
\end{definition}
For example, if $\Delta$ is an open subset of a larger manifold,
with compact closure, whose boundary is a finite union of smooth
hypersurfaces in general position, and 
$\Leb$ is obtained by restricting to $\Delta$ a smooth
measure defined in a neighborhood of $\overline{\Delta}$, then
$(\Delta,\Leb)$ is a John domain. 

\begin{definition} \label {markovmap}
Let $L$ be a finite or countable set, 
let $\Delta$ be a John domain, and let $\{\Delta^{(l)}\}_{l\in L}$ be a
partition into open sets of a full measure subset of $\Delta$. A map $T:
\bigcup_{l} \Delta^{(l)} \to \Delta$ is a \emph{uniformly expanding
Markov map} if
\begin{enumerate}
\item For each $l$, $T$ is a $C^1$ diffeomorphism between
$\Delta^{(l)}$ and $\Delta$, and there exist constants
$\expansion>1$ (independent of $l$) and $C_{(l)}$ such that, for all
$x\in \Delta^{(l)}$ and all $v\in T_x \Delta$, $\expansion \norm{v}
\leq \norm{ DT(x)\cdot v} \leq C_{(l)} \norm{v}$.
\item Let $J(x)$ be the inverse of the Jacobian of $T$ with respect to $\Leb$.
Denote by $\HH$ the set of inverse branches of $T$. The function
$\log J$ is $C^1$ on each set $\Delta^{(l)}$ and 
there exists $C>0$ such that, for all $h\in \HH$, $\norm{ D((\log
J)\circ h)}_{C^0(\Delta)}\leq C$.
\end{enumerate}
\end{definition}
Such a map $T$ preserves a unique absolutely continuous measure
$\mu$. Its density is bounded from above and from below and is
$C^1$. This measure is ergodic and even mixing
(see e.g.\ \cite{aaronson:book}). Notice that $\Leb$ is not assumed to
be absolutely continuous with respect to Lebesgue measure. Although
this will be the case in most applications, this definition covers
also e.g.\ the case of maximum entropy measures when $L$ is finite
(in which case $\log J$ is constant, which yields $D( (\log J)\circ h)=0$).

\begin{definition}
\label{def:goodroof}
Let $ T:\bigcup_{l} \Delta^{(l)} \to \Delta$ be
a uniformly expanding Markov map on a John domain. A function $r:
\bigcup \Delta^{(l)} \to \R_+$ is a \emph{good roof function} if
\begin{enumerate}
\item There exists $\epsilon_1>0$ such that $r \geq \epsilon_1$.
\item There exists $C>0$ such that, for all $h\in \HH$, $\norm{
D(r\circ h)}_{C^0}\leq C$.
\item It is not possible to write $r=\psi + \phi\circ T-\phi$ 
on $\bigcup \Delta^{(l)}$,
where $\psi:\Delta \to \R$ is constant on each set $\Delta^{(l)}$
and $\phi:\Delta\to \R$ is $C^1$.
\end{enumerate}
\end{definition}
If $r$ is a good roof function for $T$, we will write
$r^{(n)}(x)=\sum_{k=0}^{n-1}r(T^k x)$.

\begin{definition}
A good roof function $r$ as above has
\emph{exponential tails} if there exists $\sigma_0>0$ such that
$\int_{\Delta} e^{\sigma_0 r} \dLeb <\infty$.
\end{definition}

If $\widehat \Delta$ is a Finsler manifold, we will denote by
$C^1(\widehat{\Delta})$ the set of functions $u:\widehat{\Delta}\to
\R$ which are bounded, continuously differentiable, and such that
$\sup_{x\in \widehat{\Delta}} \norm{Du(x)}< \infty$. Let
  \begin{equation}
  \label{define_C1}
  \norm{u}_{C^1(\widehat{\Delta})} = \sup_{x\in \widehat{\Delta}} |u(x)| + \sup_{x\in
  \widehat{\Delta}} \norm{Du(x)}
  \end{equation}
be the corresponding norm.

\begin{definition} \label {skew-product map}
Let $T:\bigcup_{l} \Delta^{(l)} \to \Delta$ be a uniformly expanding
Markov map, preserving an absolutely continuous
measure $\mu$. An \emph{hyperbolic
skew-product} over $T$ is a map $\widehat{T}$ from a dense open subset
of a bounded connected
Finsler
manifold $\widehat{\Delta}$, to $\widehat{\Delta}$, satisfying the following
properties:
\begin{enumerate}
\item There exists a continuous
map $\pi : \widehat{\Delta} \to \Delta$
such that $T\circ \pi = \pi \circ \widehat{T}$ whenever both members
of this equality are defined.
\item There exists a probability measure $\nu$ on $\widehat{\Delta}$,
giving full mass to the domain of definition of $\widehat{T}$, which is
invariant under $\widehat{T}$.
\item There exists a family of probability measures $\{\nu_x\}_{x\in
\Delta}$ on $\widehat{\Delta}$ which is a disintegration of $\nu$ over
$\mu$ in the
following sense: $x\mapsto \nu_x$ is measurable, $\nu_x$ is supported
on $\pi^{-1}(x)$ and, for all measurable set $A \subset
\widehat{\Delta}$, $\nu(A)=\int \nu_x(A) \dd\mu(x)$.

Moreover, this disintegration satisfies the following property: 
there exists a constant $C>0$ such that,
for any open subset $O\subset \bigcup \Delta^{(l)}$, for any
$u\in C^1(\pi^{-1}(O))$, the
function $\bar u : O \to \R$ given by $\bar u(x)=\int u(y)
\dd\nu_x(y)$ belongs to $C^1(O)$ and
satisfies the inequality
  \begin{equation}
  \sup_{x\in O} \norm{D\bar u(x)}
  \leq C \sup_{y\in \pi^{-1}(O)} \norm{Du(y)}.
  \end{equation}
\item
There exists $\expansion>1$ such that, for all $y_1,y_2 \in
\widehat{\Delta}$ with $\pi(y_1)=\pi(y_2)$, holds
  \begin{equation}
  d(\widehat{T} y_1, \widehat{T} y_2) \leq \expansion^{-1} d(y_1,y_2).
  \end{equation}
\end{enumerate}
\end{definition}

Let $\widehat{T}$ be an hyperbolic skew-product over a uniformly
expanding Markov map $T$. Let $r$ be a good roof function for $T$,
with exponential tails. It is then possible to define a space
$\widehat{\Delta}_r$ and a semiflow $\widehat{T}_t$ over
$\widehat{T}$ on $\widehat{\Delta}$, using the roof function $r\circ
\pi$, in the following way.
 Let $\widehat{\Delta}_r =
\{(y,s) \tq y\in \bigcup_l \widehat{\Delta}_l, 0\leq s<r(\pi y )\}$.
For almost
all $y\in \widehat{\Delta}$, all $0\leq s<r(\pi y)$ and all $t\geq 0$, there
exists a unique $n\in \N$ such that $r^{(n)}(\pi y)\leq
t+s<r^{(n+1)}(\pi y)$. Set $\widehat{T}_t (y,s)=(\widehat{T}^n y,
s+t -r^{(n)}(\pi y))$. This is a semiflow defined almost everywhere
on $\widehat{\Delta}_r$. It preserves the probability measure
$\nu_r=\nu \otimes \Leb/ (\nu\otimes \Leb)(\widehat{\Delta}_r)$.
Using the canonical Finsler metric on $\widehat{\Delta}_r$, namely
the product metric given by 
$\|(u,v)\|:=\|u\|+\|v\|$, we define the space $C^1(\widehat{\Delta}_r)$ as
in \eqref{define_C1}.  Notice that $\widehat \Delta_r$ is not connected, and
the distance between points in different connected components is infinite.

\begin{definition}
A semiflow $\widehat{T}_t$ as above is called an \emph{excellent
hyperbolic semiflow}.
\end{definition}

The main motivations for this definition are that the Teichm\"{u}ller
flow is isomorphic to an excellent hyperbolic semiflow  -- the proof
of this isomorphism will take a large part of this article -- and
that such a flow has exponential decay of correlations:

\begin{thm}
\label{main_thm_hyperbolic} Let $\widehat{T}_t$ be an excellent
hyperbolic semi-flow on a space $\widehat{\Delta}_r$, preserving the
probability measure $\nu_r$. There exist constants $C>0$ and
$\delta>0$ such that, for all functions $U,V\in
C^1(\widehat{\Delta}_r)$, for all $t\geq 0$,
  \begin{equation}
  \left| \int U\cdot V \circ \widehat{T}_t \dd\nu_r -\left( \int U
  \dd\nu_r\right) \left( \int V \dd\nu_r \right) \right| \leq C
  \norm{U}_{C^1} \norm{V}_{C^1}e^{-\delta t}.
  \end{equation}
\end{thm}

We will see the consequences of this theorem in the next sections.
The proof of Theorem \ref{main_thm_hyperbolic} will be deferred to
Sections \ref{par:define_expanding} and \ref{par:define_hyperbolic}.

\subsection{The Teichm\"{u}ller flow}

\subsubsection{Teichm\"{u}ller space, moduli space and the Teichm\"{u}ller
flow}

Let $g\in \N^*$ and $s\in \N^*$ be positive integers. Take $M$ a
compact orientable $C^\infty$ surface of genus $g$, and let
$\Sigma=\{A_1,\dots, A_s\}$ be a subset of $M$. Let
$\kappa=(\kappa_1,\dots, \kappa_s)\in (\mathbb{N}^*)^s$ be such that $\sum
(\kappa_i -1)=2g-2$.

A \emph{translation structure} on $(M,\Sigma)$ with singularities
type $\kappa$ is an atlas on $M \moins \Sigma$ for which the
coordinate changes are translations, and such that each singularity
$A_i$ has a neighborhood which is isomorphic to the $\kappa_i$-fold
covering of a neighborhood of $0$ in $\R^2 \moins \{0\}$. The
Teichm\"{u}ller space $\QQ_{g,\kappa}=\QQ( M, \Sigma, \kappa)$ 
is the set of such
structures modulo isotopy rel.\ $\Sigma$. It has a canonical
structure of manifold.

Let us describe this manifold structure by introducing charts
through the \emph{period map} $\Theta$. Let $\xi$ be a translation
structure on $(M,\Sigma)$. If $\gamma \in C^0( [0,T], M)$ is a path
on $M$, then it is possible to lift it in $\R^2$, starting from $0$:
this lifting is possible locally outside of the singularities, and
the local form of the translation structure close to the
singularities implies that this lifting is also possible at the
singularities. Taking the value of the lifting at $T$, we get a
\emph{developing map}
  \begin{equation}
  D_\xi : C^0( [0,T], M) \to \R^2.
  \end{equation}
This map yields a linear map $H_1(M,\Sigma ; \Z) \to
\R^2$, i.e., an element of $H^1(M,\Sigma; \R^2)$. It is invariant
under isotopy rel.\ $\Sigma$. Hence, it defines a map $\Theta :
\QQ(M,\Sigma, \kappa) \to H^1(M,\Sigma; \R^2)$.

This map is a local diffeomorphism for the canonical manifold
structure of $\QQ(M,\Sigma, \kappa)$, and gives in particular local
coordinates. It even endows $\QQ(M,\Sigma,\kappa)$ with a complex
affine manifold structure.

A translation structure on $(M,\Sigma)$ defines a volume form on $M
\moins \Sigma$ (namely, the pullback of the standard volume form on
$\R^2$ by any translation chart). The manifold $M$ has finite area
for this volume form. Let $\QQ^{(1)}(M,\Sigma,\kappa)$ be the smooth
hypersurface of $\QQ(M,\Sigma, \kappa)$ given by area $1$
translation structures.

The space $H^1(M,\Sigma ; \R^2)$ has a standard volume form (the
Lebesgue form giving covolume $1$ to the integer lattice). Pulling
it back locally with $\Theta$, we obtain a smooth measure $\mu$ on
$\QQ(M,\Sigma,\kappa)$. It induces a smooth measure $\mu^{(1)}$ on
the hypersurface $\QQ^{(1)}(M,\Sigma, \kappa)$.

The group $\SL(2,\R)$ acts on $\QQ(M,\Sigma,\kappa)$ by
postcomposition in the charts. It preserves the hypersurface
$\QQ^{(1)}(M,\Sigma,\kappa)$ and leaves invariant the measures $\mu$
and $\mu^{(1)}$. In particular, the action of $\TF_t:= \left(
\begin{array}{cc} e^t & 0 \\0& e^{-t}\end{array}\right)$ is a
measure preserving flow, called the \emph{Teichm\"{u}ller flow}.

The modular group of $(M,\Sigma)$ is the group of diffeomorphisms of
$M$ fixing $\Sigma$, modulo isotopy rel.\ $\Sigma$. It acts on the
Teichm\"{u}ller space $\QQ(M,\Sigma, \kappa)$. The quotient is denoted
by $\MM_{g,\kappa}=\MM(M,\Sigma, \kappa)$ and is called 
the \emph{moduli space}.
The action of the modular group on $\QQ(M,\Sigma,\kappa)$ is proper
and faithful, but it is not free. Hence, $\MM(M,\Sigma, \kappa)$ has
a complex affine orbifold structure.

Since the action of the modular group preserves the measure $\mu$
and the hypersurface $\QQ^{(1)}$, we also obtain a measure $\nu$ on
the moduli space, as well as a codimension $1$ hypersurface
$\MM^{(1)}(M,\Sigma, \kappa)$ of area $1$ translation structures,
and a measure $\nu^{(1)}$ on it. Moreover, the action of $\SL(2,\R)$
commutes with the action of the modular group, whence $\SL(2,\R)$
still acts on $\MM(M,\Sigma, \kappa)$ and $\MM^{(1)}(M,\Sigma,
\kappa)$, preserving respectively $\nu$ and $\nu^{(1)}$. In
particular, the action of $\TF_t$ defines a flow on $\MM(M,\Sigma,
\kappa)$, that we still call the Teichm\"{u}ller flow.

\begin{thm}[Masur, Veech]
\label{thm_masur_veech} The measure $\nu^{(1)}$ has finite mass.
Moreover, on each connected component of $\MM^{(1)}(M,\Sigma,
\kappa)$, the Teichm\"{u}ller flow is ergodic, and even mixing.
\end{thm}

Our goal in this paper is to estimate the speed of mixing of the
Teichm\"{u}ller flow. Our estimates will in particular give a new proof
of Theorem \ref{thm_masur_veech}.

\subsubsection{A Finsler metric on the Teichm\"{u}ller space} \label
{metricfinsler}

For a general dynamical system, the exponential decay of
correlations usually only holds at best for sufficiently regular
functions. In our case, ``regular'' will mean H\"{o}lder continuous, for some
natural metric. This metric will be a Finsler metric on the
Teichm\"{u}ller space, invariant under the action of the modular group.

Let $\xi$ be a translation structure on $(M,\Sigma)$ with
singularities type $\kappa$. The \emph{saddle connections} of $\xi$
are the unit speed geodesic paths $\gamma : [0,T] \to M$ such that
$\gamma^{-1}(\Sigma)= \{0,T\}$. Equivalently, these are straight
lines (for the translation structure) connecting two singularities,
and without singularity in their interiors. If $\gamma$ is a saddle
connection, then $D_\xi(\gamma)$ is a complex number measuring the
holonomy of the translation structure along $\gamma$. If $[\gamma]$
is the class of $\gamma$ in $H_1(M,\Sigma; \Z)$, then
$D_\xi(\gamma)= \Theta(\xi) ([\gamma])$ by definition of $\Theta$.

The saddle connections define in particular elements of
$H_1(M,\Sigma ; \Z)$. They are invariant under isotopy, and depend
only on the class of $\xi$ in $\QQ(M,\Sigma,\kappa)$. The following
lemma is well known (see e.g.\ \cite{EM}).
\begin{lem}
\label{lem:gluing}
Any translation surface $\xi$ admits a triangulation whose vertices
are the singularities $\Sigma$ and whose edges are saddle
connections. In particular, the saddle connections generate the
homology $H_1(M,\Sigma; \R)$.
\end{lem}

\begin{prop}
Let $q \in \QQ(M,\Sigma,\kappa)$, and let $\xi$ be a translation
surface representing $q$. Let $\{\gamma_n\}$ be the set of its
saddle connections. Define a function $\norm{\cdot}_q$ on
$H^1(M,\Sigma; \C)$ by
  \begin{equation}
  \norm{ \omega}_q = \sup_{n\in \N} \left|
  \frac{\omega([\gamma_n])}{\Theta(q)([\gamma_n])}\right|.
  \end{equation}
This function defines a norm on $H^1(M,\Sigma; \C)$.
\end{prop}
\begin{proof}
Let $\norm{\cdot}$ be any norm on $H_1(M,\Sigma; \R)$. We will prove
the existence of $C>0$ such that, for any saddle connection
$\gamma$, $C^{-1} \norm{[\gamma]} \leq | \Theta(q)([\gamma])| \leq C
\norm{[\gamma]}$.
Since the
saddle connections generate the homology, this will easily imply the
result of the proposition.

Since $\gamma \mapsto \Theta(q)([\gamma])$ is linear, the inequality
$ | \Theta(q)([\gamma])| \leq C \norm{[\gamma]}$ is trivial. For the
converse inequality, let $L>0$ be such that any point of $M$ can be
joined to a point of $\Sigma$ by a path of length at most $L$. The
inequality $C^{-1} \norm{\gamma} \leq |\Theta(q)([\gamma])|$ is
trivial for the (finite number of) saddle connections of length
$\leq L$. Consider now a saddle connection $\gamma$ with length
$\geq L$, and let $n\geq 2$ be such that $(n/2)L \leq
|\Theta(q)([\gamma])| \leq nL$. We can subdivide $\gamma$ in $n$
segments $[x_i,x_{i+1}]$ of length at most $L$. Joining each $x_i$
to a singularity, we obtain a decomposition in homology
$[\gamma]=\sum_{i=1}^n [ \gamma_i]$, where $\gamma_i$ is a path of
length at most $3L$. There exists a constant $C$ such that any such
path $\gamma_i$ satisfies $\norm{ [\gamma_i]} \leq C$, and we obtain
$\norm{ [\gamma]} \leq n C \leq \frac{2C}{L} |
\Theta(q)([\gamma])|$.
\end{proof}

\begin{prop}
The map from $\QQ(M,\Sigma,\kappa)$ to the set of norms on
$H^1(M,\Sigma; \C)$ given by $q\mapsto \norm{\cdot}_q$ is continuous.
\end{prop}
\begin{proof}
Let $\epsilon>0$. By compactness of the unit ball, there exists a
finite number of saddle connections $\gamma_1,\dots, \gamma_N$ such
that, for any $\omega \in H^1(M, \Sigma ; \C)$,
  \begin{equation}
  \norm{\omega}_q \leq
  (1+\epsilon) \sup_{1\leq n\leq N} \left|
  \frac{\omega([\gamma_n])}{\Theta(q)([\gamma_n])}\right|.
  \end{equation}
If $q'$ is close enough to $q$, the saddle connections $\gamma_i$
survive in $q'$, and we get
  \begin{equation}
  \norm{\omega}_{q'} \geq \sup_{1 \leq n \leq N}
  \left|
  \frac{\omega([\gamma_n])}{\Theta(q')([\gamma_n])}\right|
  \geq (1-\epsilon) \sup_{1 \leq n \leq N}
  \left|
  \frac{\omega([\gamma_n])}{\Theta(q)([\gamma_n])}\right|
  \geq \frac{1-\epsilon}{1+\epsilon} \norm{\omega}_q.
  \end{equation}
For the converse inequality, we have to prove that the new saddle
connections appearing in $q'$ do not increase the norm too much. Let
$\xi$ be a translation surface representing $q$. By Lemma 
\ref{lem:gluing}, $\xi$ is obtained by gluing a finite number of triangles
along some parallel edges. A translation surface $\xi'$ close to $\xi$ is
obtained by modifying slightly the sides of these
triangles in $\R^2$ and then gluing them along the same pattern. Hence, we get
a map $\phi_{\xi \xi'} : \xi \to \xi'$ which is affine in each
triangle of the triangulation. Moreover, if $\xi'$ is close enough
to $\xi$, the differential of $\phi_{\xi \xi'}$ is $\epsilon$-close
to the identity

Let $\gamma'$ be a saddle connection in $\xi'$. The path $\phi_{\xi
\xi'}^{-1}(\gamma')$ is a union of a finite number of segments in
$\xi$, and its length is at most $(1+\epsilon) |D_{\xi'}(\gamma')|$.
It is homotopic to a unique geodesic path $\gamma$ in $\xi$. This
path is a union of a finite number of saddle connections
$\gamma_1,\dots, \gamma_N$, with $\sum |D_\xi(\gamma_i)| \leq
(1+\epsilon)|D_{\xi'}(\gamma')|$. For $\omega \in H^1(M,\Sigma ;
\C)$,  we get
  \begin{equation}
  \left| \frac{ \omega( [\gamma'])}{ \Theta(q')([\gamma'])} \right|
  = \frac{ \left| \sum_{i=1}^N \omega( [ \gamma_i])
  \right|}{|D_{\xi'}(\gamma')| }
  \leq (1+\epsilon) \frac{ \sum_{i=1}^N |\omega( [ \gamma_i])|}
  { \sum_{i=1}^N | D_{\xi}( \gamma_i)|}
  \leq (1+\epsilon) \sup_{1\leq i \leq N} \frac{ |\omega( [
  \gamma_i])|}{ | D_{\xi}( \gamma_i)|}
  \leq (1+\epsilon) \norm{\omega}_q.
  \end{equation}
Hence, we obtain $\norm{\omega}_{q'}\leq (1+\epsilon) \norm{\omega}_q$.
\end{proof}

Since the tangent space of $\QQ(M,\Sigma,\kappa)$ is everywhere
identified through $\Theta$ with $H^1(M, \Sigma; \C)$, the norm
$\norm{\cdot}_q$ gives a Finsler metric on $\QQ(M,\Sigma,\kappa)$.
It defines a distance (which is infinite for points in different
connected components) on $\QQ(M,\Sigma,\kappa)$ as follows: the
distance between two points $x,x'$ is the infimum of the length
(measured with the Finsler metric) of a $C^1$ path joining $x$ and
$x'$.

Let $\sys : \QQ(M,\Sigma, \kappa) \to \R_+$ be the systole
function, i.e., the shortest length of a saddle connection. It is
bounded on $\QQ^{(1)}(M,\Sigma,\kappa)$.

\begin{lem}
\label{systole_Lpischitz}
The function $q\mapsto \log( \sys(q))$ is
$1$-Lipschitz on $\QQ(M,\Sigma, \kappa)$.
\end{lem}
\begin{proof}
We will prove that, for any $C^1$ path $\rho : (-1,1) \to \QQ(M,\Sigma,
\kappa)$ with $\rho(0)=q$ and $\rho'(0)=\omega \in H^1(M,\Sigma; \C)$
holds
  \begin{equation}
  \limsup_{t\to 0} \frac{| \log \sys( \rho(t)) - \log \sys(q)|}{|t|}
  \leq \norm{\omega}_q.
  \end{equation}
This will easily imply the result.

In a translation surface representing $q$, there is a finite number of
saddle connections $\gamma_1,\dots, \gamma_N$ with minimal length. For
small enough $t$, $\sys(\rho(t))=\min_{1\leq i \leq N} |\Theta(
\rho(t))([\gamma_i])|$. Moreover,
  \begin{equation}
  \log   |\Theta( \rho(t))([\gamma_i])| - \log(\sys(q))
  = \log \left | \frac{\Theta(q)([\gamma_i]) +
  t\omega([\gamma_i])+o(t)}{ \Theta(q)([\gamma_i])} \right|
  =
  t \Re\left(
  \frac{ \omega([\gamma_i])}{\Theta(q)([\gamma_i])} \right)+o(t).
  \end{equation}
Hence,
  \begin{equation}
  | \log \sys( \rho(t)) - \log \sys(q)| \leq |t| \max_{1\leq i \leq N}
  \left| \frac{\omega([\gamma_i])}{\Theta(q)([\gamma_i])} \right|+o(t)
  \leq |t| \norm{\omega}_q + o(t).
  \qedhere
  \end{equation}
\end{proof}

By construction, the norm $\norm{\cdot}_q$ is invariant under the
action of the modular group. As a consequence, the modular group
acts by isometries on $\QQ(M,\Sigma,\kappa)$. Hence, the distance on
$\QQ(M,\Sigma,\kappa)$ induces a distance on the quotient
$\MM(M,\Sigma,\kappa)$. It is Finsler outside of the singularities
of this orbifold. Notice that the systole is also invariant under
the modular group, and passes to the quotient. We will still denote
by $\sys$ this new function. The function $\log\circ \sys$
is still $1$-Lipschitz on $\MM(M,\Sigma,\kappa)$.

The systole plays an important role
in the topology of $\MM^{(1)}(M,\Sigma, \kappa)$ since, for all
$\epsilon>0$, the set $\{ q \in \MM^{(1)}(M,\Sigma, \kappa) \tq
\sys(q) \geq \epsilon\}$ is compact. To say it differently, a
sequence $q_n \in \MM^{(1)}(M,\Sigma,\kappa)$ diverges to infinity
if and only if $\sys(q_n)\to 0$.

\begin{cor}
The distance on $\QQ^{(1)}(M,\Sigma,\kappa)$ is complete.
\end{cor}
\begin{proof}
It is sufficient to prove the same statement in the quotient
$\MM^{(1)}(M,\Sigma,\kappa)$.
If $q_n$ is a Cauchy sequence in $\MM^{(1)}(M,\Sigma,\kappa)$, the sequence
$\log \sys(q_n)$ is also Cauchy by Lemma
\ref{systole_Lpischitz}. Hence, $\sys(q_n)$ is bounded away from
$0$. In particular, the sequence $q_n$ remains in a compact subset of
$\MM^{(1)}(M,\Sigma,\kappa)$, and converges to any of its cluster values.
\end{proof}

Any element $\omega\in H^1(M,\Sigma; \C)$ can be written uniquely as
$\omega=a+ib$ where $a,b\in H^1(M,\Sigma; \R)$. Let $\overline \omega=a-ib$.
In this notation, the
differential of the action of the Teichm\"uller flow is given by 
  \begin{equation}
  \left.\frac{\dd \TF_t(q)}{\dd t}\right|_{t=0}= \overline{ \Theta(q)}.
  \end{equation}
Hence, $\norm{ \left.\frac{\dd \TF_t(q)}{\dd t}\right|_{t=0}}_q \leq 1$.
In particular, the Teichm\"{u}ller flow satisfies $d(\TF_t(q),q)\leq
|t|$. The same inequality holds in the quotient space
$\MM(M,\Sigma,\kappa)$.

If $q\in \QQ(M,\Sigma,\kappa)$ and $\omega=a+ib \in H^1(M,\Sigma; \C)$
(identified through $\Theta$ with the tangent space of
$\QQ(M,\Sigma,\kappa)$ at $q$), then the differential of the
Teichm\"uller flow is given by
  \begin{equation}
  D \TF_t(q) \cdot \omega= e^t a+i e^{-t}b.
  \end{equation}
This implies the inequality
  \begin{equation}
  \label {e2t}
  e^{-2|t|} \|\omega\|_q\leq 
  \|D \TF_t(q) \cdot \omega\|_{\TF_t(q)}\leq e^{2|t|} \|\omega\|_q,
  \end{equation}
which corresponds to the 
classical fact that the extreme Lyapunov exponent of the
Teichm\"uller flow are $-2$ and $2$.  

\subsubsection{Exponential decay of correlations}

Let ${\CC^{(1)}}$ be a connected component of $\MM^{(1)}(M,\Sigma,\kappa)$.
It is an orbifold, and is endowed with a finite mass measure
$\nu_{\CC^{(1)}}$ (which we will assume to be normalized so that it is a
probability measure), and a distance $d_{\CC^{(1)}}$. The Teichm\"{u}ller
diagonal flow $\TF_t$ acts ergodically on ${\CC^{(1)}}$ and preserves the
measure $\nu_{\CC^{(1)}}$.

For $0<\alpha \leq 1$ and $f : {\CC^{(1)}} \to \R$, we will denote by
$\omega_\alpha(f,x)$ the local H\"{o}lder constant of $f$ at $x$, i.e.
  \begin{equation}
  \omega_\alpha(f,x)=\sup_{\substack{y\in
  B(x,1) \\ y\not=x}} \frac{ |f(y)-f(x)|}{d_{\CC^{(1)}}(y,x)^\alpha}.
  \end{equation}
For $k\in \N$ and $0<\alpha \leq 1$, let
$\DD_{k,\alpha}$ be the set of functions $f : {\CC^{(1)}} \to \R$ such that the
norm
  \begin{equation}
  \norm{f}_{\DD_{k,\alpha}}:= \sup_{x\in {\CC^{(1)}}} |f(x)| \sys(x)^k + \sup_{x\in
  {\CC^{(1)}}} \omega_\alpha(f,x) \sys(x)^k
  \end{equation}
is finite. This is the set of functions which are locally
$\alpha$-H\"{o}lder at each point and do not behave worse than
$\sys(x)^{-k}$ at infinity. When $f$ is compactly supported, this
condition reduces to the fact that $f$ is $\alpha$-H\"{o}lder, but it is
much more permissive in general.

For example, if a function $f:{\CC^{(1)}} \to \R$ is compactly supported and
$C^1$ (meaning that its lift to the manifold $\QQ^{(1)}(M,\Sigma,
\kappa)$ is $C^1$), then it belongs to all spaces $\DD_{k,\alpha}$.

The main result of this article is the following theorem:

\begin{thm}
\label{decay_teichmuller} Let $k\in \N$ and $0<\alpha \leq 1$. Let
$p,q \in \R_+$ be such that $1/p + 1/q <1$. Then there exist
constants $\delta>0$ and $C>0$ (depending on $k,\alpha,p,q$) such
that, for all functions $f: {\CC^{(1)}} \to \R$ belonging to
$\DD_{k,\alpha} \cap L^p(\nu_{\CC^{(1)}})$ and $g:{\CC^{(1)}} \to \R$ belonging to
$\DD_{k,\alpha} \cap L^q(\nu_{\CC^{(1)}})$, for all $t\geq 0$, holds
  \begin{equation*}
  \left| \int f \cdot g\circ \TF_t \dd\nu_{\CC^{(1)}} - \left( \int f
  \dd\nu_{\CC^{(1)}} \right) \left( \int g \dd\nu_{\CC^{(1)}} \right) \right|
  \leq C \bigl( \norm{f}_{\DD_{k,\alpha}}+\norm{f}_{L^p} \bigr)
  \bigl( \norm{g}_{\DD_{k,\alpha}} + \norm{g}_{L^q} \bigr)
  e^{-\delta t}.
  \end{equation*}
\end{thm}

An important ingredient in the course of the proof will be
recurrence estimates to a given compact set. We give here a consequence
of these estimates, which is of independent interest:

\begin{thm}
\label{compact_return} Let $\delta>0$. Then there exist a compact
set $K\subset {\CC^{(1)}}$ and a constant $C>0$ such that, for all $t\geq
0$,
  \begin{equation}
  \nu_{\CC^{(1)}}\{ x \in {\CC^{(1)}} \tq \forall s\in [0,t],\;  \TF_s (x) \not\in
  K\} \leq C e^{-(1-\delta) t}.
  \end{equation}
\end{thm}

This result easily implies the following corollary:
\begin{cor}
\label{cor:estim_systole} For all $\delta>0$, there exists $C>0$
such that, $\forall \epsilon\geq 0$,
  \begin{equation}
  \label{eq:estim_systole}
  \nu_{\CC^{(1)}}\{x\in {\CC^{(1)}} \tq \sys(x)<\epsilon\} \leq C \epsilon^{2 -
  \delta}.
  \end{equation}
\end{cor}
\begin{proof}
Let $K$ be a compact subset as in Theorem \ref{compact_return}. On
$K$, the systole is larger than a constant $\epsilon_0$. If
$\sys(x)<\epsilon<\epsilon_0$, then $\TF_t x\not \in K$ for $|t|
\leq \log ( \epsilon_0/ \epsilon)$ since $\log \circ \sys$ is
$1$-Lipschitz and $d(\TF_t x,x)\leq |t|$. Hence,
  \begin{align*}
  \nu_{\CC^{(1)}}\bigl\{x\in {\CC^{(1)}} \tq \sys(x)<\epsilon\bigr\}&
  \leq \nu_{\CC^{(1)}}\bigl\{ x\in {\CC^{(1)}} \tq \forall s\in [-\log
  (\epsilon_0/ \epsilon), \log (\epsilon_0/
  \epsilon)], \TF_s(x)
  \not\in K\bigr\}
  \\&
  \!\!\!\!\!\!\!\!\!\!\!\!\!
  = \nu_{\CC^{(1)}}\bigl\{ x\in {\CC^{(1)}} \tq \forall  s \in [0, 2\log
  (\epsilon_0/ \epsilon)], \TF_s(x)
  \not\in K\bigr\}
  \leq C \left( \frac{\epsilon}{\epsilon_0} \right)^{2(1-\delta)}.
  \qedhere
  \end{align*}
\end{proof}

This estimate is known not to be optimal: by the Siegel-Veech
formula (see e.g.\ \cite{EM}), there exists a constant $C>0$ such
that
  \begin{equation}
  \label{eq:siegelveech}
  \nu_{\CC^{(1)}}\{x\in {\CC^{(1)}} \tq \sys(x)<\epsilon\} \sim C\epsilon^{2}.
  \end{equation}
Notice however that the proof of this result relies heavily on the
$\SL(2,\R)$ action, while our estimates involve only the Teichm\"{u}ller
flow. Since the loss between \eqref{eq:siegelveech} and
\eqref{eq:estim_systole} is arbitrarily small,
Theorem \ref{compact_return} is quite sharp. In particular, the
combinatorial estimates we will develop in Section
\ref{sec:distortion_estimate} for the proofs of Theorems
\ref{decay_teichmuller} and \ref{compact_return} are quasi-optimal.

\begin{rem}
As a consequence of Corollary \ref{cor:estim_systole} (or of
Equation \eqref{eq:siegelveech}), the function $\phi : x \mapsto
1/\sys(x)$ belongs to $L^p$ for all $p<2$. Moreover, Lemma
\ref{systole_Lpischitz} shows that $\phi \in \DD_{1,1}$. 
\end{rem}

\section{The Veech flow}

In this section we introduce the Veech flow, and discuss its basic
combinatorics, related to interval exchange transformations.  The
Veech flow is a finite cover of the Teichm\"{u}ller flow,
and it will be shown in the next section that our results for the
Teichm\"{u}ller flow follow from corresponding results for the Veech flow.

We follow the presentation of \cite {MMY}.

\subsection{Rauzy classes and interval exchange transformations}

\subsubsection{Interval exchange transformations}

An interval exchange transformation is defined as follows.
Let $\AA$ be some fixed alphabet on $d\ge 2$ letters.
\begin{enumerate}
\item Take an interval $I \subset \R$ (all
intervals will be assumed to be closed at the left and open at the right),
\item Break it into $d \geq 2$ intervals $\{I_\alpha\}_{\alpha \in \AA}$,
\item Rearrange the intervals in a new order (via translations) inside $I$.
\end{enumerate}

Modulo translations, we may always assume that the left endpoint of $I$ is
$0$.  Thus the interval exchange transformation is entirely defined by the
following data:
\begin{enumerate}
\item The lengths of the intervals $\{I_\alpha\}_{\alpha \in \AA}$,
\item Their orders before and after rearranging.
\end{enumerate}
The first are called length data, and are given by a vector $\lambda \in
\R^\AA_+$ (here and henceforth $\R_+=(0,\infty)$).  The second are called combinatorial data, and are given by a pair
of bijections $\pi=(\pi_t,\pi_b)$ from $\AA$ to $\{1,\ldots,d\}$ (we will
sometimes call such a pair of bijections a permutation).
We denote the set of all such pairs of bijections by $\sssigma(\AA)$.
The bijections $\pi_\varepsilon:\AA \to \{1,\ldots,d\}$ can be viewed
as rows where the elements of $\AA$ are displayed in the order
$(\pi_\varepsilon^{-1}(1),...,\pi_\varepsilon^{-1}(d))$.
Thus we can see an element of
$\sssigma(\AA)$ as a pair of rows, the top (corresponding to $\pi_t$) and
the bottom (corresponding to
$\pi_b$) of $\pi$.  The interval exchange transformation associated to these
data will be denoted $f=f(\lambda,\pi)$.

Notice that if the combinatorial data are such that the set of the first $k$
elements in the top and bottom of $\pi$ coincide for some $1 \leq k<d$ then,
irrespective of the
length data, the interval exchange transformation splits into two simpler
transformations.
We are mostly interested in combinatorial data for which this does not
happen,
which we will call ${\it irreducible}$.  Let $\ssigma(\AA) \subset
\sssigma(\AA)$
be the set of irreducible combinatorial data.

\subsubsection{Rauzy classes} \label {rauzy c}

A {\it diagram} (or directed graph) consists of two kinds of objects,
vertices and (oriented) arrows joining two vertices. Thus, an arrow has a
start and an end.
A {\it path} of length $m \geq 0$
in the diagram is a finite sequence $v_0, \ldots,
v_m$ of vertices and a sequence of arrows
$a_1, \ldots, a_m$ such that $a_i$ starts at
$v_{i-1}$ and ends in $v_i$.  A path is said to start at $v_0$, end in
$v_m$, and pass through $v_1,...,v_{m-1}$.
If $\gamma_1$ and $\gamma_2$ are paths such
that the end of $\gamma_1$ is the start of $\gamma_2$, their concatenation
is also a path, denoted by $\gamma_1 \gamma_2$.
We can identify paths of length zero with vertices and paths of length one
with arrows.  Paths of length zero are called trivial.  We introduce a
partial order on paths: $\gamma_s \leq \gamma$ if and only if
$\gamma$ starts by $\gamma_s$.

Given $\pi \in \ssigma(\AA)$ we consider two operations.  Let $\alpha$ and
$\beta$ be the last elements of the top and bottom rows.
The {\it top} operation keeps the top row unchanged, and it changes the
bottom row by moving $\beta$ to the position immediately to
the right of the position occupied by $\alpha$.
When applying this operation to $\pi$, we will say that $\alpha$
{\it wins} and $\beta$ {\it loses}.
The  {\it bottom} operation is defined in a similar way, just interchanging
the words top and bottom, and the roles of $\alpha$ and $\beta$.  In this
case we say that $\beta$ wins and $\alpha$ loses. Notice that both
operations preserve the first elements of both the top and the bottom row.

It is easy to see that each of these operations gives a bijection of
$\ssigma(\AA)$.
A Rauzy class $\RR$ is a minimal non-empty subset of $\ssigma(\AA)$ which is
invariant under the top and bottom operations.  Given a Rauzy class $\RR$,
we define a diagram, called {\it Rauzy diagram}.  Its vertices are the
elements of $\RR$ and for each vertex
$\pi \in \RR$ and each of the operations considered above, we define an
arrow joining $\pi$ to the image of $\pi$ by the corresponding operation.
Notice that every vertex is the start and end of two arrows, one top and one
bottom.  Every arrow has a start, an end, a type (top or bottom),
a winner and a loser.
The set of all paths is denoted by $\Pi(\RR)$.

\subsubsection{Linear action}

Let $\RR \subset \ssigma(\AA)$ be a Rauzy class.  To each path $\gamma \in
\Pi(\RR)$, we associate a linear map $B_\gamma \in \SL(\AA,\Z)$ as follows.
If $\gamma$ is trivial, then $B_\gamma=\id$.  If $\gamma$ is an arrow with
winner $\alpha$ and loser $\beta$ then $B_\gamma \cdot e_\xi=e_\xi$ for
$\xi \in \AA \setminus \{\alpha\}$ and
$B_\gamma \cdot e_\alpha=e_\alpha+e_\beta$,
where $\{e_\xi\}_{\xi \in \AA}$ is the canonical basis of $\R^\AA$.  We extend
the definition to paths so that $B_{\gamma_1 \gamma_2}=B_{\gamma_2} \cdot
B_{\gamma_1}$.

\subsection{Rauzy induction}

Let $\RR \subset \ssigma(\AA)$ be a Rauzy class, and define
$\Delta^0_\RR=\R^\AA_+ \times \RR$.
Given $(\lambda,\pi)$ in $\Delta^0_\RR$, we say that we can apply Rauzy
induction to $(\lambda,\pi)$ if $\lambda_\alpha \neq \lambda_\beta$, where
$\alpha, \beta \in \AA$ are the last elements of the top and bottom rows of
$\pi$, respectively. Then we define $(\lambda',\pi')$ as follows:
\begin{enumerate}
\item Let $\gamma=\gamma(\lambda,\pi)$ be a top or bottom arrow on the Rauzy
diagram starting at $\pi$, according to whether
$\lambda_\alpha>\lambda_\beta$ or $\lambda_\beta>\lambda_\alpha$.
\item Let $\lambda'_\xi=\lambda_\xi$ if $\xi$ is not the winner of
$\gamma$, and
$\lambda_\xi=|\lambda_\alpha-\lambda_\beta|$ if $\xi$
is the winner of $\gamma$.
\item Let $\pi'$ be the end of $\gamma$.
\end{enumerate}
We say that $(\lambda',\pi')$ is obtained from $(\lambda,\pi)$ by
applying Rauzy induction, of type top or bottom depending on
whether the type of $\gamma$ is top or bottom. We have that $\pi'
\in \RR$ and $\lambda' \in \R^\AA_+$. The interval
exchange transformations $f:I \to I$ and $f':I' \to I'$ specified
by the data $(\lambda,\pi)$ and $(\lambda',\pi')$ are related as
follows.  The map $f'$ is the first return map of $f$ to a
subinterval of $I$, obtained by cutting from $I$ a subinterval
with the same right endpoint and of length $\lambda_\xi$, where $\xi$
is the loser of $\gamma$.  The map $Q:(\lambda,\pi) \mapsto
(\lambda',\pi')$ is called {\it Rauzy induction map}.  Its domain
of definition, the set of all $(\lambda,\pi) \in \Delta^0_\RR$
such that $\lambda_\alpha \neq \lambda_\beta$ (where $\alpha$ and $\beta$
are the last letters in the top and bottom rows of $\pi$),
will be denoted by $\Delta^1_\RR$.

The connected components $\Delta_\pi=\R^\AA_+ \times \{\pi\}$ of
$\Delta^0_\RR$ are naturally labeled by the elements of $\RR$, or
equivalently, by paths in $\Pi(\RR)$ of length $0$. The connected
components $\Delta_\gamma$ of $\Delta^1_\RR$ are naturally labeled
by arrows, that is, paths in $\Pi(\RR)$ of length $1$. One easily
checks that each connected component of $\Delta^1_\RR$ is mapped
homeomorphically to some connected component of $\Delta^0_\RR$.

Let $\Delta^n_\RR$ be the domain of $Q^n$, $n \geq 2$. The connected
components of $\Delta^n_\RR$ are naturally labeled by paths in
$\Pi(\RR)$ of length $n$: if $\gamma$ is obtained by following a
sequence of arrows $\gamma_1,...,\gamma_n$, then $\Delta_\gamma=\{x
\in \Delta^0_\RR \tq Q^{k-1}(x) \in \Delta_{\gamma_k},\,
1 \leq k \leq n\}$.  Notice that if $\gamma$ starts at $\pi$ then
$\Delta_\gamma=(B^*_\gamma \cdot \R^\AA_+) \times \{\pi\}$ (here and in the
following we will use $A^*$ to denote the transpose of a matrix $A$).
Indeed for arrows this follows from the definitions, and the extension to
paths is then immediate.

If $\gamma$ is a path in $\Pi(\RR)$
of length $n$ ending at $\pi \in \RR$, let
\be
Q^\gamma=Q^n:\Delta_\gamma \to \Delta_\pi.
\ee
This map is a homeomorphism.

Let $\Delta^\infty_\RR=\bigcap_{n \geq 0}
\Delta^n_\RR$.  A sufficient condition for $(\lambda,\pi)$ to belong to
$\Delta^\infty_\RR$ is for the coordinates of $\lambda$ to be independent
over $\Q$.

\subsubsection{Complete and positive paths}

\begin{definition} \label {as defined}

Let $\RR \subset \ssigma(\AA)$ be a Rauzy class.
A path $\gamma \in \Pi(\RR)$ is called {\it complete} if every
$\alpha \in \AA$ is the winner of some arrow composing $\gamma$.

\end{definition}

\begin{lemma}[\cite {MMY}, \S 1.2.3, Proposition]
\label{Lemme_MMY}
Let $(\lambda,\pi) \in \Delta^\infty_\RR$, and let
$\Delta_{\gamma(n)}$ be the connected component of $(\lambda,\pi)$ in
$\Delta^n_\RR$.  Then $\gamma(n)$ is complete for all $n$ large enough.

\end{lemma}

In particular any Rauzy diagram contains complete paths.

We say that $\gamma \in \Pi(\RR)$ is $k$-complete if it is a concatenation
of $k$ complete paths.  We say that $\gamma \in \Pi(\RR)$ is {\it positive}
if $B_\gamma$ is given, in the canonical basis of $\R^\AA_+$,
by a matrix with all entries positive.

\begin{lemma}[\cite {MMY}, \S 1.2.4, Lemma] \label {complete 2d-3}

If $\gamma$ is a $k$-complete path with $k \geq 2 \# \AA-3$, then $\gamma$
is positive.

\end{lemma}

\subsection{Zippered rectangles} \label {recall}

Let $\RR \subset \ssigma(\AA)$ be a Rauzy class.
Let $\pi=(\pi_t,\pi_b) \in \RR$.
Let $\Ga_\pi \subset \R^\AA$ be the set of all
$\tau$ such that
\be
\sum_{\pi_t(\xi) \leq k} \tau_\xi>0 \quad\text{and}\quad
\sum_{\pi_b(\xi) \leq k} \tau_\xi<0 \quad\text{for all } 1 \leq k \leq
d-1.
\ee
Notice that $\Ga_\pi$ is an open convex polyhedral cone.  It is non-empty,
since the vector $\tau$ with coordinates $\tau_\xi=\pi_b(\xi)-\pi_t(\xi)$
belongs to $\Ga_\pi$.

From the data $(\lambda,\pi,\tau)$, it is possible to define a marked
translation surface $S=S(\lambda,\pi,\tau)$ in some
$\QQ_{g,\kappa}$, where $g$ and $\kappa$ depend only on $\pi$ (see \cite
{MMY}, \S 3.2).  It is obtained (in the {\it zippered rectangles}
construction) by gluing rectangles of horizontal sides
$\lambda_\alpha$ and vertical sides $h_\alpha$, where the height vector $h
\in \R^\AA_+$ is given by $h=-\Omega(\pi) \cdot \tau$, and $\Omega(\pi)$ is
the linear operator on $\R^\AA$,
\be\label{defOmega}
\langle \Omega(\pi) \cdot e_x,e_y \rangle=
\left \{ \begin{array}{ll}
1, & \pi_t(x)>\pi_t(y), \pi_b(x)<\pi_b(y),\\[5pt]
-1, & \pi_t(x)<\pi_t(y), \pi_b(x)>\pi_b(y),\\[5pt]
0, & \text {otherwise}.
\end{array}
\right .
\ee
In particular, the area of the translation surface $S$ is
$A(\lambda,\pi,\tau)=-\langle \lambda,\Omega \cdot \tau \rangle$.

\subsubsection{Extension of induction to the space of zippered rectangles}

If $\gamma \in \Pi(\RR)$ is a path starting at $\pi$, let $\Ga_\gamma
\subset \R^\AA$ be defined by the condition
\be
B_\gamma^* \cdot \Ga_\gamma=\Ga_\pi.
\ee
If $\gamma$ is
a top arrow ending at $\pi'$, then $\Ga_\gamma$ is the set of
all $\tau \in \Ga_{\pi'}$ such that $\sum_{x \in \AA}
\tau_x<0$, and if $\gamma$ is a bottom arrow ending at $\pi'$, then
$\Ga_\gamma$ is the set of all $\tau \in \Ga_{\pi'}$ such
that $\sum_{x \in \AA} \tau_x>0$.  Thus, the map
\be
\widehat Q^\gamma:\Delta_\gamma \times \Ga_\pi \to
\Delta_{\pi'} \times \Ga_\gamma,
\quad \widehat Q^\gamma(\lambda,\pi,\tau)=
(Q(\lambda,\pi),(B_\gamma^*)^{-1} \cdot \tau)
\ee
is invertible.
Now we can define an invertible map
by putting together the $\widehat Q^\gamma$ for every arrow
$\gamma$.  This is a map from $\bigcup \Delta_\gamma \times \Ga_\pi$
(where the union is taken over all $\pi \in \RR$ and all arrows $\gamma$
starting at $\pi$) to $\bigcup \Delta_{\pi'} \times \Ga_\gamma$ (where the
union is taken over all $\pi' \in \RR$ and all arrows ending at $\pi'$).
We let $\widehat \Delta_\RR=\bigcup_{\pi \in \RR} \Delta_\pi \times \Ga_\pi$.
The map $\widehat Q$ is a skew-product over $Q$: $\widehat
Q(\lambda,\pi,\tau)=(Q(\lambda,\pi),\tau')$ where $\tau'$ depends on
$(\lambda,\pi,\tau)$.

The translation surfaces $S$ and $S'$ corresponding to $(\lambda,\pi,\tau)$
and $Q(\lambda,\pi,\tau)$ are obtained by appropriate cutting and pasting,
so they correspond to the same element in the moduli space $\MM_{g,\kappa}$
(the marking on the homology is however not preserved), see \cite {MMY}, \S
4.1.  We have thus a well defined map $\proj:\widehat \Delta_\RR \to \CC$
satisfying
\be \label {commutat}
\proj \circ \widehat Q=\proj,
\ee
where $\CC=\CC(\RR)$ is a connected
component of $\MM_{g,\kappa}$ (the connectivity of the image of $\proj$ is
due to the relation (\ref {commutat})).  In particular
$g$ and $\kappa$ only depend on $\RR$.

\begin{thm}[Veech]

If $\CC$ is a connected component of $\MM_{g,\kappa}$ then there exists a
Rauzy class $\RR$ such that $\CC=\CC(\RR)$.

\end{thm}

\begin{thm}[Veech]

The image of $\proj:\widehat \Delta_\RR \to \CC$ has full Lebesgue measure
in $\CC$.

\end{thm}

The action of $\widehat Q$ on $\widehat \Delta_\RR$ admits a nice
fundamental domain.  Let $\phi(\lambda,\pi,\tau)=\|\lambda\|=
\sum_{\alpha \in \AA}
\lambda_\alpha$.
Let $\Om_\RR \subset \widehat \Delta_\RR$ be the set of all $x$ such
that either
\begin{enumerate}
\item $\widehat Q(x)$ is defined and $\phi(\widehat Q(x))<1 \leq \phi(x)$,
\item $\widehat Q(x)$ is not defined and $\phi(x) \geq 1$,
\item $\widehat Q^{-1}(x)$ is not defined and $\phi(x)<1$.
\end{enumerate}
It is a fundamental domain for the action of $\widehat Q$: each orbit
of $\widehat Q$ intersects $\Om_\RR$ in exactly one point.
The fibers of the map $\proj:\Om_\RR \to \CC$ are almost everywhere
finite (with constant cardinality).
The projection of the standard Lebesgue measure on
$\Om_\RR$ is (up to scaling) the standard volume form on $\CC$.

\subsubsection{The Veech flow}

There is a natural flow $\TV_t:\widehat \Delta_\RR \to \widehat \Delta_\RR$,
$\TV_t(\lambda,\pi,\tau)=(e^t \lambda,\pi,e^{-t} \tau)$, which lifts the
Teichm\"{u}ller flow in $\MM_{g,\kappa}$.  This flow commutes
with $\widehat Q$.  The {\it Veech flow} $\VT_t:\Om_\RR \to \Om_\RR$ is
defined by $\VT_t(x)=\widehat Q^n(\TV_t(x))$ where $n$ is
the unique value such that $\widehat Q^n(\TV_t(x)) \in \Om_\RR$.  It lifts the
Teichm\"{u}ller flow on $\CC$:
\be
\proj \circ \VT_t=\TF_t \circ \proj.
\ee
Since both the flow $\TV_t$ and the map $\widehat Q$ trivially preserve the
standard Lebesgue measure on $\widehat \Delta_\RR$, the Veech flow $\VT_t$
preserves the standard Lebesgue measure on $\Om_\RR$.

Let $\Om^{(1)}_\RR=\proj^{-1}(\CC^{(1)})$
be the set of all $(\lambda,\pi,\tau)$ such that $A(\lambda,\pi,\tau)=1$.
The Veech flow leaves invariant $\Om^{(1)}_\RR$.  It follows that its
restriction $\VT_t:\Om^{(1)}_\RR \to \Om^{(1)}_\RR$ leaves invariant a
smooth volume form $\dd\omega$ 
(such that $\dd\omega \wedge \dd A=\dLeb$), whose
projection is, up to scaling, the standard volume form on $\CC^{(1)}$.

\begin{rem}

Veech's proof of the fact that
the standard volume form on $\CC^{(1)}$ is finite actually first establishes
finiteness of the lift measure on $\Om^{(1)}_\RR$.  A different proof of
finiteness follows from our recurrence estimates.

\end{rem}

\begin{rem}

Finiteness is a crucial
step in Veech's proof of
conservativity of an absolutely continuous invariant measure
for the Rauzy renormalization (which is itself the center of Veech's proof
of unique ergodicity for typical interval exchange transformations \cite
{V1}).  A different proof of conservativity for the Rauzy
renormalization follows immediately from our recurrence estimates
(the proof of which does not depend on the zippered rectangle construction).

\end{rem}

\section{Reduction to recurrence estimates}

\subsection{Measurable models}

\subsubsection{The Veech flow as suspension over the Rauzy renormalization}

Let $\widehat \Up_\RR \subset \Om_\RR$ be the set of all
$(\lambda,\pi,\tau)$ with $\phi(\lambda,\pi,\tau)=\|\lambda\|=
\sum_{\alpha \in \AA} \lambda_\alpha=1$.  The connected components of
$\widehat \Up_\RR$ are naturally denoted $\widehat \Up_\pi$.  Let
$\widehat \Up^{(1)}_\RR=\Om^{(1)}_\RR \cap \widehat \Up_\RR$,
$\widehat \Up^{(1)}_\pi=\Om^{(1)}_\RR \cap \widehat \Up_\pi$.  Let
$\Up^n_\RR \subset \Delta^n_\RR$ be the set of
$(\lambda,\pi)$ with $\|\lambda\|=1$.  We let $\widehat m$ denote
the induced Lebesgue measure to $\widehat \Up^{(1)}_\RR$.

Notice that $\widehat \Up^{(1)}_\RR$ is transverse to the Veech flow on
$\Om^{(1)}_\RR$.
We are interested in the first return map $\widehat R$ to
$\widehat \Up^{(1)}_\RR$.  Its domain is the intersection of $\widehat
\Up^{(1)}_\RR$ with the domain of definition of $\widehat Q$, and we have
\be
\widehat R(\lambda,\pi,\tau)=(e^r \lambda',\pi',e^{-r} \tau'),
\ee
where $(\lambda',\pi',\tau')=\widehat Q(\lambda,\pi,\tau)$ and
$r=r(\lambda,\pi)=-\ln \|\lambda'\|=-\ln \phi(\lambda,\pi,\tau)$
is the first return time.  The map $\widehat R$ is a skew-product: $\widehat
R(\lambda,\pi,\tau)=(R(\lambda,\pi),e^{-r} \tau')$.  The map
$R:\Up^1_\RR \to \Up^0_\RR$ is called
the {\it Rauzy renormalization map}.  The measure $\widehat m$ is invariant
under $\widehat R$.

The Veech flow can thus be seen as a special suspension over the map
$\widehat R$, which is itself an ``invertible extension'' of a
non-invertible map $R$.  This ``suspension model'' loses control of some
orbits (the ones that do not return to $\widehat \Up^{(1)}_\RR$),
but those have zero
Lebesgue measure, and will not affect further considerations.

\subsubsection{Precompact sections}

In the above suspension model for the Veech flow, the underlying discrete
transformation $\widehat R$ is only very weakly hyperbolic.  This is related
to the fact that the section $\widehat
\Up^{(1)}_\RR$ is too large (for instance, it
has infinite area).  Zorich \cite {Z}
has introduced an alternative section with
finite area, but such a section is still somewhat too large, so that there
is not a good control on distortion.  In the following we will introduce
a class of suitably small (precompact in $\widehat \Up^{(1)}_\RR$)
sections with good distortion estimates.

The section we will choose will be the intersection of $\widehat
\Up^{(1)}_\RR$ with
(finite unions of) sets of the form $\Delta_\gamma \times \Ga_{\gamma'}$. 
Precompactness in the $\lambda$ direction is equivalent to having
$B^*_\gamma \cdot (\overline \R^\AA_+ \setminus \{0\}) \subset
\R^\AA_+$, which is equivalent to $\gamma$ being a positive path.  To take
care of both the $\lambda$ and
the $\tau$ direction, we introduce the following notion.

\begin{definition}

A path $\gamma$, starting in $\pi_s$ and ending in $\pi_e$,
is said to be strongly positive if it is positive and
$(B^*_\gamma)^{-1} \cdot (\overline \Ga_{\pi_s} \setminus \{0\}) \subset
\Ga_{\pi_e}$.

\end{definition}

\begin{rem}

According to Bufetov (personal communication), a positive path is
automatically strongly positive, but we will not need this fact.

\end{rem}

\begin{lemma} \label {strongly positive}

Let $\gamma$ be a $k$-complete path with $k \geq 3 \#\AA-4$.
Then $\gamma$ is strongly positive.

\end{lemma}

\begin{proof}

Let $d=\#\AA$.  Fix $\tau \in \overline \Ga_{\pi_s} \setminus \{0\}$.
Write $\gamma$ as a concatenation of arrows $\gamma=\gamma_1...\gamma_n$,
and let $\pi^{i-1}$ and $\pi^i$ denote the start and the end of $\gamma_i$. 
Let $\tau^0=\tau$,
$\tau^i=(B^*_{\gamma_i})^{-1} \cdot \tau^{i-1}$.  We must show that $\tau^n
\in \Ga_{\pi^n}$.

Let $h^i=-\Omega(\pi^i) \cdot \tau^i$.  Notice that $\tau \in \overline
\Ga_{\pi^0} \setminus \{0\}$ implies that $h^0 \in \overline \R^\AA_+
\setminus \{0\}$.  Indeed, since $\tau \in \overline \Ga_{\pi^0}$,
for every $\xi \in \AA$, we have 
$\sum_{\pi^0_t(\alpha)<\pi^0_t(\xi)} \tau_\alpha \geq 0$,   
$\sum_{\pi^0_b(\alpha)<\pi^0_b(\xi)} \tau_\alpha \leq 0$.  
Moreover, since $\tau
\neq 0$, there exists $1 \leq k^t,k^b \leq d$ minimal such that
$\tau_{(\pi^0_t)^{-1}(k^t)} \neq 0$ and $\tau_{(\pi^0_b)^{-1}(k^b)} \neq 0$. 
Since $\pi^0$ is irreducible, $\min\{k^t,k^b\}<d$.
Noticing that
\be
h^0_\xi=\sum_{\pi^0_t(\alpha)<\pi^0_t(\xi)} \tau_\alpha-
\sum_{\pi^0_b(\alpha)<\pi^0_b(\xi)} \tau_\alpha,
\ee
we see that $h^0_\xi \geq 0$ for all $\xi$, and the inequality is strict if
$\pi^0_t(\xi)=k^t+1$ (if $k^t<d$) or if $\pi^0_b(\xi)=k^b+1$ (if $k^b<d$).

Notice that $h^i=B_{\gamma_i} \cdot h^{i-1}$, so
if $\gamma_1 ... \gamma_i$ is a positive path then $h^i \in \R^\AA_+$.

Let $0 \leq k^t_i,k^b_i \leq d-1$ be maximal such that
\be
\sum_{\pi^i_t(\xi) \leq k} \tau^i_\xi>0 \quad\text{for all }
1 \leq k \leq k^t_i,
\ee
\be
\sum_{\pi^i_b(\xi) \leq k} \tau^i_\xi<0 \quad\text{for all } 1 \leq k
\leq k^b_i,
\ee
where $\pi^i_t$ and $\pi^i_b$ are the top and the bottom of $\pi^i$.
We claim that
\begin{enumerate}
\item If $h^{i-1} \in \R^\AA_+$ then $k^t_i \geq k^t_{i-1}$ and $k^b_i \geq
k^b_{i-1}$,
\item If $h^{i-1} \in \R^\AA_+$ and the winner of $\gamma_i$ is one of the
first $k^t_{i-1}+1$ letters in the top of $\pi^{i-1}$ then $k^t_{i} \geq \min
\{d-1,k^t_{i-1}+1\}$,
\item If $h^{i-1} \in \R^\AA_+$ and the winner of $\gamma_i$ is one of the
first $k^b_{i-1}+1$ letters in the bottom of $\pi^{i-1}$ then $k^b_{i} \geq
\min \{d-1,k^b_{i-1}+1\}$.
\end{enumerate}

Let us see that (1), (2) and (3) imply the result, which is equivalent to
the statement that $k^t_n=d-1$ and $k^b_n=d-1$.  We will show that
$k^t_n=d-1$, the other estimate being analogous.  Let us write
$\gamma=\gamma_{(1)}...\gamma_{(3d-4)}$ where $\gamma_{(j)}$ is complete. 
Write $\gamma_{(j)}=\gamma_{s_j}...\gamma_{e_j}$.
By Lemma \ref {complete 2d-3}, $h^k \in \R^\AA_+$ for $k \geq
e_{2d-3}$.  From the definition of a complete path, for each $j>2d-3$, there
exists $e_{j-1}<i \leq e_j$ such that the winner of $\gamma_i$ is one of the
first $k^t_{e_{j-1}}+1$ letters in the top of $\pi^{i-1}$.  It follows that
$k^t_{e_j} \geq \min \{d-1,k^t_{e_{j-1}}+1\}$, and so $k^t_n=k^t_{e_{3d-4}}
\geq \min \{d-1,k^t_{e_{2d-3}}+d-1\}=d-1$.

We now check (1), (2) and (3).  Assume that $h^{i-1} \in  \R^\AA_+$, and
that $\gamma_i$ is a top, the other case being analogous.
In this case $\pi^i_t=\pi^i_{t-1}$ and $\tau^i_\alpha=\tau^{i-1}_\alpha$ for
$\pi^i_t(\alpha)<d$, hence $k_i^t \geq k_{i-1}^t$.  This shows that the
first claim of (1) holds.  Moreover, (2) also holds since its hypothesis can
only be satisfied if $k_{i-1}^t=d-1$.

If the winner
of $\gamma_i$ is not one of the $k^b_{i-1}+1$ first letters in the bottom of
$\pi^{i-1}$, then for every $\alpha \in \AA$ such that $1 \leq
\pi^{i-1}_b(\alpha) \leq k^b_{i-1}$,
we have $\pi^{i-1}_b(\alpha)=\pi^i_b(\alpha)$,
$\tau^{i-1}_\alpha=\tau^i_\alpha$, so $k^b_i \geq k^b_{i-1}$.  

If the winner
$\beta$ of $\gamma_i$ appears in the $k$-th position in the bottom of
$\pi^{i-1}$ with $1 \leq k \leq k^b_{i-1}+1$, then
\be
\sum_{\pi^i_b(\xi) \leq j} \tau^i_\xi=
\sum_{\pi^{i-1}_b(\xi) \leq j} \tau^{i-1}_\xi<0 \quad\text{for all }
1 \leq j \leq k-1,
\ee
\be
\sum_{\pi^i_b(\xi) \leq j} \tau^i_\xi=
\sum_{\pi^{i-1}_b(\xi) \leq j-1} \tau^{i-1}_\xi<0 \quad\text{for all }
k+1 \leq j \leq k^b_{i-1}+1,
\ee
\be
\sum_{\pi^i_b(\xi) \leq k} \tau^i_\xi=
\sum_{\pi^{i-1}_b(\xi) \leq d-1} \tau^{i-1}_\xi-h^{i-1}_\beta 
\leq -h^{i-1}_\beta<0,
\ee
which implies that $k^b_i \geq \min \{d-1,k^b_{i-1}+1\}$.

This shows that both (3) and the second claim of (1) must hold.
\end{proof}

\subsubsection{A better model}

We will now choose a specific precompact section, adapted for the problem of
exponential mixing (Theorem \ref {decay_teichmuller}).  Our particular
choice aims to simplify the combinatorial
description of the first return map.  We will later consider a different
choice for the recurrence problem (Theorem \ref {compact_return}).

Let $\gamma_* \in \Pi(\RR)$ be a strongly positive
path starting and ending in the same
$\pi \in \RR$.  Assume further
that if $\gamma_*=\gamma_s\gamma=\gamma \gamma_e$ then either
$\gamma=\gamma_*$ or $\gamma$
is trivial.\footnote {Notice that if $\gamma_*$ ends by
a bottom arrow and starts by a sufficiently long (at least half the length
of $\gamma_*$) sequence of top arrows then this
last condition is automatically satisfied.}
We will say that $\gamma_*$ is {\it neat}.

Let $\widehat \De=\widehat \Up^{(1)}_\RR \cap
(\Delta_{\gamma_*} \times \Ga_{\gamma_*})$, and let $\De=\Up^0_\RR \cap
\Delta_{\gamma_*}$.
We are interested in the first return map $T_{\widehat \De}$ to $\widehat
\De$ under the Veech flow.
The connected components of its domain are given by
$\widehat \Up^{(1)}_\RR \cap (\Delta_{\gamma\gamma^*} 
\times \Ga_{\gamma_*})$, where
$\gamma$ is either $\gamma_*$, or a minimal path of the form 
$\gamma_* \gamma_0 \gamma_*$ not beginning by $\gamma_*\gamma_*$.  
The restriction of $T_{\widehat \De}$ to such a component is
given by
\be
T_{\widehat \De}(\lambda,\pi,\tau)=\left (\frac
{(B_{\gamma}^*)^{-1} \cdot \lambda}
{\|(B_{\gamma}^*)^{-1} \cdot \lambda\|},\pi,
\|(B_{\gamma}^*)^{-1} \cdot \lambda\|
(B_{\gamma}^*)^{-1} \cdot \tau \right ).
\ee
The return time function is just
\be
r_{\widehat \De}(\lambda,\pi,\tau)=r_\De(\lambda,\pi)=
-\ln \|(B_{\gamma}^*)^{-1} \cdot \lambda\|.
\ee
The map $T_{\widehat \De}(\lambda,\pi,\tau)=(\lambda',\pi,\tau')$
is a skew-product over a
non-invertible transformation $T_\De(\lambda,\pi)=(\lambda',\pi)$.

The Veech flow can be seen as a suspension over $T_{\widehat \De}$,
with roof function $r_{\widehat \De}$.
In this suspension model, many more orbits escape
control (the ones that do not come back to $\widehat \De$).  Still, due to
ergodicity of the Veech flow, almost every orbit is captured by the
suspension model.

\subsection{Hyperbolic properties} \label {hyperbolic proper}

The transformation $T_{\widehat \De}$ turns out to have much better
hyperbolic properties than $\widehat R$.

\begin{lemma} \label {skew-product}

$T_{\widehat \De}$ is a hyperbolic skew-product over $T_\De$.

\end{lemma}

Implicit in the above statement is the choice of
probability measure $\nu$ and Finsler
metric $\| \cdot \|_{\widehat \De}$ which are part of the definition of a
hyperbolic skew-product.  The choice of $\nu$ is
clear (the normalized restriction of $\widehat m$ to $\widehat \De$) but
there is some freedom in the choice of the Finsler metric.  In order to
enforce the hyperbolicity properties we want from $T_{\widehat \De}$,
we will introduce a particular complete Finsler metric on
$\widehat \Up^{(1)}_\pi$, and then take $\| \cdot \|_{\widehat \De}$ as its
restriction.  By strong positivity of $\gamma_*$,
$\widehat \De$ is a precompact open subset of
$\widehat \Up^{(1)}_\RR$, so $\widehat \De$ will
have bounded diameter with respect to such metric.

\subsubsection{Hilbert metric}

The Hilbert pseudo-metric on $\R^2_+$ is $\dist_{\R^2_+}
(x,y)=\ln \max_{1 \leq i,j \leq 2} \frac {x_i y_j} {x_j y_i}$.
One easily checks that if $B \in \GL(2,\R)$ is a linear map such that
$B \cdot \R^2_+ \subset \R^2_+$ then $B$
contracts weakly the Hilbert pseudo-metric: $\dist_{\R^2_+}
(B \cdot x,B \cdot y) \leq \dist_{\R^2_+}(x,y)$.  In particular, the Hilbert
pseudo-metric is invariant under linear isomorphisms of $\R^2_+$.

More generally, if $C \subset \R^\AA \setminus \{0\}$ is an open
convex cone whose closure does not contain any one-dimensional subspace of
$\R^\AA$, one defines a Hilbert pseudo-metric on $C$ as follows.
If $x$ and $y$ are colinear then $\dist_C(x,y)=0$.  Otherwise, $C$
intersects the subspace generated by $x$ and $y$ in a cone
isomorphic to $\R^2_+$.  We let
$\dist_C(x,y)=\dist_{\R^2_+}(\psi(x),\psi(y))$ where $\psi$ is any
such isomorphism.  If $C=\R^\AA_+$ then we have
$\dist_C(x,y)=\max_{\alpha,\beta \in \AA} \ln \frac {x_\alpha y_\beta}
{x_\beta y_\alpha}$.

If $C' \subset C$ is a smaller cone then the inclusion $C' \to C$ is a weak
contraction of the respective Hilbert pseudo-metrics:
$\dist_C(x,y) \leq \dist_{C'}(x,y)$.  Moreover, if the diameter of $C'$ with
respect to $\dist_C$ is bounded by some $M$ then the contraction is
definite: $\dist_C(x,y) \leq \delta \dist_{C'}(x,y)$ where
$\delta=\delta(M)<1$.

We notice that the Hilbert pseudo-metric on a cone $C$ induces the Hilbert
metric on the space of rays $\{t x \tq t \in \R_+\}$ contained in $C$
(which is a projective manifold).  It is a complete Finsler metric.

\subsubsection{Uniform expansion and contraction}

Recall that $\widehat \Up^{(1)}_\pi$ is contained in
$\Delta_\pi \times \Ga_\pi$, which is a product of two cones.  In
$\Delta_\pi \times \Ga_\pi$, we have the product Hilbert pseudo-metric
$\dist((\lambda,\pi,\tau),
(\lambda',\pi,\tau'))=\dist_{\Delta_\pi}((\lambda,\pi),(\lambda',\pi))+
\dist_{\Ga_\pi}(\tau,\tau')$.
Each product of rays $\{(a \lambda,\pi,b \tau) \tq a,b \in \R_+\} \subset
\Delta_\pi \times \Ga_\pi$ intersects transversely $\Up^{(1)}_\RR$ in a
unique point.  It follows that the product Hilbert pseudo metric induces a
metric $\dist$ on $\widehat \Up^{(1)}_\pi$.  It is a complete Finsler
metric.

\smallskip

\noindent{\it {Proof of Lemma \ref {skew-product}}.}
Let us first show that $T_\De$ is a uniformly expanding Markov map (the
underlying Finsler metric being the restriction of
$\dist_{\Delta_\pi}$, and the underlying measure $\Leb$ being the
induced Lebesgue measure) .  It is clear that
$\De$ is a John domain.

Condition (1) of Definition \ref {markovmap}
is easily verified, except for the definite contraction of inverse branches. 
To check this property, we notice that an inverse branch can be written as
$h(\lambda,\pi)=\left (\frac {B_\gamma^* \cdot \lambda} {\|B_\gamma^*
\cdot \lambda\|},\pi \right )$.  Since $\gamma_*$ is neat, we can write
$B^*_\gamma=B^*_{\gamma_*} B^*_{\gamma_0}$ for some $\gamma_0$.  Thus $h$
can be written as (the restriction of) the composition of two maps
$\Delta_\pi \to \Delta_\pi$, $h=h_* \circ h_0$, where
$h_0$ is weakly contracting and $h_*$ is definitely
contracting by precompactness of $\De$ in $\Delta_\pi$
(which is a consequence of positivity of $\gamma_*$).

To check condition (2) of Definition \ref {markovmap},
let $h(\lambda,\pi)=\left (\frac {B_\gamma^* \cdot \lambda} {\|B_\gamma^*
\cdot \lambda\|},\pi \right )$ be an inverse branch of $T_\De$.
The Jacobian of $h$ at
$(\lambda,\pi)$ is $J \circ h(\lambda,\pi)=\left (\frac
{1} {\|B_\gamma^* \cdot \lambda\|} \right )^d$, where $d=\# \AA$. 
It follows that
\be
\frac {J \circ h(\lambda,\pi)} {J \circ h(\lambda',\pi)} \leq
\sup_{\alpha \in \AA} \left (\frac {\lambda_\alpha}
{\lambda'_\alpha} \right )^d \leq e^{d
\dist_{\Delta_\pi}((\lambda,\pi),(\lambda',\pi))},
\ee
so that $\ln J \circ h$ is $d$-Lipschitz with respect to
$\dist_{\Delta_\pi}$.

To see that $T_{\widehat \De}$ is a hyperbolic skew-product over $T_\De$,
one checks the conditions (1-4) of Definition \ref {skew-product map}.
Condition (1) is
obvious, and condition (4) follows from precompactness of $\widehat \De$ in
$\Delta_\pi \times \Ga_\pi$ as before.
Since $T_{\widehat \De}$ is a first return map,
the restriction of $\widehat m$
to $\widehat \De$ is $T_{\widehat \De}$-invariant.  Its normalization is the
probability measure $\nu$ of condition (2).  In order to check condition
(3), it is convenient to trivialize $\widehat \De$
to a product (via the natural diffeomorphism
$\widehat \De \to \De \times \P\Ga_{\gamma_*}$).  Since
$\nu$ has a smooth density with respect to the product of the Lebesgue
measure on the factors, condition (3) follows by the Leibniz rule.
\qed

\subsection{Basic properties of the roof function}

\comm{
\begin{lemma}
There exists $\epsilon>0$ such that $r_\De(\lambda,\pi)>\epsilon$ for any
$(\lambda,\pi)$ in the domain of $T_\De$.  Moreover, there exists no
$C^1$ function $\psi:\De \to \R$ such that $r_\De+\psi \circ T_\De-\psi$
is constant in each domain of definition of $T_\De$.

\end{lemma}

\begin{proof}
The first statement is obvious, since $r_\De(\lambda,\pi) \geq -\ln
|1-\min_{\alpha \in \AA} \lambda_\alpha|$ and $\min_{\alpha \in \AA}
\lambda_\alpha$ is bounded away from zero for $(\lambda,\pi)$ in the
domain of $T_\De$.
For the second statement,
}

Let $H(\pi)=\Omega(\pi) \cdot \R^\AA$.  Recall (from \ref {recall})
that if $\tau \in \Ga_\pi$
then $-\Omega(\pi) \cdot \tau \in \R^\AA_+$, and that $\Ga_\pi$ is
non-empty, so $H(\pi) \cap \R^\AA_+ \neq
\emptyset$.

\begin{lemma}

Let $\Gamma \subset \Pi(\RR)$ be the set of all $\gamma$ such that $\gamma
$ is either $\gamma_*$, or a minimal path of the 
form $\gamma_* \gamma_0 \gamma_*$ not beginning by $\gamma_*\gamma_*$.
Let $K \subset \P H(\pi)$ be a closed set such that $B_\gamma
\cdot K=K$ for every
$\gamma \in \Gamma$.  Then either $K=\emptyset$ or $K=\P H(\pi)$.

\end{lemma}

\begin{proof}

Let $\Pi(\pi) \subset \Pi(\RR)$ be
the set of all paths that start and end in $\pi$.
Then any element of $\gamma_* \Pi(\pi) \gamma_*$ is a concatenation of
elements of $\Gamma$.  It
follows that if $K$ is invariant under all $B_\gamma$, $\gamma \in \Gamma$,
then $K$ is invariant under all $B_\gamma$, $\gamma \in \Pi(\pi)$: indeed
$B_\gamma \cdot K=
B_{\gamma_*}^{-1} \cdot 
B_{\gamma_* \gamma \gamma_*} \cdot 
B_{\gamma_*}^{-1} \cdot K=K$, since $\gamma_*$ and
$\gamma_*\gamma\gamma_*$
are
concatenation of elements of $\Gamma$. 
According to Corollary 3.6 of \cite {AV}, this implies that $K$ is either
empty or equal to $\P H(\pi)$.
\end{proof}

\begin{lemma} \label {bla}

The roof function $r_\De$
is good (in the sense of Definition \ref {def:goodroof}).

\end{lemma}

\begin{proof}

We check conditions (1-3) of Definition \ref {def:goodroof}.
Let $\Gamma \subset \Pi(\RR)$ be the set defined in the previous lemma. 
Notice that $\Gamma$ consists of positive paths.

The set $\HH$ of inverse branches $h$ of $T_\De$ is in bijection
with $\Gamma$, since each inverse branch is of the form
$h(\lambda,\pi)=((B_{\gamma_h})^* \cdot \lambda,\pi)$ for some
$\gamma_h \in \Gamma$.

Let $h \in \HH$.
Then $r_\De(h(\lambda,\pi))=
\ln \|(B_{\gamma_h}^*) \cdot \lambda\|$.  Since $\gamma_h$ is positive,
$r_\De \geq \ln 2$, which implies condition (1).  Notice that
$r_\De \circ h=\frac {1} {d} \ln J \circ h$, where $J$ is as in the
condition (2) of Definition \ref {markovmap}, so (2) follows
(by the previous discussion, it even follows that $r_\De \circ h$ is
$1$-Lipschitz with respect to $\dist_{\Delta_\pi}$).

Let us check condition (3).  We identify the tangent space to $\De$ at a
point $(\lambda,\pi) \in \De$ with $V=\{\lambda \in \R^\AA \tq \sum
\lambda_\alpha=0\}$.  Assume that we can write $r_\De=
\psi+\phi \circ T_\De-\phi$ with $\phi$ $C^1$, $\psi$ locally constant. 
Write $r^{(n)}(x)=\sum_{j=0}^{n-1} r_\De(T_\De^j(\lambda,\pi))$.
Then $D(r^{(n)} \circ h^n)=D \phi-D(\phi \circ
h^n)$, which can be rewritten
\be
\frac {\|(B^*_{\gamma_h})^n \cdot v\|} {\|(B^*_{\gamma_h})^n \cdot
\lambda\|}=D\phi(\lambda,\pi) \cdot v-D(\phi \circ h^n)(\lambda,\pi)
\cdot v, \quad (\lambda,\pi) \in \De,\, v \in V,
\ee
or
\be
\frac {\langle v,B^n_{\gamma_h} \cdot (1,...,1) \rangle}
{\langle \lambda,B^n_{\gamma_h}
\cdot (1,...,1) \rangle}=
D\phi(\lambda,\pi) \cdot v-D(\phi \circ h^n)(\lambda,\pi) \cdot v, \quad
(\lambda,\pi) \in \De,\, v \in V.
\ee
Since $Dh^n \to 0$, we conclude that
$[B^n_{\gamma_h} \cdot (1,...,1)] \in \P \R^\AA$ converges to a limit $[w]
\in \P \R^\AA$ independent of $h$.  This obviously implies
that $[w]$ is invariant by all $B_{\gamma_h}$, $h \in \HH$.  Since $w$ is a
limit of positive vectors (vectors with positive coordinates), by the
Perron-Frobenius Theorem, $w$ is colinear with the (unique) positive
eigenvector of $B_{\gamma_h}$, which also corresponds to the largest
eigenvalue.  Recalling that $H(\pi)$ is invariant
under $B_{\gamma_h}$, and intersects $\R^\AA_+$, it follows that $w \in
H(\pi)$.  According to the previous lemma, $K=\{[w]\} \subset \P H(\pi)$
should be either empty or equal to the whole $\P H(\pi)$, so $H(\pi)$ should
be one-dimensional.  This gives a contradiction since $H(\pi)$ is even
dimensional (since $H(\pi)$ is the image of the antisymmetric operator
$\Omega(\pi)$).
\end{proof}

\subsection{A recurrence estimate and exponential mixing}

We will show later (in Section \ref{sec:proof_rec})
the following recurrence estimate.

\begin{thm} \label {exponentialtails}

The roof function $r_\De$ has exponential tails.

\end{thm}

We will now show how to conclude exponential mixing for the Teichm\"{u}ller
flow, Theorem \ref {decay_teichmuller},
assuming the above recurrence estimate and the abstract result on
exponential mixing for hyperbolic skew-product flows.

The map $T_{\widehat \De}$ and the roof function
$r_\De$ define together a flow $\widehat T_t$ on the space $\wde=
\{(x,y,s)\,:\, (x,y) \in \widehat \De,\,T_{\widehat\De}(x,y)\text{ is defined
and } 0 \leq s <r_\De(x)\}$.
Since $T_{\widehat \De}$ is
a hyperbolic skew-product (Lemma \ref {skew-product}),
and $r_\De$ is a good roof
function (Lemma \ref {bla}) with exponential tails (Theorem \ref
{exponentialtails}), $\widehat
T_t$ is an excellent hyperbolic semi-flow.  By Theorem \ref
{main_thm_hyperbolic}, we get
exponential decay of correlations
\be
C_t(\tilde f,\tilde g)=\int \tilde f \cdot \tilde g \circ \widehat T_t
\dd\nu-\int \tilde f \dd\nu \int
\tilde g \dd\nu,
\ee
for $C^1$ functions
$\tilde f$, $\tilde g$, that is
\be \label {tilde f}
|C_t(\tilde f,\tilde g)| \leq C e^{-3 \delta t} \|\tilde f\|_{C^1}
\|\tilde g\|_{C^1},
\ee
for some $C>0$, $\delta>0$. This estimate holds for $C^1$ functions on
$\widehat{\Delta}_r$, while Theorem \ref {decay_teichmuller} deals
with H\"older functions on $\CC^{(1)}$. Hence, one needs an additional
lifting and smoothing argument, provided by the following technical lemma.

Let $P:\widehat \Up^{(1)}_\RR \times \R \to \CC^{(1)}$
be given by $P(z,s)=\TF_s(\proj(z))$,
where $\proj:\widehat \Delta_\RR \to \CC$ is the natural projection.

\begin{lemma}
\label{lem:regularisation}
For every $k \in \N$, $0<\alpha \leq 1$, $p>p' \geq 1$, $\delta>0$,
there exist $C>0$, $\epsilon_0>0$ with the following property.
Let $f:\CC^{(1)} \to \R$ be a function belonging to $\DD_{k,\alpha} \cap
L^p(\nu_{\CC^{(1)}})$.  For every $t>0$, there exists a
$C^1$ function $f^{(t)}:\wde \to \R$,
such that $\|f \circ P-f^{(t)}\|_{L^{p'}(\nu)} \leq
C (\|f\|_{\DD_{k,\alpha}}+\|f\|_{L^p(\nu_{\CC^{(1)}})}) e^{-\epsilon_0 t}$ and
$\|f^{(t)}\|_{C^1(\wde)} \leq
C \|f\|_{\DD_{k,\alpha}} e^{\delta t}$.

\end{lemma}

\begin{proof}

We identify $\widehat \Up^{(1)}_\RR \cap \Delta_\pi \times \Ga_\pi$ with a
subset $U$ of $\R^{2d-2}$ via a map
$(\lambda,\pi,\tau) \mapsto (x,y)$, where $x,y \in \R^{d-1}$ are defined
by $x_i=\frac {\lambda_{\pi_t^{-1}(i+1)}} {\lambda_{\pi_t^{-1}(1)}}$,
$y_i=\frac {\tau_{\pi_t^{-1}(i+1)}} {\tau_{\pi_t^{-1}(1)}}$, $1 \leq i \leq 
d-1$ (here $\pi_t$ is the top of $\pi$).  In this way, $\widehat \De$
becomes a precompact subset of $\R^{2d-2}$.  Using this identification, we
will write $r_\De(x)$ for $r_\De(\lambda,\pi)$.

This also provides an identification of $\wde$ with a subset of
$U \times [0,\infty) \subset \R^{2d-1}$ via the map
$(\lambda,\pi,\tau,s) \mapsto (x,y,s)$.
We will use $\| \cdot \|$ to denote the usual norm in $\R^{2d-1}$, and
$\dist$ for the corresponding distance.

Let $\| \cdot \|_F$ be the
Finsler metric on $U \times \R$ obtained by pullback via $P$ of the
Finsler metric on $\CC^{(1)}$ defined in \S \ref {metricfinsler}.
At a point $(x,y,s) \in \widehat \De \times \R$, we have the
estimate $C^{-1} e^{-2 |s|} \|w\| \leq \|w\|_F \leq C
e^{2|s|} \|w\|$ where $w$ is a vector tangent to $(x,y,s)$.
This follows from precompactness of $\widehat \De$ when $s=0$, and the
general case follows from this one by applying the
Teichm\"uller flow, see \eqref {e2t}.  We let $\dist_F$ be the metric in
$\wde$ corresponding to $\| \cdot \|_F$.  We recall that $\wde$ is
disconnected, so the $\dist_F$ distance between two points of $\wde$
is sometimes infinite.

There is another Finsler metric $\| \cdot \|_\wde$ over $\wde$, which is
the product of $\| \cdot \|_{\widehat \De}$ (introduced in section \ref
{hyperbolic proper}) in the $(x,y)$ direction and the usual metric in the
$s$ direction.  We recall that it is with respect to this metric that the
$C^1(\wde)$ norm is defined.  One easily checks that
$C^{-1} \|w\| \leq \|w\|_\wde \leq C \|w\|$.

We may assume that $\|f\|_{\DD_{k,\alpha}} \leq 1$.
This implies that for $z_0=(x_0,y_0,s_0) \in \wde$,
$|f \circ P(z_0)| \leq C e^{ks_0}$ and if $\dist_F(z,z_0)=r<1$ then
$|f \circ P(z)-f \circ P(z_0)| \leq C e^{ks} r^\alpha$.

Let $\epsilon>0$.
Let $\phi^{(t)}:\R^{2d-1} \to [0,\infty)$ be
a $C^\infty$ function supported in
$\{z \in \R^{2d-1} \tq \|z\|<e^{-\epsilon t}/10\}$, such that
$\int_{\R^{2d-1}} \phi^{(t)}(z) \dd z=1$ and such that
$\|\phi^{(t)}\|_{C^1(\R^{2d-1})} \leq C e^{2 d \epsilon
t}$.  Let $\psi^{(t)}:\R^{2d-1} \to \R$ be given by
$\psi^{(t)}(x,y,s)=f \circ P(x,y,s)$ if $(x,y,s) \in \wde$ and $0 \leq s
\leq \epsilon t$, and $\psi^{(t)}(x,y,s)=0$ otherwise.
We will show that if $\epsilon$ is small enough then
one can take $f^{(t)}=\phi^{(t)} * \psi^{(t)}|\wde$, where $*$ denotes
convolution.

Let us first check the assertion $\|f^{(t)}\|_{C^1(\wde)} \leq C e^{\delta
t}$.  It is immediate to check that, by choosing $\epsilon>0$ small, we have
indeed $\|f^{(t)}\|_{C^1(\wde)} \leq C \|\psi^{(t)} *
\phi^{(t)}\|_{C^1(\R^{2d-1})} \leq C \|\psi^{(t)}\|_{L^1(\dd z)}
\|\phi^{(t)}\|_{C^1(\R^{2d-1})} \leq C e^{\delta t}$.

We will now check the other assertion $\|f \circ P-f^{(t)}\|_{L^{p'}(\nu)}
\leq C e^{-\epsilon_0 t}$, assuming
$\|f\|_{\DD_{k,\alpha}}+\|f\|_{L^p(\nu_{\CC^{(1)}})} \leq 1$.

Choose $C_0>\max\{4,2k/\alpha\}$.
Let $Y \subset \wde$ be the union of the connected components of $\wde$ which
intersect $\{(x,y,s) \in \R^{2d-1} \tq s>C_0^{-1} \epsilon t\}$.
Let $X \subset \wde \setminus Y$ be the set of points with
$\dist((x,y,s),\partial \wde) \leq 4 e^{-\epsilon t}$.  Thus $\wde \setminus
(X \cup Y)$ consists of points well inside the connected components of
$\wde$ with not so long (maximal) return time.

\begin{lemma} \label {nuy}

We have $\nu(X \cup Y) \leq C e^{-\epsilon t/C}$.

\end{lemma}

\begin{proof}

Since $r_\De$ is a good roof function, by condition (2) of
Definition \ref {def:goodroof} we
have $r_\De(x)>(C C_0)^{-1} \epsilon t$ for every $(x,y,s) \in Y$.
By Theorem \ref
{exponentialtails}, $\nu(Y) \leq C e^{-\epsilon t/C}$.

The boundary of each connected component of $\wde$ can be split in three
parts: a floor (containing points $(x,y,s)$ with $s=0$), a roof (containing
points $(x,y,s)$ such that $r_\De(x)=s$) and a remaining lateral part.

Points $(x,y,s) \in X$ are at distance at most $4 e^{-\epsilon
t}$ of either the floor, the roof, or the lateral part of the boundary of
their connected component in $\wde$: we can thus write $X=X_\floor \cup
X_\roof \cup X_\lat$ (there is non-trivial intersection of $X_\floor$ and
$X_\roof$ with $X_\lat$).  We will now show that each of those three sets
have $\nu$-measure at most $C e^{-\epsilon t/C}$.
Clearly $\nu(X_\floor) \leq C e^{-\epsilon t}$.

Using (\ref {e2t}) and condition (2) of Definition \ref {def:goodroof},
we see that if $(x,y)$ is in the domain of $T_{\widehat \De}$ then
$\|DT_{\widehat \De}(x,y)\|,\|DT_{\widehat \De}(x,y)^{-1}\|
\leq C e^{2 r_\De(x)}$.  Using condition (2) of
Definition \ref {def:goodroof} again we get $\|Dr_\De(x)\| \leq
C \|DT_{\widehat \De}(x,y)\| \leq C e^{2 r_\De(x)}$.
Thus if $(x,y,s) \in X$ then
$r_\De$ is $C e^{2\epsilon t/C_0}$-Lipschitz restricted
to the connected component
of the domain of $T_{\widehat \De}$ containing $(x,y)$, and
we conclude that if
$(x,y,s) \in X_\roof$ then $s \geq r_\De(x)-C e^{-\epsilon t/2}$, so
$\nu(X_\roof) \leq C e^{-\epsilon t/2}$.

Projecting $X_\lat$ on $(x,y)$,
we obtain a set $Z \subset \widehat \De$.
By Theorem \ref {exponentialtails}, $\nu(X_\lat) \leq C
e^{-\epsilon t/C}$ follows from
$\widehat m(Z) \leq C e^{-\epsilon t/C}$.  Let us show the latter estimate.
Using that $T_{\widehat \De}$, restricted to a connected component of its
domain intersecting $Z$, is
$C e^{2 \epsilon t/C_0}$-Lipschitz,
we get that $T_{\widehat \De}(Z)$ is
contained in a $C e^{2 \epsilon t/C_0} e^{-\epsilon t} \leq C e^{-\epsilon
t/2}$ neighborhood (with respect to the metric $\dist$) of the
boundary of $\widehat \De$.  Since $\widehat m$ is invariant and smooth,
and the boundary of $\widehat \De$ is
piecewise smooth, it follows that
$\widehat m(Z) \leq C e^{-\epsilon t/2}$.
\end{proof}

Notice that $\ln \frac {\dd\nu} {\dd z}$ is bounded over $\wde$, so
$\|f^{(t)}\|_{L^p(\nu)} \leq C \|f^{(t)}\|_{L^p(\dd z)} \leq C \|f \circ
P\|_{L^p(\dd z)} \leq C \|f \circ P\|_{L^p(\nu)}=C \|f\|_{L^p(\nu_{\CC^{(1)}})}
\leq C$.  Hence $\|f \circ P-f^{(t)}\|_{L^p(\nu)} \leq C$ and using Lemma
\ref {nuy} we conclude
that $\|\chi_{X \cup Y} (f \circ P-f^{(t)})\|_{L^{p'}(\nu)} \leq
C e^{-\epsilon t/C}$,
where $\chi_{X \cup Y}$ is the characteristic function of $X \cup Y$.
On the other hand, if $z_0 \in \wde \setminus (X \cup Y)$ and $\|z-z_0\|
\leq e^{-\epsilon t}/10$ then $\dist_F(z_0,z) \leq C e^{2 \epsilon t/C_0}
e^{-\epsilon t} \leq C e^{-\epsilon t/2}$. 
It follows that $|\psi^{(t)}(z)-f \circ P(z_0)|=|f \circ P(z)-f \circ
P(z_0)|<C e^{-\alpha\epsilon t/2}e^{k\epsilon t/ C_0}$.  Thus
$|f^{(t)}(z_0)-f \circ P(z_0)| \leq C e^{-\alpha \epsilon t/4}$.
This implies that
$\|\chi_{\wde \setminus (X \cup Y)} (f \circ P-f^{(t)})\|_{L^\infty(\nu)}
\leq C e^{-\epsilon t/C}$.  The result follows.
\end{proof}

Let now $k$, $\alpha$, $p$, $q$, $f$ and $g$ be as in Theorem \ref
{decay_teichmuller}. Let $\delta$ satisfy \eqref {tilde f}, and 
let $\epsilon_0$ be given by Lemma
\ref{lem:regularisation}.
Choose $p>p'>1$, $q>q'>1$ such that
$\frac {1} {p'}+\frac {1} {q'}=1$.  For $t>0$,
let $f^{(t)}$ and $g^{(t)}$ satisfy
\be \label {a}
\|f \circ P-f^{(t)}\|_{L^{p'}} \leq
C (\|f\|_{\DD_{k,\alpha}}+\|f\|_{L^p}) e^{-\epsilon_0 t},
\ee
\be \label {b}
\|f^{(t)}\|_{C^1(\wde)} \leq
C (\|f\|_{\DD_{k,\alpha}}+\|f\|_{L^p}) e^{\delta t},
\ee
\be \label {c}
\|g \circ P-g^{(t)}\|_{L^{q'}} \leq
C (\|g\|_{\DD_{k,\alpha}}+\|g\|_{L^q}) e^{-\epsilon_0 t},
\ee
\be \label {d}
\|g^{(t)}\|_{C^1(\wde)} \leq
C (\|g\|_{\DD_{k,\alpha}}+\|g\|_{L^q}) e^{\delta t}.
\ee
Then (\ref {tilde f}), (\ref {b}) and (\ref {d})
imply
\be \label {decay}
\left|\int f^{(t)} \cdot g^{(t)} \circ \widehat T_t \dd\nu
-\int f^{(t)} \dd\nu \int g^{(t)} \dd\nu\right| 
\leq C e^{-3 \delta t} \|f^{(t)}\|_{C^1} \|g^{(t)}\|_{C^1} \leq
C e^{-\delta t}.
\ee
We have
\begin{align}
\int f \cdot g \circ \TF_t \dd\nu_{\CC^{(1)}}-\int f \dd\nu_{\CC^{(1)}}\int g
\dd\nu_{\CC^{(1)}} &=\int f \circ P \cdot g \circ P
\circ \widehat T_t \dd\nu-\int f \circ P \dd\nu \int g \circ P \dd\nu\\
\nonumber
&=C_t(f \circ P,g \circ P).
\end{align}
Using (\ref {a}), (\ref {c}) and (\ref {decay}) we get
\begin{align}
|C_t(f \circ P,g \circ P)| &\leq |C_t(f^{(t)},g^{(t)})|+|C_t(f \circ
P-f^{(t)},g \circ P)|+|C_t(f \circ P,g \circ P-g^{(t)})|\\
\nonumber
&\phantom {\leq} +|C_t(f \circ P-f^{(t)},g \circ P-g^{(t)})|\\
\nonumber
&\leq |C_t(f^{(t)},g^{(t)})|+
2\|f \circ P-f^{(t)}\|_{L^{p'}} \|g\|_{L^{q'}}+
2\|f\|_{L^{p'}} \|g \circ P-g^{(t)}\|_{L^{q'}}\\
\nonumber
&\phantom {\leq}+
2\|f \circ P-f^{(t)}\|_{L^{p'}} \|g \circ P-g^{(t)}\|_{L^{q'}}\\
\nonumber
&\leq C e^{-\min(\delta,\epsilon_0) t} (\|f\|_{\DD_{k,\alpha}}+\|f\|_{L^p})
(\|g\|_{\DD_{k,\alpha}}+\|g\|_{L^q}).
\end{align}
This concludes the proof of Theorem \ref {decay_teichmuller}, modulo
Theorem \ref{exponentialtails} which will be proved in Sections
\ref{sec:distortion_estimate} and \ref{sec:proof_rec}. \qed

\subsection{A better recurrence estimate and the complement of large balls}

In the formulation of Theorem \ref {exponentialtails}, the particular
recurrence estimate is not necessarily good because we were more concerned
in obtaining not only a precompact transversal, but one for which the
combinatorics of the first return map is particularly simple (it is in
particular conjugate to a horseshoe on infinitely many symbols).  By
considering slightly more complicated combinatorics, one can get
considerably better estimates:

\begin{thm} \label {largeballs}

For every $\delta>0$, there exists a finite union $\widehat
Z=\bigcup \Delta_{\gamma_s} \times \Gamma_{\gamma_e}$ such that $\widehat
Z^{(1)}=\widehat
Z \cap \Up^{(1)}_\RR$ is precompact in
$\Up^{(1)}_\RR$, and the first return time $r_{\widehat Z}$ to
$\widehat Z$ under the Veech flow satisfies
\be
\int_{\widehat Z} e^{(1-\delta) r_{\widehat Z}} \dd\widehat m<\infty.
\ee

\end{thm}

This result easily implies Theorem \ref {compact_return}
(taking $K=\proj(\widehat Z^{(1)})$). It will be proved at the end of
Section \ref{sec:proof_rec}, by using a similar argument to the proof
of Theorem \ref{exponentialtails}.

\section{A distortion estimate}

\label{sec:distortion_estimate}

The proof of the recurrence estimates is based on the analysis of the Rauzy
renormalization map $R$.  The key step involves a control on the measure of
sets which present big distortion after some long (Teichm\"{u}ller) time.  In
order to obtain nearly optimal estimates, we will need to carry on a more
elaborate combinatorial analysis of Rauzy diagrams.

\subsection{Degeneration of Rauzy classes}

Let $\RR \subset \ssigma(\AA)$ be a Rauzy class.
Let $\AA' \subset \AA$ be a non-empty proper subset.

\begin{definition}

An arrow is called $\AA'$-colored if its winner belongs to
$\AA'$.  A path $\gamma \in \RR(\pi)$ is $\AA'$-colored if
it is a concatenation of $\AA'$-colored arrows.

\end{definition}

We call $\pi \in \RR$ {\it $\AA'$-trivial} if the last letters on both the
top and the bottom rows of $\pi$ do not belong to $\AA'$,
{\it $\AA'$-intermediate} if exactly one of those letters belong to $\AA'$
and {\it $\AA'$-essential} if both letters
belong to $\AA'$.  Alternatively, $\pi \in \RR$ is
trivial/intermediate/essential if it is the beginning (and ending) of
exactly $0$/$1$/$2$ $\AA'$-colored arrows.

An {\it $\AA'$-decorated Rauzy class} $\RR_* \subset \RR$ is a
maximal subset whose elements can be joined by an $\AA'$-colored
path.  We let $\Pi_*(\RR_*)$ be the set of all $\AA'$-colored
paths starting (and ending) at permutations in $\RR_*$.  We will
sometimes write $\Pi_*$ for $\Pi_*(\RR_*)$.

A decorated Rauzy class is called trivial if it
contains a trivial element $\pi$.  In this case $\RR_*=\{\pi\}$
and $\Pi_*(\RR_*)=\{\pi\}$ (recall that vertices are identified with trivial
(zero-length) paths).

A decorated Rauzy class is called essential if it contains an essential
element.

Since $\Pi_*(\RR_*) \neq \Pi(\RR)$ (for instance,
$\Pi_*(\RR_*)$ does not contain complete paths), any essential decorated
Rauzy class contains intermediate elements.

\subsubsection{Essential decorated Rauzy classes}

Let $\RR_*$ be an essential decorated Rauzy class.
Let $\RR^\ess_* \subset \RR_*$ be the set of essential elements of
$\RR_*$.  Let $\Pi^\ess_*(\RR_*) \subset \Pi_*(\RR_*)$ be the set of paths
which start and end at an element of $\RR^\ess_*$.

An {\it arc}
$\gamma \in \Pi_*(\RR_*)$ is a minimal non-trivial path in $\Pi^\ess_*$.
All arrows in an arc are of the same type and have the same winner, so the
type and winner of an arc are well defined.  Any
element of $\RR^\ess_*$ is thus the start (and end) of one top arc and
one bottom arc.  The losers in an arc are all distinct, moreover the first
loser is in $\AA'$ (and the others are not).

If $\gamma \in \Pi_*(\RR_*)$ is an arrow, then there exist
unique paths $\gamma_s,\gamma_e \in \Pi_*$ such that $\gamma_s \gamma
\gamma_e$ is an arc, called the {\it completion} of $\gamma$.  If $\pi$ is
intermediate, there is a single arc passing through $\pi$, the completion of
the arrow starting (or ending) at $\pi$.

If $\pi \in \RR_*$ we define $\pi^\ess \in \RR^\ess_*$ as follows.
If $\pi$ is essential then $\pi^\ess=\pi$.  If $\pi$ is intermediate,
let $\pi^\ess$ be the end of the arc passing through $\pi$.

To $\gamma \in \Pi_*$ we associate an element $\gamma^\ess \in \Pi^\ess_*$
as follows.  For a trivial path $\pi \in \RR_*$, we use
the previous definition of $\pi^\ess$.  Assuming that $\gamma$ is an arrow,
we distinguish two cases:
\begin{enumerate}
\item If $\gamma$ starts in an essential element, we let $\gamma^\ess$ be
the completion of $\gamma$,
\item Otherwise, we let $\gamma^\ess$ be the endpoint of the
completion of $\gamma$.
\end{enumerate}
We extend the definition to paths $\gamma \in \Pi_*$ by concatenation.
Notice that if $\gamma \in \Pi^\ess_*$ then $\gamma^\ess=\gamma$.

\subsubsection{Reduction}

We will now generalize the notion of simple reduction of \cite {AV}.  We
will need the following concept.

\begin{definition}

Given $\pi \in \sssigma(\AA)$ whose top and bottom rows end with
different letters, we obtain the {\it admissible end} of $\pi$ by deleting
as many letters from the beginning of the top and bottom rows of $\pi$ as
necessary to obtain an admissible permutation.  The resulting permutation
$\pi'$ belongs then to $\ssigma(\AA')$ for some $\AA' \subset \AA$.

\end{definition}

Let $\RR_*$ be an essential decorated Rauzy class, and let $\pi \in
\RR^\ess_*$.  Delete all letters not belonging to $\AA'$ from the top and
bottom rows of $\pi$.  The resulting permutation $\pi' \in \sssigma(\AA')$
is not necessary admissible, but since $\pi$ is essential the letters in
the end of the top and bottom rows of $\pi'$ are distinct.  Let $\pi^\red$
be the admissible end of $\pi'$.
We call $\pi^\red$ the {\it reduction} of $\pi$.

We extend the operation $\pi \mapsto \pi^\red$
of reduction from $\RR^\ess_*$ to the whole $\RR_*$ by taking the reduction
of an element $\pi \in \RR_*$ as the reduction of $\pi^\ess$.

If $\gamma \in \Pi^\ess_*$
is an arc, starting at $\pi_s$ and ending at $\pi_e$, then the reductions of
$\pi_s$ and $\pi_e$ belong to the same Rauzy class, and are joined by an
arrow $\gamma^\red$ (called the reduction of $\gamma$) of the
same type, same winner, and whose loser is the first loser of the arc
$\gamma$.  Thus the set of reductions of all
$\pi \in \RR_*$ is a Rauzy class $\RR^\red_* \subset
\ssigma(\AA'')$ for some $\AA'' \subset \AA'$.

We define the reduction of a path $\gamma \in \Pi_*$ as follows.  If
$\gamma$ is a trivial path or an arc, it is defined as above.  We extend the
definition to the case $\gamma \in \Pi^\ess_*$ by concatenation.  In general
we let the reduction of $\gamma$ to be equal to the reduction of
$\gamma^\ess$.

Notice that the reduction map $\RR^\ess_* \to \RR^\red_*$ is a bijection.
The reduction map $\Pi^\ess_* \to \Pi(\RR^\red_*)$ is a bijection compatible
with concatenation.

\subsection{Further combinatorics}

Let $\AA' \subset \AA$ be a non-empty proper subset.

\subsubsection{Drift in essential decorated Rauzy classes}

Let $\RR_* \subset \RR$ be an essential $\AA'$-decorated Rauzy class.

For $\pi \in \RR_*$, let $\alpha_t(\pi)$ (respectively, $\alpha_b(\pi)$)
be the rightmost letter in the top (respectively, bottom) row
of $\pi$ that belongs to $\AA \setminus \AA'$.  Let $d_t(\pi)$
(respectively, $d_b(\pi)$) be the position of $\alpha_t(\pi)$ (respectively,
$\alpha_b(\pi)$) in the top (respectively, bottom) of $\pi$.  Let
$d(\pi)=d_t(\pi)+d_b(\pi)$.

An essential element of $\RR_*$ is thus some $\pi$ such that
$d_t(\pi),d_b(\pi)<d$.
If $\pi_s$ is an essential element of $\RR_*$ and $\gamma \in \Pi_*(\RR_*)$
is an arrow starting at $\pi_s$ and ending at $\pi_e$, then
\begin{enumerate}

\item
$d_t(\pi_e)=d_t(\pi_s)$ or $d_t(\pi_e)=d_t(\pi_s)+1$, the second possibility
happening if and only if $\gamma$ is a bottom whose winner precedes
$\alpha_t(\pi_s)$ in the top of $\pi_s$.

\item
$d_b(\pi_e)=d_b(\pi_s)$ or $d_b(\pi_e)=d_b(\pi_s)+1$, the second possibility
happening if and only if $\gamma$ is a top whose winner precedes
$\alpha_b(\pi_s)$ in the bottom of $\pi_s$.

\end{enumerate}
In particular $d(\pi_e)=d(\pi_s)$ or $d(\pi_e)=d(\pi_s)+1$.  In the second
case, we say that $\gamma$ is {\it drifting}.

Let $\RR^\red_*$ be the reduction of $\RR_*$, so that $\RR^\red_* \subset
\ssigma(\AA'')$ for some $\AA'' \subset \AA'$.  If $\pi \in \RR_*$ is
essential then there exists $\alpha \in \AA''$ that either precedes
$\alpha_t(\pi)$ in the top of $\pi$ or precedes $\alpha_b(\pi)$ in the
bottom of $\pi$ (we call such an $\alpha$ {\it good} for $\pi$).  Indeed, if
$\gamma \in \Pi_*(\RR_*)$ is a path starting at $\pi$, ending with
a drifting arrow and minimal with this property then the winner of the last
arrow of $\gamma$ belongs to $\AA''$ and either precedes $\alpha_t(\pi)$
in the top of $\pi$ (if the drifting arrow is a bottom) or precedes
$\alpha_b(\pi)$ in the bottom of $\pi$ (if the drifting arrow is a top).

Notice that if $\gamma \in \Pi_*$ is an arrow starting and ending at
essential elements $\pi_s,\pi_e$, then a good letter for $\pi_s$ is also a
good letter for $\pi_e$.  Moreover, if $\gamma$ is not drifting then the
winner of $\gamma$ is not a good letter for $\pi_s$.

\subsubsection{Standard decomposition of separated paths}

\begin{definition}

An arrow is called $\AA \setminus \AA'$-separated if both its winner and
loser belong to $\AA'$.  A path $\gamma \in \RR$ is
$\AA \setminus \AA'$-separated if
it is a concatenation of $\AA \setminus \AA'$-separated arrows.

\end{definition}

If $\gamma \in \Pi(\RR)$ is a non-trivial maximal
$\AA \setminus \AA'$-separated path, then there
exists an essential $\AA'$-decorated Rauzy class
$\RR_* \subset \RR$ such that $\gamma \in
\Pi_*(\RR_*)$.  Moreover, if
$\gamma=\gamma_1...\gamma_n$, then each $\gamma_i$
starts at an essential element $\pi_i \in \RR_*$ (and $\gamma_n$ ends at an
intermediate element of $\RR_*$ by maximality).

Let $r=d(\pi_n)-d(\pi_1)$.  Let $\gamma=\gamma^{(1)} \gamma^1...\gamma^{(r)}
\gamma^r$ where the $\gamma^i$ are drifting arrows
and $\gamma^{(i)}$ are (possibly trivial)
concatenations of non-drifting arrows.  If $\alpha$ is a good letter for
$\pi_1$, then it follows that $\alpha$ is not the winner of any arrow in any
$\gamma^{(i)}$.  The reduction of the $\gamma^{(i)}$ are thus non-complete
paths in $\Pi(\RR^\red_*)$, according to Definition \ref {as defined}.

\subsection{The distortion estimate} \label {distortion}

The distortion argument will involve not only the study of Lebesgue
measure, but also of its forward images under the renormalization
map. Technically, this is most
conveniently done by introducing a class of measures which is
invariant as a whole. For $q\in \R_+^\AA$, let 
$\e_q=\{\lambda \in \R^\AA_+ \tq \langle \lambda, q \rangle<1\}$. If
$\nu_q$ denotes the measure on $\P\R_+^\AA$ given by
$\nu_q(A)=\Leb(\R_+A \cap \Lambda_q)$, then
  \begin{equation}
  \nu_q(B_\gamma^* A)= \Leb( B_\gamma^* (\R_+A) \cap
  \Lambda_q)
  =\Leb( \R_+A \cap \Lambda_{B_\gamma\cdot q}) = \nu_{B_\gamma\cdot q}(A).
  \end{equation}
An important point of the discussion to follow is that we will work at
the same time with all parameters $q$, and obtain estimates which are
\emph{uniform in $q$}. In fact, we will not really study $\nu_q$,
rather the quantities $\nu_q(\Delta_\gamma)$ for $\gamma\in \Pi(\RR)$.
The following notations make it possible to do so in a formalism where
conditioning is more or less transparent.

Let $\RR \subset \ssigma(\AA)$ be a Rauzy class and let
$\gamma \in \Pi(\RR)$.  We let $\e_{q,\gamma}=B_\gamma^*
\cdot \e_{B_\gamma \cdot q}$.  This definition is such that $\{\pi\} \times
\e_{q,\gamma}=(\{\pi\} \times \e_q) \cap \Delta_\gamma$ where
$\pi$ is the start of $\gamma$.

For $\AA' \subset \AA$, $q \in \R^\AA_+$, let
$\n_{\AA'}(q)=\prod_{\alpha \in \AA'} q_\alpha$, $\n(q)=\n_\AA(q)$.
Then $\Leb(\e_q)=\frac {1} {(\# \AA)! \n(q)}$, where $\Leb$ is the usual
Lebesgue measure on $\R^\AA$.  This
gives $\frac {\Leb(\e_{q,\gamma})} {\Leb(\e_q)}=\frac {\n(q)}
{\n(B_\gamma \cdot q)}$.

For $\AA' \subset \AA$ non-empty, let $\M_{\AA'}(q)=\max_{\alpha \in \AA'}
q_\alpha$.  Let $\M(q)=\M_\AA(q)$.

If $\Gamma \subset \Pi(\RR)$ is a set of paths starting at the same $\pi
\in \RR$, let
$\e_{q,\Gamma}=\bigcup_{\gamma \in \Gamma} \e_{q,\gamma}$.

Given $\Gamma \subset \Pi(\RR)$, $\gamma_s \in
\Pi(\RR)$, let $\Gamma_{\gamma_s} \subset \Gamma$ be the set of paths
starting by $\gamma_s$, and let $\Gamma^{\gamma_s}$ be the collection of
ends $\gamma_e$ of paths $\gamma=\gamma_s \gamma_e \in \Gamma$.

Let $P_q(\Gamma \di \gamma)=\frac {\Leb(\e_{q,\Gamma_\gamma})}
{\Leb(\e_{q,\gamma})}$.  If $\pi$ is the end
of $\gamma$, we have
$P_q(\Gamma \di \gamma)=P_{B_\gamma \cdot q}(\Gamma^\gamma \di \pi)$.  If
$\gamma$ is an arrow starting at $\pi$
with winner $\alpha$ and loser $\beta$, we have
\be
P_q(\gamma \di \pi)=\frac {q_\beta} {q_\alpha+q_\beta}.
\ee

A family $\Gamma_s \subset \Pi(\RR)$ is called disjoint if no two elements
are comparable (for the partial order defined in \S \ref {rauzy c}).
If $\Gamma_s$ is disjoint and $\Gamma \subset \Pi(\RR)$ is
a family such that any $\gamma \in \Gamma$ starts by some element $\gamma_s
\in \Gamma_s$, then for every $\pi \in \RR$
\be
P_q(\Gamma \di \pi)=\sum_{\gamma_s \in \Gamma_s} P_q(\Gamma \di \gamma_s)
P_q(\gamma_s \di \pi) \leq P_q(\Gamma_s \di \pi) \sup_{\gamma_s \in \Gamma_s}
P_q(\Gamma \di \gamma_s).
\ee

The key distortion estimate is the following.

\begin{thm} \label {4}

There exist $C>0$, $\theta>0$, depending only on $\# \AA$
with the following property.  Let $\AA' \subset \AA$ be a non-empty proper
subset, $0 \leq m \leq M$ be integers, $q \in \R^\AA_+$.
Then for every $\pi \in \RR$,
\be
P_q(\gamma \in \Pi(\RR),\, \M(B_\gamma \cdot q)>2^M \M(q) \text { and }
\M_{\AA'}(B_\gamma \cdot q)<2^{M-m} \M(q) \di \pi) \leq C (m+1)^\theta 2^{-m}.
\ee

\end{thm}

The proof is based on induction on
$\# \AA$, and will take the remaining of this section.

\subsection{Reduction estimate}

Let $\RR_*$ be a $\AA'$-decorated Rauzy class, and let
$\gamma \in \Pi_*(\RR_*)$ start at
$\pi \in \RR_*$.  If $\RR_*$ is essential, let
$\RR^\red_* \subset \ssigma(\AA'')$ be its reduction.  Let $q^\red$ be the
(canonical) projection of $q$ on $\R^{\AA''}$ (obtained by forgetting the
coordinates in $\AA \setminus \AA''$).  Then the
projection of $B_\gamma \cdot q$ on $\R^{\AA''}$ coincides with
$B_{\gamma^\red} \cdot q^\red$.  Notice also that the
projections of $q$ and $B_\gamma \cdot q$ on $\R^{\AA' \setminus \AA''}$
coincide.  This gives the formula
\be \label {reduction}
\frac {P_q(\gamma \di \pi)} {P_{q^\red}(\gamma^\red \di \pi^\red)}=\frac
{\n_{\AA \setminus \AA'}(q)}
{\n_{\AA \setminus \AA'}(B_\gamma \cdot q)}.
\ee

\begin{prop} \label {reduction estimate}

Let $\RR_*$ be an $\AA'$-decorated Rauzy class, and let $\Gamma \subset
\Pi_*(\RR_*)$ be a family of paths such that, for all $\gamma\in \Gamma$,
$\n_{\AA \setminus \AA'}(B_\gamma \cdot q)
\geq 2^M \n_{\AA \setminus \AA'}(q)$.  Then for every $\pi \in \RR_*$,
\be
P_q(\Gamma \di \pi) \leq 2^{-M}.
\ee

\end{prop}

\begin{proof}

We may assume that $\Gamma$ is the collection of all minimal paths $\gamma
\in \Pi_*(\RR_*)$ starting at $\pi$ and
satisfying $\n_{\AA \setminus \AA'}(B_\gamma \cdot q)
\geq 2^M \n_{\AA \setminus \AA'}(q)$.
If $\RR_*$ is trivial then either $\Gamma$ is empty or
$M=0$ and the estimate is obvious.  If $\RR_*$ is
neither trivial nor essential, then $\Gamma$ consists of a
single path $\gamma$, and the result follows from the definition of
$P_q(\gamma \di \pi)$.
If $\RR_*$ is essential, we notice that two distinct
paths in $\Gamma$ have disjoint reductions, so the estimate follows from
(\ref {reduction}).
\end{proof}

\subsection{The main induction scheme}

\begin{definition}

A path $\gamma \in \Pi(\RR)$
is called $\AA'$-preferring if it is a concatenation of a
$\AA'$-separated path (first) and a $\AA'$-colored path (second).
\end{definition}

A path is $\AA'$-preferring if it has no loser in $\AA'$. 
Notice that $\gamma$ is $\AA'$-preferring if and only if
$\langle B_\gamma \cdot e_\alpha,e_\beta \rangle=0$ for $\alpha \in \AA
\setminus \AA'$,
$\beta \in \AA'$ (so $B_\gamma$ is block-triangular).

Notice also that the $\AA'$-separated part or the $\AA'$-colored part in an
$\AA'$-preferring path may
very well be trivial.

\begin{prop} \label {3}

There exist $C>0$, $\theta>0$, depending only on $\# \AA$
with the following property.  Let $\AA' \subset \AA$ be a non-empty proper
subset, $M \in \N$, $q \in \R^\AA_+$.
Then for every $\pi \in \RR$,
\be
P_q(\gamma \text { is $\AA \setminus \AA'$-separated and }
\M_{\AA'}(B_\gamma \cdot q)>2^M \M(q) \di \pi)
\leq C (M+1)^\theta 2^{-M}.
\ee

\end{prop}

\begin{prop} \label {5}

There exist $C>0$, $\theta>0$, depending only on $\# \AA$
with the following property.  Let $\AA' \subset \AA$ be a non-empty proper
subset, $M \in \N$, $q \in \R^\AA_+$.
Then for every $\pi \in \RR$,
\be
P_q(\gamma \text { is $\AA'$-preferring and }
\M_{\AA'}(B_\gamma \cdot q) \leq 2^M \M(q)<\M(B_\gamma \cdot q) \di \pi)
\leq C (M+1)^\theta 2^{-M}.
\ee

\end{prop}

\begin{prop} \label {2}

There exist $C>0$, $\theta>0$, depending only on $\# \AA$ with the
following property.  Let $M \in \N$, $q \in \R^\AA_+$.
Then for every $\pi \in \RR$,
\be
P_q(\gamma \text { is not complete and }
\M(B_\gamma \cdot q)>2^M \M(q) \di \pi)
\leq C (M+1)^\theta 2^{-M}.
\ee

\end{prop}

The proof of Theorem \ref {4} and Propositions \ref {3},
\ref {5} and \ref {2} will be
carried out simultaneously in an induction argument on $d=\#\AA$.  For $d
\geq 2$, consider the statements:
\begin{enumerate}
\item [(A$_d$)] Proposition \ref {2} holds for $\#\AA=d$,
\item [(B$_d$)] Proposition \ref {3} holds for $\#\AA=d$,
\item [(C$_d$)] Proposition \ref {5} holds for $\#\AA=d$,
\item [(D$_d$)] Theorem \ref {4} holds for $\#\AA=d$.
\end{enumerate}
The induction step will be composed of four parts:
\begin{enumerate}
\item (A$_j$), $2 \leq j<d$, implies (B$_d$),
\item (B$_d$) implies (C$_d$),
\item (C$_d$) implies (D$_d$),
\item (D$_j$), $2 \leq j \leq d$, implies (A$_d$).
\end{enumerate}
Notice that the start of the induction is trivial (for $d=2$ the
hypothesis in (1) is trivially satisfied).

In what follows, $C$ and $\theta$ denote generic constants, whose actual
value may vary during the course of the proof.

\subsubsection{Proof of (1)}

Let $\Gamma$ be the set of all maximal $\AA \setminus \AA'$-separated
$\gamma$ starting at $\pi$
such that $\M_{\AA'}(B_\gamma \cdot q)>2^M \M(q)$. By Lemma \ref{Lemme_MMY},
it is sufficient to
prove 
\be
\label{foobar}
P_q(\Gamma \di \pi) \leq C (M+1)^\theta 2^{-M}.
\ee

If $\Gamma$ is non-empty then $\pi$ is essential (and if $\Gamma=\emptyset$
the statement is trivial).
Let $\RR_*$ be the $\AA'$-decorated class containing $\pi$.  We have $\Gamma
\subset \Pi_*(\RR_*)$.  Decompose $\Gamma$ into
subsets $\Gamma_{\overline M}$, $\overline M \geq M$,
containing the $\gamma \in \Gamma$ with $2^{\overline M+1} \geq
\M_{\AA'}(B_\gamma \cdot q)>2^{\overline M} \M(q)$.
Recall the decomposition of $\gamma \in \Gamma$, $\gamma=\gamma^{(1)}
\gamma^1...\gamma^{(r)} \gamma^r$ where $r=r(\gamma)<2d$.
Let $\Gamma_{\overline M,r} \subset \Gamma_{\overline M}$
collect the $\gamma$ with $r(\gamma)=r$.  Let
$\gamma_{(i)}=\gamma^{(1)} \gamma^1...\gamma^{(i)} \gamma^i$,
$1 \leq i \leq r$, and let $\gamma_{(0)}$ be the start of $\gamma$.
To $\gamma \in \Gamma_{\overline{M},r}$ we associate $\mm=(m_1,...,m_r)$ where
\be
2^{m_i} \leq \frac
{\max \{\M(q),\M_{\AA'}(B_{\gamma_{(i-1)} \gamma^{(i)}} \cdot q)\}}
{\max \{\M(q),\M_{\AA'}(B_{\gamma_{(i-1)}} \cdot q)\}}
<2^{m_i+1}.
\ee
We have $2^{\sum m_i} \M(q) \leq \M_{\AA'}(B_\gamma \cdot q)
\leq 2^{2r+\sum m_i} \M(q)$, so $\overline M+1 \geq \sum m_i \geq \overline
M-2r$.
Let $\Gamma_{\overline M,r,\mm}$ collect the $\gamma$ with the same
$\mm$.  For $0 \leq i \leq r$,
let $\Gamma_{\overline M,r,\mm,i}$ be the collection of
all possible $\gamma_{(i)}$.

Let $\RR^\red_* \subset \ssigma(\AA'')$ be the reduction of $\RR_*$.  If
$\gamma_s \in \Pi_*(\RR_*)$ is $\AA \setminus \AA'$-separated then
\be \label {mm,i}
P_q(\Gamma_{\overline M,r,\mm,i} \di \gamma_s)=
P_{q^\red}(\Gamma^\red_{\overline M,r,\mm,i} \di \gamma^\red_s),
\ee
where $q^\red$ is the orthogonal projection of $q$ on $\R^{\AA''}$,
$\Gamma^\red_{\overline M,r,\mm,i}$ is the image of
$\Gamma_{\overline M,r,\mm,i}$ by the reduction
map and $\gamma^\red_s$ is the reduction of $\gamma_s$.  If $\gamma \in
\Gamma_{\overline M,r,\mm,i}$ starts by
$\gamma_s \in \Gamma_{\overline M,r,\mm,i-1}$ then we can
write $\gamma=\gamma_s \gamma_a \gamma_b$, where $\gamma_b$ is a drifting
arrow, and $\gamma_a$ is a concatenation of non-drifting arrows.  Then
$\gamma^\red_a$ is a non-complete path (in $\Pi(\RR^\red_*)$) satisfying
$\M_{\AA''}(B_{\gamma^\red_a}\cdot B_{\gamma^\red_s} \cdot q^\red)
\geq 2^{m_i-1} 
\M_{\AA''}(B_{\gamma^\red_s} \cdot q^\red)$.
By (A$_j$) with $j=\#\AA''<d$,
\be \label {mm,i1}
P_{q^\red}(\Gamma^\red_{\overline M,r,\mm,i} \di \gamma^\red_s)
\leq C (m_i+1)^\theta 2^{-m_i}, \quad \gamma_s \in
\Gamma_{\overline M,r,\mm,i-1}.
\ee
Each family $\Gamma_{\overline M,r,\mm,i}$ is disjoint,
so (\ref {mm,i}) and (\ref
{mm,i1}) imply
\be \label {mm,i2}
P_q(\Gamma_{\overline M,r,\mm,i} \di \pi)\leq C (m_i+1)^\theta 2^{-m_i}
P_q(\Gamma_{\overline M,r,\mm,i-1} \di \pi),
\ee
which gives
\be
P_q(\Gamma_{\overline M,r,\mm} \di \pi)=
P_q(\Gamma_{\overline M,r,\mm,r} \di \pi)
\leq \prod_{i=1}^r C (m_i+1)^\theta 2^{-m_i} \leq
C(\overline M+1)^\theta 2^{-\overline M}.
\ee
Summing over the different $\mm$ (with $\sum m_i \leq \overline M+1$),
$r<2d$, and $\overline M \geq M$, we get \eqref{foobar}.

\subsubsection{Proof of (2)}

Let $\Gamma$ be the set of all $\AA'$-preferring
$\gamma$ such that $\M_{\AA'}(B_\gamma \cdot q) \leq 2^M
\M(q)<\M(B_\gamma \cdot q)$, and which are minimal with those
properties.
Any $\gamma \in \Gamma$ is of the form $\gamma=\gamma_s \gamma_e$ where
$\gamma_s$ is $\AA'$-separated and $\gamma_e$ is $\AA'$-colored.  Let
$\Gamma_s \subset \Pi(\RR)$ collect all possible $\gamma_s$.  Notice that
$\Gamma_s$ is disjoint.

Let $m=m(\gamma_s) \in [-1,M]$ be the smallest integer such that $\M_{\AA
\setminus \AA'}(B_{\gamma_s} \cdot q) \leq 2^{m+1} \M(q)$.  Notice that
$\M_{\AA'}(B_{\gamma_s} \cdot q)=\M_{\AA'}(q)
\leq \M(q)$.  Let $\Gamma_{s,m}$ collect all $\gamma_s \in \Gamma_s$ with
$m(\gamma_s)=m$.

Let us show that for $\gamma_s \in \Gamma_{s,m}$
\be \label {s}
P_q(\Gamma \di \gamma_s) \leq 2^{m+1-M}.
\ee
Let $\pi_e$ be the ending of $\gamma_s$.  Let $\Gamma^{\gamma_s}$ be the
set of all endings $\gamma_e$ of paths $\gamma=\gamma_s \gamma_e \in \Gamma$
that begin with $\gamma_s$.  Let $\RR_*$ be the $\AA'$-decorated Rauzy
class containing $\pi_e$.
Then $\Gamma^{\gamma_s} \subset \Pi_*(\RR_*)$ is a collection
of paths $\gamma_e$ satisfying $\M_{\AA \setminus \AA'}
(B_{\gamma_e} \cdot B_{\gamma_s} \cdot q)>2^M \M(q) \geq 2^{M-1-m}
\M(B_{\gamma_s} \cdot q)$.  By Proposition \ref {reduction estimate},
$P_q(\Gamma \di \gamma_s)=P_{B_{\gamma_s} \cdot q}(\Gamma^{\gamma_s} \di \pi_e)
\leq 2^{m+1-M}$.

If $m \geq 0$ then $\Gamma_{s,m}$ consists of $\AA'$-separated paths
$\gamma_s$ with $\M_{\AA \setminus \AA'}(B_{\gamma_s} \cdot q)>2^m \M(q)$.
By (B$_d$),
\be \label {s,m}
P_q(\Gamma_{s,m} \di \pi) \leq C(m+2)^\theta 2^{-m}.
\ee
Notice that (\ref {s,m}) is still satisfied (trivially) for $m=-1$.
Putting together \eqref {s,m} and \eqref {s}, and summing over $m$,
we get
\be
P_q(\Gamma \di \pi) \leq C (M+1)^\theta 2^{-M}.
\ee

\subsubsection{Proof of (3)}

The proof is by descending recurrence on $\# \AA'$.  We may
assume that $m>0$ since the case $m=0$ is trivial.  Let $\Gamma \subset
\Pi(\RR)$ be the set of $\gamma$ starting at $\pi$ and such
that $\M(B_\gamma \cdot q)>2^M \M(q)$,
$\M_{\AA'}(B_\gamma \cdot q)<2^{M-m}\M(q)$ and which are minimal with those
properties.  We want to
estimate $P_q(\Gamma \di \pi) \leq C (m+1)^\theta 2^{-m}$.

Let $\Gamma_D \subset \Gamma$ be the set of $\AA'$-preferring paths.
We have $P_q(\Gamma_D \di \pi) \leq C
(M+1)^\theta 2^{-M}$ by (C$_d$),
so we just have to prove that
$P_q(\Gamma \setminus \Gamma_D \di \pi) \leq C (m+1)^\theta 2^{-m}$.

If $\gamma \in \Gamma \setminus \Gamma_D$, then at least one of the arrows
composing $\gamma$ has as winner an element of $\AA \setminus \AA'$, and as
loser an element of $\AA'$.  Decompose $\gamma=\gamma_s \gamma_e$ with
$\gamma_s$ maximal such that no arrow composing $\gamma_e$ has as winner
an element of $\AA \setminus \AA'$, and as loser an element of $\AA'$; let
$n_0=n_0(\gamma)$ be the length of $\gamma_s$.  Let
$\beta=\beta(\gamma) \in \AA \setminus
\AA'$ be the winner of the last arrow of $\gamma_s$.

We can then write $\Gamma \setminus \Gamma_D$ as the union of
$\Gamma_\beta$, $\beta \in \AA \setminus \AA'$, where $\Gamma_\beta$
collects all $\gamma$ with $\beta(\gamma)=\beta$.  We only have to prove
that $P_q(\Gamma_\beta \di \pi) \leq C (m+1)^\theta 2^{-m}$
for any $\beta \in \AA \setminus \AA'$.

Let $\Gamma^*_\beta \subset \Gamma_\beta$
be the set of $\gamma$ such that $\M(B_{\gamma_s} \cdot q) \leq 2^{M-m}
\M(q)$.  For $\gamma \in \Gamma^*_\beta$, write $\gamma=\gamma^*_s
\gamma^*_e$ with $\gamma^*_s$ minimal with $\M(B_{\gamma^*_s} \cdot
q)>2^{M-m}\M(q)$.  In particular
$\M(B_{\gamma^*_s} \cdot q) \leq 2^{M+1-m}\M(q)$.  Let $n^*(\gamma)$ be the
length of $\gamma^*_s$.  Since $n^*>n_0$, $\gamma^*_e$ is $\AA'$-preferring.
Notice that $\gamma^*_e$ is also such that $\M(B_{\gamma^*_e} \cdot
B_{\gamma^*_s} \cdot q)>2^{m-1}\M(B_{\gamma^*_s} \cdot q)$.
By (C$_d$), it follows that
$P_q(\Gamma^*_\beta \di \gamma^*_s) \leq C m^\theta 2^{1-m}$.  Since the
collection of all possible $\gamma^*_s$ is disjoint, we get
\be
P_q(\Gamma^*_\beta \di \pi) \leq C m^\theta 2^{1-m}.
\ee
Thus we only need to show that $P_q(\Gamma_\beta \setminus
\Gamma^*_\beta \di \pi) \leq C m^\theta 2^{-m}$.

Before continuing, let us notice that if $\# \AA'=\# \AA-1$, then
$\Gamma_\beta=\Gamma^*_\beta$.  Indeed in this case
$\AA=\AA' \cup \{\beta\}$, and since $\M_\beta(B_{\gamma_s} \cdot q) \leq
\M_{\AA'}(B_{\gamma_s} \cdot q)$, we have $\M(B_{\gamma_s} \cdot q) \leq
2^{M-m} \M(q)$.  In particular, the previous argument is enough to establish
(D$_d$) in the case $\# \AA'=d-1$, which
allows us to start the
reverse induction on $\# \AA'$ used in the argument below.

For $\gamma \in \Gamma_\beta \setminus \Gamma^*_\beta$, there exists
an integer $m_0=m_0(\gamma) \in [0,m)$ such that $2^{M-m_0} \M(q) \geq
\M(B_{\gamma_s} \cdot q)>2^{M-1-m_0} \M(q)$.  We collect all $\gamma$ with
$m_0(\gamma)=m_0$ in $\Gamma_{\beta,m_0}$.  It is enough to show that
\be \label {P_q}
P_q(\Gamma_{\beta,m_0} \di \pi) \leq C (m+1)^\theta 2^{-m}.
\ee

Write $\gamma=\gamma^1_s
\gamma^1_e=\gamma^2_s \gamma^2_e$ where $\gamma^1_s$, $\gamma^2_s$ are
minimal such that $\M(B_{\gamma^1_s} \cdot q)>2^{M-m_0} \M(q)$,
$\M(B_{\gamma^2_s} \cdot q)>2^{M-1-m_0} \M(q)$.
Let $n_1=n_1(\gamma)$ and $n_2=n_2(\gamma)$ be the
lengths of $\gamma^1_s$ and $\gamma^2_s$.  We have
$n_2<n_0<n_1$.\footnote{Notice that we can not have $n_2=n_0$, since
otherwise $\M_{\AA'}(B_{\gamma_s} \cdot q)=\M(B_{\gamma_s} \cdot
q)>2^{M-m} \M(q)$, so that $\gamma \notin \Gamma$.}

Let $\Gamma^1_{\beta,m_0,s}$, $\Gamma^2_{\beta,m_0,s}$ collect all possible
paths $\gamma^1_s$, $\gamma^2_s$ as above.
The families $\Gamma^1_{\beta,m_0,s}$, $\Gamma^2_{\beta,m_0,s}$ are
disjoint.  If $\gamma=\gamma^1_s \gamma^1_e
\in \Gamma_{\beta,m_0,s}$ with $\gamma^1_s \in \Gamma^1_{\beta,m_0,s}$,
the path
$\gamma^1_e$ is $\AA'$-preferring and satisfies $\M(B_{\gamma^1_e} \cdot
B_{\gamma^1_s} \cdot q)>2^{m_0-1} \M(B_{\gamma^1_s} \cdot q)$,
$\M_{\AA'}(B_{\gamma^1_e} \cdot B_{\gamma^1_s} \cdot q)<2^{M-m} M(q)<
\M(B_{\gamma^1_s} \cdot q)$, so by (C$_d$) we have
\be \label {P_q1}
P_q(\Gamma_{\beta,m_0,s} \di \gamma^1_s) \leq C (m_0+1)^\theta 2^{-m_0}, \quad
\gamma^1_s \in \Gamma^1_{\beta,m_0,s}.
\ee
On the other hand, $\M_\beta(B_{\gamma^2_s} \cdot q)<
\M_{\AA'}(B_{\gamma_s} \cdot q)<2^{M-m} \M(q)$ so that $M_{\AA' \cup
\{\beta\}}(B_{\gamma^2_s} \cdot q)<2^{M-m} \M(q)$.
Then
\be \label {P_q2}
P_q(\Gamma^1_{\beta,m_0,s} \di \pi) \leq
P_q(\Gamma^2_{\beta,m_0,s} \di \pi) \leq C (m-m_0)^\theta 2^{m_0+1-m},
\ee
where the first
inequality is trivial and the second is by the reverse induction hypothesis
(that is, (D$_d$)
with $\AA' \cup \{\beta\}$ instead of $\AA'$).  Since
$\Gamma^1_{\beta,m_0,s}$ is disjoint,
(\ref {P_q1}) and (\ref {P_q2}) imply (\ref {P_q}).

\subsubsection{Proof of (4)}

Let $\gamma \in \Pi(\RR)$ be a non-complete path starting at $\pi$.
Let $\beta \in \AA$ be a letter which is not winner of any arrow of $\gamma$,
and let $\AA'=\AA \setminus \{\beta\}$.
If $\RR_* \subset \RR$ is the $\AA'$-decorated Rauzy class containing $\pi$
then $\gamma \in \Pi_*(\RR_*)$.
Let $\Gamma_\beta \subset \Pi_*(\RR_*)$ be the family of paths
$\gamma$ satisfying $\M(B_\gamma \cdot q)>2^M \M(q)$ and minimal with this
property.  It is enough to show that
\be \label {Gamma,pi}
P_q(\Gamma_\beta \di \pi) \leq C (M+1)^\theta 2^{-M}
\ee
for an arbitrary choice of $\beta$ and $\RR_*$.

First notice that $\RR_*$ can not be a trivial decorated Rauzy class, since
$\AA \setminus \AA'$ has a single element.  If $\RR_*$
is neither trivial nor essential, then $\Gamma_\beta$ contains a
unique path $\gamma$ starting at $\pi$.  In this case
$P_q(\Gamma_\beta \di \pi)=P_q(\gamma \di \pi)<2^{-M}$.  It is enough then to
consider the case where $\RR_*$ is essential.

Let $\Gamma^*_\beta \subset \Gamma_\beta$ be the set of all $\gamma$
such that $\M_\beta(B_\gamma \cdot q) \leq \M(q)$.  By (D$_d$) applied
to $\{\beta\}$, we have
\be \label {beta,M}
P_q(\Gamma^*_\beta \di \pi) \leq C (M+1)^\theta 2^{-M}.
\ee

For $\gamma \in \Gamma_\beta \setminus \Gamma^*_\beta$,
there is at least one arrow composing $\gamma$ with $\beta$ as loser.
Let $\alpha=\alpha(\gamma)$ be the winner of the last such arrow.
Let $m_0=m_0(\gamma) \in
[0,M]$ be such that $2^{m_0} \M(q)<\M_\beta(B_\gamma \cdot q) \leq 2^{m_0+1}
\M(q)$.
Write $\gamma=\gamma_s \gamma_e$ where $\gamma_s$ is minimal with
$\M_\beta(B_{\gamma_s} \cdot q)>2^{m_0} \M(q)$.
Let $M_0=M_0(\gamma) \in [m_0,M]$ be such that
$2^{M_0} \M(q)<\M(B_{\gamma_s} \cdot q) \leq 2^{M_0+1} \M(q)$.
Let $\Gamma \subset \Gamma_\beta \setminus \Gamma^*_\beta$ collect
the $\gamma$ with the same $\alpha$, $m_0$ and $M_0$.  It is enough to show
that
\be \label {Gamma,M}
P_q(\Gamma \di \pi) \leq C (M+1)^\theta 2^{-M}.
\ee

Let $\Gamma_s$ be the family of possible $\gamma_s$ for $\gamma \in \Gamma$.
By (D$_d$) applied to $\{\beta\}$,
\be \label {M-M0}
P_q(\Gamma \di \gamma_s) \leq C (M+1-M_0)^\theta 2^{M_0-M},
\quad \gamma_s \in \Gamma_s.
\ee

Let $\RR^\red_* \subset \ssigma(\AA'')$ be the reduction of $\RR_*$.
Notice that two distinct
paths in $\Gamma_s$ have disjoint reductions.  Let
$\Gamma^\red_s \subset \Pi(\RR^\red_*)$ be the image of
$\Gamma_s$ by the reduction map.  Let $q^\red$ be the
canonical projection of $q$ on $\R^{\AA''}$.  Then by (\ref {reduction}),
\be \label {-m0}
P_q(\Gamma_s \di \pi) \leq
P_q(\Gamma^\red_s \di \pi^\red) \sup_{\gamma_s \in
\Gamma_s} \frac {\n_\beta(q)} {\n_\beta(B_{\gamma_s} \cdot
q)} \leq P_q(\Gamma^\red_s \di \pi^\red) 2^{-m_0}.
\ee
Notice that if $\gamma^\red_s \in \Gamma^\red_s$ then
$\M_\alpha(B_{\gamma^\red_s} \cdot q^\red) \leq 2^{m_0+1} \M(q)$, and if
$M_0>m_0$ we also have
$\M(B_{\gamma^\red_s} \cdot q^\red)>2^{M_0} \M(q)$.  Thus, if $M_0>m_0$,
by (D$_j$) with $j=\#\AA''<d$,
\be \label {m0=M0}
P_q(\Gamma^\red_s \di \pi^\red) \leq C (M_0+1-m_0)^\theta
2^{m_0-M_0},
\ee
and we notice that (\ref {m0=M0}) also holds, trivially, if $m_0=M_0$.
Putting together (\ref {m0=M0}), (\ref {-m0}) and (\ref {M-M0}) we get
(\ref {Gamma,M}).

\section{Proof of the recurrence estimates}

\label{sec:proof_rec}

\begin{lemma}
\label{lem61}

For every $\hat \gamma \in \Pi(\RR)$, there exist $M \geq 0$, $\rho<1$
such that for every
$\pi \in \RR$, $q \in \R^\AA_+$,
\be
P_q(\gamma \text { can not be written as }
\gamma_s \hat \gamma \gamma_e
\text { and } \M(B_\gamma \cdot q)>2^M \M(q) \di \pi)
\leq \rho.
\ee

\end{lemma}

\begin{proof}

Fix $M_0 \geq 0$ large and let $M=2 M_0$.
Let $\Gamma$ be the set of all minimal paths $\gamma$ starting at $\pi$
which can not be written as $\gamma_s \hat \gamma \gamma_e$ and such that
$\M(B_\gamma \cdot q)>2^M \M(q)$.  Any path
$\gamma \in \Gamma$ can be written as
$\gamma=\gamma_1 \gamma_2$ where $\gamma_1$ is minimal with $\M(B_\gamma
\cdot q)>2^{M_0} \M(q)$.  Let $\Gamma_1$ collect the possible $\gamma_1$. 
Then $\Gamma_1$ is disjoint.  Let $\tilde \Gamma_1 \subset \Gamma_1$ be the
set of all $\gamma_1$ such that $\M_{\AA'}(B_{\gamma_1} \cdot q) \geq \M(q)$
for all $\AA' \subset \AA$ non-empty.  By Theorem \ref {4},
if $M_0$ is sufficiently large we have
\be
P_q(\Gamma_1 \setminus \tilde \Gamma_1 \di \pi)<\frac {1} {2}.
\ee

For $\pi_e \in \RR$, let $\gamma_{\pi_e}$ be a shortest possible path
starting at $\pi_e$ with $\gamma_{\pi_e}=\gamma_s \hat \gamma$.
If $M_0$ is sufficiently large then
$\norm{B_{\gamma_{\pi_e}}}<
\frac {1} {d} 2^{M_0-1}$.  It follows that if $\gamma_1
\in \Gamma_1$ ends at $\pi_e$ then
\be
P_q(\Gamma \di \gamma_1) \leq 1-P_{B_{\gamma_1} \cdot q}(\gamma_{\pi_e} \di
\pi_e).
\ee
If furthermore $\gamma_1 \in \tilde \Gamma_1$ then
\be
P_{B_{\gamma_1} \cdot q}(\gamma_{\pi_e} \di \pi_e)=\frac {\n(B_{\gamma_1}
\cdot q)} {\n(B_{\gamma_{\pi_e}} \cdot B_{\gamma_1} \cdot q)} \geq
\frac {\M(q)^d} {(2^{2 M_0} \M(q))^d}=2^{-2 d M_0}.
\ee
The result follows with $\rho=1-2^{-2d M_0-1}$.
\end{proof}

\begin{prop}

For every $\hat \gamma \in \Pi(\RR)$,
there exist $\delta>0$, $C>0$ such that for every $\pi \in \RR$, $q \in
\R^\AA_+$ and for every $T>1$
\be
P_q(\gamma \text { can not be written as }
\gamma_s \hat \gamma
\gamma_e \text { and } \M(B_\gamma \cdot q)>T \M(q) \di \pi)
\leq C T^{-\delta}.
\ee

\end{prop}

\begin{proof}

Let $M$ and $\rho$ be as in the previous lemma.  Let $k$ be maximal with $T
\geq 2^{k(M+1)}$.
Let $\Gamma$ be the set of minimal paths
$\gamma$ such that $\gamma$ is not of the form
$\gamma_s \hat \gamma \gamma_e$ and $\M(B_\gamma \cdot q)>2^{k(M+1)} \M(q)$.
Any path $\gamma \in
\Gamma$ can be written as $\gamma_1...\gamma_k$ where
$\gamma_{(i)}=\gamma_1...\gamma_i$ is minimal with $\M(B_{\gamma_{(i)}}
\cdot q)>2^{i (M+1)} \M(q)$.
Let $\Gamma_{(i)}$ collect the $\gamma_{(i)}$.  Then
the $\Gamma_{(i)}$ are disjoint. Moreover, by Lemma \ref{lem61}, 
for all $\gamma^{(i)} \in
\Gamma_{(i)}$,
\be
P_q(\Gamma_{(i+1)} \di \gamma^{(i)}) \leq \rho.
\ee
This implies that $P_q(\Gamma \di \pi) \leq \rho^k$.  The result follows.
\end{proof}

\noindent{\it Proof of Theorem \ref {exponentialtails}.}
Let $\pi$ be the start of $\gamma_*$.  
The push-forward under radial projection
of the Lebesgue measure on $\e_{q_0}$ onto $\Delta_\pi \cap \Up^{(1)}_\RR$
yields a smooth measure $\tilde \nu$.  It is enough to show that $\tilde
\nu\{x \in \De \tq r_\De(x) \geq T\} \leq C T^{-\delta}$, for some $C>0$,
$\delta>0$.
A connected component of the domain of $T_\De$ that intersects the set
$\{x \in \De \tq r_\De(x) \geq T\}$ is of the form
$\Delta_\gamma \cap \Up^{(1)}_\RR$
where $\gamma$ can not be written as $\gamma_s \hat \gamma \gamma_e$ with
$\hat \gamma=\gamma_* \gamma_* \gamma_* \gamma_*$ and $\M(B_\gamma \cdot
q_0) \geq C^{-1} T$, where $q_0=(1,...,1)$ and
$C$ is a constant depending on $\gamma_*$.  Thus
\be
\tilde \nu \{x \in \De \tq r_\De(x) \geq T\} \leq P_{q_0}(\gamma \text { can
not be written as } \gamma_s \hat \gamma \gamma_e \text { and }
\M(B_\gamma \cdot q_0) \geq C^{-1} T \di \pi).
\ee
The result follows from the previous proposition.
\qed

\begin{lemma}

For every $k_0 \geq 1$ there exist $C>0$, $\theta>0$, depending only on
$\# \AA$ and $k_0$ with the following property.
Let $M \in \N$, $q \in \R^\AA_+$.
Then for every $\pi \in \RR$,
\be
P_q(\gamma \text { is not $k_0$-complete and }
\M(B_\gamma \cdot q)>2^M \M(q) \di \pi)
\leq C (M+1)^\theta 2^{-M}.
\ee

\end{lemma}

\begin{proof}

The proof is by induction on $k_0$.  For $k_0=1$,
it is Proposition \ref {2}.
Assume it holds for some $k_0 \geq 1$.  Let $\Gamma$ be the set of minimal
paths which are not $k_0+1$-complete and such that $\M(B_\gamma \cdot q)>2^M
\M(q)$.  Let $\Gamma_- \subset \Gamma$ be the set of paths which are not
$k_0$-complete.  Then
$P_q(\Gamma_- \di \pi) \leq C (M+1)^\theta 2^{-M}$ by the induction hypothesis.
Every $\gamma \in \Gamma \setminus \Gamma_-$ can be written as
$\gamma=\gamma_s\gamma_e$ with $\gamma_s$ minimal $k_0$-complete.  Let
$m=m(\gamma_s) \in [0,M]$
be such that $2^m \M(q)<\M(B_{\gamma_s} \cdot q) \leq 2^{m+1} \M(q)$.  Let
$\Gamma_m$ collect the $\gamma_s$ with $m(\gamma_s)=m$.  Then $\Gamma_m$ is
disjoint.  By the induction hypothesis $P_q(\Gamma_m \di \pi) \leq C 
(m+1)^\theta
2^{-m}$ and by Proposition \ref {2},
$P_q(\Gamma \di \gamma_s) \leq (M+1-m)^\theta 2^{m-M}$, $\gamma_s \in \Gamma_m$.
The result follows by summing over $m$.
\end{proof}

\begin{prop}

For every $k_0 \geq 2 \#\AA-3$,
$\delta>0$, there exist $C>0$ and a finite
disjoint set $\Gamma_0 \subset \Pi(\RR)$ with the following properties:
\begin{enumerate}
\item If $\gamma \in \Gamma_0$ then $\gamma$ is minimal $k_0$-complete,
\item For every $\pi \in \RR$, $q \in \R^\AA_+$, $T \geq 0$,
\be
\label{blah69}
P_q(\gamma \text { can not be written as }
\gamma_s \gamma_0 \gamma_e \text { with }
\gamma_0 \in \Gamma_0 \text { and } \M(B_\gamma \cdot q)>T \M(q) \di \pi)
\leq C T^{(\delta-1)}.
\ee
\end{enumerate}

\end{prop}

\begin{proof}

Fix some $M \geq 0$.  Let $\Gamma_0$ be the set of all minimal
paths which are $k_0$-complete
and such that $\|B_\gamma\| \leq 2^{M+2}$.  Obviously $\Gamma_0$
satisfies condition (1).  Let us show that if $M$ is large then it also
satisfies condition (2). It is sufficient to prove \eqref{blah69} for
times $T$ of the form $2^{k(M+1)}$.

For $k \geq 0$, let $\Gamma$ be the set of paths
$\gamma$ such that $\gamma$ is not of the form
$\gamma_s \gamma_0 \gamma_e$ with $\gamma_0
\in \Gamma_0$ and $\M(B_\gamma \cdot q)>2^{k(M+1)} \M(q)$.
Any path $\gamma \in
\Gamma$ can be written as $\gamma_1...\gamma_k$ where
$\gamma_{(i)}=\gamma_1...\gamma_i$ is minimal with $\M(B_{\gamma_{(i)}}
\cdot q)>2^{i (M+1)} \M(q)$.
Let $\Gamma_{(i)}$ collect the $\gamma_{(i)}$.  Then
the $\Gamma_{(i)}$ are disjoint.

Notice that the $\gamma_i$ are not $2k_0$-complete.  Otherwise,
$\gamma_i=\gamma_s\gamma_e$ with $\gamma_s$ and $\gamma_e$
$k_0$-complete.  By Lemma \ref {complete 2d-3},
all coordinates of
$B_{\gamma_s} \cdot B_{\gamma_{(i-1)}} \cdot q$ are larger than
$\M(B_{\gamma_{(i-1)}} \cdot q)>2^{(i-1) (M+1)}$.  It follows that
$\|B_{\gamma_e}\| \leq 2^{M+2}$, so $\gamma_e \in \Gamma_0$, contradiction.

By the previous lemma,
$P_q(\Gamma_{(i)} \di \gamma_s) \leq C (M+1)^\theta 2^{-M}$,
$\gamma_s \in \Gamma_{(i-1)}$.  This implies that $P_q(\Gamma \di \pi) \leq
(C (M+1)^\theta 2^{-M})^k$.  If $M$ is large enough, this gives
$P_q(\Gamma \di \pi) \leq 2^{(\delta-1) k(M+1)}$.
\end{proof}

\noindent{\it Proof of Theorem \ref {largeballs}.}
Let $\Gamma_0$ be as in the previous proposition, with $k_0=6\#\AA-8$.
We let $\widehat Z=\bigcup \Delta_{\gamma_e} \times \Ga_{\gamma_s}$ where
$\gamma_s$ is minimal $4\#\AA-6$ complete, $\gamma_e$ is minimal
$2 \#\AA-3$-complete and there exists $\gamma \in \Gamma_0$ that starts by
$\gamma_s \gamma_e$.
Its intersection with $\widehat \Up^{(1)}$ is precompact by Lemmas \ref
{complete 2d-3} and \ref {strongly positive}.

Fix some component $\Delta_{\gamma_{e_0}} \times \Ga_{\gamma_{s_0}}$
of $\widehat Z$ and let us estimate $\widehat m\{x \in
\Delta_{\gamma_{e_0}} \times \Ga_{\gamma_{s_0}}
\cap \widehat \Up^{(1)}_\RR \tq r_{\widehat Z}(x)>T\}$.
Let $\pi$ be the start of $\gamma_{s_0}$.  If $\Delta_{\gamma_1} \times
\Ga_{\gamma_2}$ is a component of the domain of the first return map to
$\widehat Z$ that intersects $\{x \in
\Delta_{\gamma_{e_0}} \times \Ga_{\gamma_{s_0}} \tq r_{\widehat Z}(x)>T\}$
then $\gamma_1$ can not be written as
$\gamma_s \gamma_0 \gamma_e$ with $\gamma_0 \in \Gamma_0$.  The projection
of $\hat m|\Delta_{\gamma_{e_0}} \times \Ga_{\gamma_{s_0}} \cap \widehat
\Up^{(1)}_\RR$ on $\Up^{(1)}_\RR$ is absolutely continuous with a bounded
density, so we conclude as in the proof of Theorem \ref {exponentialtails}
that
\begin{align}
\widehat m\{x \in
&\Delta_{\gamma_{e_0}} \times \Ga_{\gamma_{s_0}}
\cap \widehat \Up^{(1)}_\RR \tq r_{\widehat Z}(x)>T\} \leq\\
\nonumber
&C P_{q_0}(\gamma \text { can not be written as }
\gamma_s \gamma_0 \gamma_e \text { with }
\gamma_0 \in \Gamma_0 \text { and } \M(B_\gamma \cdot q)>T \M(q) \di \pi),
\end{align}
where $q_0=(1,...,1)$.  The result follows from the previous proposition.
\qed

\section{Exponential mixing for expanding semiflows}

\label{par:define_expanding}

In this section and the next, our goal is to prove Theorem
\ref{main_thm_hyperbolic}. As a first step, we will prove in this
section an analogous result concerning expanding semi-flows.

Let $T: \bigcup \Delta^{(l)} \to \Delta$ be a uniformly expanding
Markov map on a John domain $(\Delta,\Leb)$, with expansion constant
$\expansion>1$, and let $r: \Delta \to \R_+$ be a good roof function
with exponential tails (as defined in Paragraph
\ref{par:define_exp}). Let $\Delta_r=\{ (x,t) \tq x\in \Delta, 0\leq
t<r(x)\}$, we define a semi-flow $T_t : \Delta_r \to \Delta_r$, by
$T_t(x,s)= (T^n x, s+t-r^{(n)}(x))$ where $n$ is the unique integer
satisfying $r^{(n)}(x)\leq t+s<r^{(n+1)}(x)$. Let $\mu$ be the
absolutely continuous probability measure on $\Delta$ which is
invariant under $T$, then the flow $T_t$ preserves the probability
measure $\mu_r= \mu \otimes \Leb / (\mu\otimes \Leb)(\Delta_r)$. We
will also use the finite measure $\Leb_r=\Leb \otimes \Leb$ on
$\Delta_r$. In this section, we will be interested in the mixing
properties of $T_t$. \emph{Unless otherwise specified, all
the integrals will be
taken with respect to the measures $\Leb$ or $\Leb_r$}.

Let us first define the class of functions for which we can prove
exponential decay of correlations:
\begin{definition}
A function $U: \Delta_r \to \R$ belongs to $\BB_0$ if it is bounded,
continuously differentiable on each set $\Delta^{(l)}_r:=\{ (x,t)
\tq x\in \Delta^{(l)}, 0<t<r(x)\}$, and $\sup_{(x,t) \in \bigcup
\Delta^{(l)}_r} \norm{DU(x,t)}<\infty$. Write then
  \begin{equation}
  \norm{U}_{\BB_0}=\sup_{(x,t)\in \bigcup \Delta_r^{(l)}} |U(x,t)|
  + \sup_{(x,t)\in \bigcup \Delta_r^{(l)}} \norm{DU(x,t)}.
  \end{equation}
\end{definition}
Notice that such a function is not necessarily continuous on the
boundary of $\Delta^{(l)}_r$.
\begin{definition}
A function $U:\Delta_r \to \R$ belongs to $\BB_1$ if it is bounded
and there exists a constant $C>0$ such that, for all fixed $x\in
\bigcup_l \Delta^{(l)}$, the function $t \mapsto U(x,t)$ is of
bounded variation on the interval $(0,r(x))$ and its variation is
bounded by $Cr(x)$. Let
  \begin{equation}
  \norm{U}_{\BB_1}=\sup_{(x,t)\in \bigcup \Delta_r^{(l)}} |U(x,t)| +
  \sup_{x\in \bigcup \Delta^{(l)}}
  \frac{\Var_{(0,r(x))}( t\mapsto U(x,t))}{r(x)}.
  \end{equation}
\end{definition}
This space $\BB_1$ is very well suited for further extensions to the
hyperbolic case. In this paper, the notation $C^1(X)$ for some space
$X$ always denotes the space of bounded continuous functions on $X$
which are everywhere continuously differentiable and such that the
norms of the differentials are bounded. Then the following
inclusions hold:
  \begin{equation}
  C^1 \subset \BB_0 \subset \BB_1.
  \end{equation}

\begin{thm}
\label{main_thm_expanding} There exist constants $C>0$ and
$\delta>0$ such that, for all functions $U\in \BB_0$ and $V\in
\BB_1$, for all $t\geq 0$,
  \begin{equation}
  \left|\int U\cdot V\circ T_t \dLeb_r - \left(\int U \dLeb_r\right)
  \left( \int V \dd\mu_r\right) \right|
  \leq C \norm{U}_{\BB_0} \norm{V}_{\BB_1}  e^{-\delta t}.
  \end{equation}
\end{thm}

\begin{rem}
\label{main_rmk_expanding}
Applying the previous theorem to the function $U(x,t)\cdot
\frac{\dd\mu}{\dLeb}(x)$, we also obtain
  \begin{equation}
  \left|\int U\cdot V\circ T_t \dd\mu_r - \left(\int U \dd\mu_r\right)
  \left( \int V \dd\mu_r\right) \right|
  \leq C \norm{U}_{\BB_0} \norm{V}_{\BB_1}  e^{-\delta t}.
  \end{equation}
\end{rem}

Notation: when dealing with a uniformly expanding Markov map $T$, we
will always denote by $\HH_n$ the set of inverse branches of $T^n$.

The proof of Theorem \ref{main_thm_expanding} will take the rest of
this section.

\subsection{Discussion of the aperiodicity condition}

In this paragraph, we discuss several conditions on the return time
$r$ which turn out to be equivalent to the aperiodicity condition 
(3) in Definition~\ref{def:goodroof}.

\begin{prop}
\label{equiv_UNI} Let $T$ be a uniformly expanding Markov map for a
partition $\{\Delta^{(l)}\}$. Let $r: \Delta \to \R$ be $C^1$ on
each set $\Delta^{(l)}$, with $\sup_{h\in \HH} \norm{D(r\circ
h)}_{C^0}<\infty$. Then the following conditions are equivalent:
\begin{enumerate}
\item There exists $C>0$ such that there exists an arbitrarily large
$n$, there exist $h,k\in \HH_n$, there exists a continuous unitary
vector field $x\mapsto y(x)$ such that, for all $x\in \Delta$,
  \begin{equation}
  \label{eq:C0}
  \bigl| D(r^{(n)}\circ h)(x)\cdot y(x) - D(r^{(n)}\circ k)(x)\cdot y(x) \bigr|
  > C.
  \end{equation}
\item There exists $C>0$ such that there exists an arbitrarily large
$n$, there exist $h,k\in \HH_n$, there exists $x\in \Delta$ and
$y\in T_x\Delta$ with $\norm{y}=1$ such that
  \begin{equation}
  \label{eq:C}
  \bigl| D(r^{(n)}\circ h)(x)\cdot y - D(r^{(n)}\circ k)(x)\cdot y \bigr|
  > C.
  \end{equation}
\item It is not possible to write $r=\psi + \phi\circ T-\phi$ on $\bigcup
\Delta^{(l)}$,
where $\psi:\Delta \to \R$ is constant on each set $\Delta^{(l)}$
and $\phi\in C^1(\Delta)$.
\item It is not possible to write  $r=\psi + \phi\circ T-\phi$ almost
everywhere, where
$\psi:\Delta \to \R$ is constant on each set $\Delta^{(l)}$ and
$\phi : \Delta \to \R$ is measurable.
\end{enumerate}
\end{prop}

The first condition is the (UNI) condition as given in \cite{BV} in
their one-dimensional setting.

\begin{proof}
The implication $(1)\Rightarrow (2)$ is trivial. Let us prove
$(2)\Rightarrow (1)$. Notice that there exists a constant $c_0$ such
that, for any inverse branch $\ell \in \HH_p$ of any iterate $T^p$
of $T$, for any $x\in \Delta$ and any $y\in T_x\Delta$, $|D(r^{(p)}
\circ \ell)(x)\cdot y|\leq c_0 \norm{y}$: for instance, take
$c_0=\frac{\sup_{h\in \HH} \norm{D(r\circ
h)}_{C^0}}{1-\expansion^{-1}}$.

Let $C>0$ be such that \eqref{eq:C} is satisfied for infinitely many
$n$. It is then possible to choose $n$ large enough so that $c_0
\expansion^{-n}\leq C/4$, $h,k\in \HH_n$, $x_0\in \Delta$ and $y_0\in
T_x\Delta$ such that \eqref{eq:C} holds. Let $y_0(x)$ be a unitary
vector field on a neighborhood $U$ of $x$ such that \eqref{eq:C}
still holds for $y_0(x)$. Fix a branch $l\in \HH_m$ for some $m$
such that $l(\Delta) \subset U$. Define a vector field $y_1$ on
$\Delta$ by $y_1(x)=Dl(x)^{-1}y_0(lx)$. For any inverse branch $\ell
\in \HH_p$ for some $p\geq 1$, we have
  \begin{align*}
  \left|D( r^{(m+n+p)} \circ \ell \circ h\circ l)(x)\cdot y_1(x) - D(r^{(m+n)} \circ
  h\circ l)(x)\cdot y_1(x)\right|&
  =\left|D( r^{(p)}\circ \ell)( hlx)Dh(lx)\cdot y_0(lx)\right|
  \\&
  \leq c_0 \norm{Dh(lx)} \leq c_0 \expansion^{-n} \leq C/4.
  \end{align*}
The same estimate applies to $k$. Since
  \begin{multline*}
  |D(r^{(m+n)} \circ h\circ l)(x)\cdot y_1(x) - D(r^{(m+n)} \circ k\circ
  l)(x)\cdot y_1(x)|\\
  =|D(r^{(n)}\circ h)(lx)\cdot y_0(lx) -D(r^{(n)}\circ k)(lx)\cdot y_0(lx)|
  \geq C,
  \end{multline*}
we get
  \begin{equation}
  \left|D( r^{(m+n+p)} \circ \ell \circ h\circ l)(x)\cdot y_1(x)-
  D( r^{(m+n+p)} \circ \ell \circ k\circ l)(x)\cdot y_1(x)\right| \geq
  C/2.
  \end{equation}
Finally, take $y(x)=y_1(x)/\norm{y_1(x)}$. This proves $(1)$.

The implication $(2)\Rightarrow (3)$ is easy: if it is possible to
write $r=\psi+\phi\circ T-\phi$, then for all $h\in \HH_n$,
$r^{(n)}\circ h(x)=S_n \psi(h(x)) + \phi(x) - \phi(hx)$. Hence, if
$\norm{y}=1$,
  \begin{align*}
  \bigl| D(r^{(n)}\circ h)(x)\cdot y - D(r^{(n)}\circ k)(x)\cdot y\bigr|
  &
  = \bigl| D( \phi\circ h)(x)\cdot y - D(\phi \circ k)(x)\cdot y \bigr|
  \\ &
  \leq 2 \norm{\phi}_{C^1} \expansion^{-n}.
  \end{align*}
This quantity tends to $0$ when $n\to \infty$, which is not
compatible with $(2)$.

\newcommand{\h}{\underline{h}}

Let us prove $(3)\Rightarrow (2)$. Assume that (2) does not hold, we
will prove that $r$ can be written as $\psi+\phi\circ T-\phi$. Let
$\h=(h_1,h_2,\dots)$ be a sequence of $\HH$. Write
$\h_n=h_n\circ \dots \circ h_1$. Then
  \begin{equation}
  D(r^{(n)}\circ \h_n)(x)\cdot y= \sum_{k=1}^n D(r\circ h_k)( \h_{k-1} x) D
  \h_{k-1}(x)\cdot y.
  \end{equation}
Since the derivative of $r\circ h_k$ is uniformly bounded
 by assumption and
$\norm{D \h_{k-1}(x)} \leq \expansion^{-k+1}$, this series is uniformly
converging. Since $(2)$ is not satisfied, its limit is independent
of the sequence of inverse branches $\h$, and defines a continuous
$1$-form $\omega(x)\cdot y$ on $\Delta$. It satisfies, for all $h\in
\HH$,
  \begin{equation}
  \label{eq:omega_cobord}
  \omega(x)\cdot y=D(r\circ h)(x)\cdot y+ \omega(hx)Dh(x)\cdot y.
  \end{equation}

Take a branch $h\in \HH$, and let $\h=(h,h,\dots)$. Let $x_0\in
\Delta$.
The series of functions
$\sum_{k=1}^\infty (r\circ \h_k - r\circ \h_k(x_0))$
is then summable in $C^1(\Delta)$,
let us denote its sum by $\phi$. By construction, $\omega(x)\cdot
y=D\phi(x)\cdot y$ for all $x\in \Delta$ and $y\in T_x\Delta$. By
\eqref{eq:omega_cobord}, $D(r+\phi-\phi\circ T)=0$. Hence,
$r+\phi-\phi\circ T$ is constant on each $\Delta^{(l)}$, which
concludes the proof.

\renewcommand{\h}{\overline{h}}

The implication $(4) \Rightarrow (3)$ is trivial, we just have to
prove $(3)\Rightarrow (4)$ to conclude. Assume that $r=\psi +
\phi\circ T- \phi$ where $\psi$ is constant on each set
$\Delta^{(l)}$ and $\phi$ is measurable. We will prove that $\phi$
has a version which is $C^1$. Let $\FF_n$ be the $\sigma$-algebra
generated by the sets $h(\Delta)$ for $h\in \HH_n$. It is an
increasing sequence of $\sigma$-algebras. For almost all $x\in
\Delta$, there exists a well defined sequence $\h=(h_1,h_2,\dots)\in
\HH^\N$ such that the element $F_n(x)$ of $\FF_n$ containing $x$ is
given by $F_n(x)=h_1 \circ \dots \circ h_n(\Delta)$. Equivalently,
$h_n$ is the unique element of $\HH$ such that $T^{n-1}(x)\in h_n(\Delta)$.
Since $T$ is ergodic, almost every
$x$ is normal in the sense that, for any finite sequence $k_1,\dots,
k_p$ of elements of $\HH$, there exist infinitely many $n$ such
that, for all $1 \leq i\leq p$, $h_{n+i}=k_i$.

The martingale convergence theorem shows that, for almost all $x\in
\Delta$, for all $\epsilon>0$,
  \begin{equation}
  \label{martingale}
  \frac{\Leb\{ x' \in F_n(x) \tq
  |\phi(x')-\phi(x)|>\epsilon\}}{\Leb(F_n(x))} \to 0.
  \end{equation}
Take a point $x_0$ such that this convergence holds and which is
normal. Replacing $\phi$ by $\phi-\phi(x_0)$, we can assume that
$\phi(x_0)=0$. Let $\h=(h_1,h_2,\dots)$ be the corresponding
sequence of $\HH$ and write $\h_n=h_1\circ \dots \circ h_n$, so that
$F_{n}(x_0)=\h_n(\Delta)$. Then \eqref{martingale} and distortion
controls give, for all $\epsilon>0$,
  \begin{equation}
  \label{tend0}
  \Leb \{ x\in \Delta \tq |\phi(\h_n x)|>\epsilon\} \to 0.
  \end{equation}
Define a strictly increasing sequence $m_k$ as follows: start from
$m_1=1$. If $m_k$ has been defined then, by normality of $x_0$,
there exists $m_{k+1}>m_k$ such that $(h_1,\dots, h_{m_{k+1}})$
finishes with $(h_1,\dots, h_{m_k})$. By \eqref{tend0}, we can
choose a subsequence $n_k$ of $m_k$ such that
  \begin{equation}
  \forall \epsilon>0,\
  \sum_{k=1}^\infty \Leb \{ x\in \Delta \tq |\phi(\h_{n_k} x)|>\epsilon\}
  <\infty.
  \end{equation}
In particular, for almost all $x$, $\phi(\h_{n_k} x) \to 0$. Notice
that $\phi(x)= \phi(\h_n x)+r^{(n)}(\h_n x)-S_n \psi(\h_n x)$. For
almost all $x$, we get $\phi(x)= \lim_{k\to \infty}
r^{(n_k)}(\h_{n_k} x)-S_{n_k} \psi(\h_{n_k} x)$. Moreover, the
choice of $m_k$ ensures that the sequence $D( r^{(n_k)}\circ
\h_{n_k})$ is Cauchy. Hence, $\phi$ coincides almost everywhere with
the $C^1$ function $\lim _{k\to \infty} r^{(n_k)}\circ \h_{n_k}
-S_{n_k} \psi\circ \h_{n_k}$, which concludes the proof.
\end{proof}

\subsection{Existence of bump functions}

The following technical lemma will prove useful later.
\begin{lem}
\label{lem_exists_bump}
There exist constants $\Cd>1$ and $\Ce>0$
satisfying the following property: for any ball $B(x,r)$ compactly
included in $\Delta$, there exists a $C^1$ function $\rho : \Delta
\to [0,1]$ such that $\rho=0$ on $\Delta \moins B(x,r)$, $\rho=1$
on $B(x, r/\Cd)$ and $\norm{\rho}_{C^1} \leq \Ce/r$.
\end{lem}
Notice that this property is not true for any John domain, and
uses the existence of the uniformly expanding Markov map $T$ on
$\Delta$.
\begin{proof}
Let $x_0\in \Delta$ be in the domain of definition of all iterates
of $T$. Let $\norm{\cdot}'$ be a flat Riemannian metric on a
neighborhood of $x$. By compactness, there exists a constant $K>0$
such that, on a small neighborhood $U$ of $x_0$,
$K^{-1}\norm{\cdot}' \leq \norm{\cdot} \leq K \norm{\cdot}'$.

For large enough $n$, the inverse branch $h\in \HH_n$ such that $x_0
\in h(\Delta)$ satisfies $h(\Delta) \subset U$, since $\diam
(h(\Delta)) \leq C \expansion^{-n}$. The set $h(\Delta)$ endowed
with the distance given by $\norm{\cdot}'$ is flat. Hence, there exists a
constant $C>0$ such that, given any ball $B'=B'(x,r)$ for this
Euclidean distance,
which is compactly included in $h(\Delta)$, there exists a $C^1$ function
$\rho$ supported in $B'$, equal to $1$ on $B'(x,r/2)$ and with
$\norm{\rho}_{C^1} \leq C/r$.

Since $h$ and its inverse have
uniformly bounded derivatives (with respect to $\norm{\cdot}$ and
$\norm{\cdot}'$), this easily implies the lemma.
\end{proof}

The same compactness argument also implies the following lemma:
\begin{lem}
\label{lem:precompact}
For all $\epsilon>0$,
  \begin{equation*}
  \sup\{ k \in \N \tq \exists x_1,\dots,x_k \in \Delta \text{ with
  }d(x_i,x_j) \geq \epsilon \text{ whenever }i\not=j\} < \infty.
  \end{equation*}
\end{lem}
\subsection{A Dolgopyat-like spectral estimate}

The main step of the proof of Theorem \ref{main_thm_expanding} is
the study of the spectral properties of weighted transfer operators
$L_s$. Let $\sigma_0>0$ be such that $\int e^{\sigma_0
r}\dLeb<\infty$, which is possible since $r$ has exponential tails.
For $s\in \C$ with $\Re s>-\sigma_0$, define
  \begin{equation}
  L_s u(x)=\sum_{Ty=x} e^{-sr(y)} J(y)u(y).
  \end{equation}

For $s=\sigma+it$ with $\Re s>-\sigma_0$ and $t\in \R$, define a
norm on $C^1(\Delta,\C)$ by
  \begin{equation}
  \norm{u}_{1,t}= \sup_{x\in \Delta} |u(x)| + \frac{1}{\max(1,|t|)} \sup_{x\in
  \Delta} \norm{Du(x)}.
  \end{equation}

The main spectral estimate concerning the operators $L_s$ is the
following Dolgopyat-like estimate:
\begin{prop}
\label{thm_UNI} There exist $\sigma'_0 \leq \sigma_0$, $T_0>0$, $C>0$
and $\beta<1$ such that, for all $s=\sigma+it$ with $|\sigma|\leq
\sigma'_0$ and $|t|\geq T_0$, for all $u\in C^1(\Delta)$, for all $k
\in \N$,
  \begin{equation}
  \norm{L_s^k u}_{L^2} \leq C \beta^k \norm{u}_{1,t}.
  \end{equation}
\end{prop}

This paragraph will be entirely devoted to the proof of Proposition
\ref{thm_UNI}. The proof will follow very closely the arguments in
\cite{BV}, with small complications due to the general dimension.

For $s=0$, $L_s$ is the usual transfer operator. It acts on the
space of $C^1$ functions, has a spectral gap, and a simple
isolated eigenvalue at $1$ (the corresponding eigenfunction will
be denoted by $f_0$ and is the density of the invariant measure
$\mu$). For $\sigma\in \R$ close enough to $0$, $L_\sigma$ acting
on $C^1(\Delta)$ is a continuous perturbation of $L_0$, by a
straightforward computation. Hence, it has a unique eigenvalue
$\lambda_\sigma$ close to $1$, and the corresponding eigenfunction
$f_\sigma$ (normalized so that $\int f_\sigma=1$) is $C^1$,
strictly positive, and tends to $f_0$ in the $C^1$ topology when
$\sigma \to 0$.

Let $0<\sigma_1\leq \min(\sigma_0,1)$ be such that $f_\sigma$ is
well defined and positive for $\sigma \in [-\sigma_1,\sigma_1]$.
For $s=\sigma+it$ with $|\sigma|\leq \sigma_1$ and $t\in \R$, define
a modified transfer operator $\tilde L_s$ by
  \begin{equation}
  \tilde L_s(u) =\frac{ L_s( f_\sigma u)}{\lambda_\sigma f_\sigma}.
  \end{equation}
It satisfies $\tilde L_\sigma 1=1$, and $|\tilde L_s u| \leq \tilde
L_\sigma |u|$.

\begin{lem}
\label{define_Ca} There exists a constant $\Ca$ such that $\forall
n\geq 1$, $\forall s=\sigma+it$ with $\sigma \in
[-\sigma_1,\sigma_1]$ and $t\in \R$, $\forall u \in C^1(\Delta)$,
holds for all $x\in \Delta$
  \begin{equation}
  \norm{ D( \tilde L_s^n u)(x)} \leq \Ca( |t|+1) \tilde L_\sigma^n
  (|u|)(x) + \expansion^{-n} \tilde L_\sigma^n( \norm{Du}) (x).
  \end{equation}
\end{lem}
\begin{proof}
We have
  \begin{equation}
  \tilde L_s^n u(x)=\sum_{h\in \HH_n} \frac{ (f_\sigma u)(h x)
  J^{(n)}(hx) e^{-s r^{(n)}(hx)}}{\lambda_\sigma^n f_\sigma(x)},
  \end{equation}
where $r^{(n)}(x)=\sum_{k=0}^{n-1} r(T^k x)$ and
$J^{(n)}(x)=\prod_{k=0}^{n-1} J(T^k x)$. Differentiating this
expression, we obtain a sum of $5$ terms: we can differentiate
$f_\sigma$, or $u$, or $J^{(n)}$, or $r^{(n)}$, or $1/f_\sigma$.

Since $f_\sigma$ is bounded in $C^1$ and uniformly bounded from
below, and any inverse branch of $T$ is contracting, there exists a
constant $C>0$ such that $\norm{D (f_\sigma \circ h)(x)} \leq C
f_\sigma(x)$. Hence, if we differentiate $f_\sigma$, the resulting
term is bounded by $C \tilde L_\sigma^n(|u|)(x)$.

In the same way, distortion controls give $\norm{D (J^{(n)} \circ
h)(x)} \leq C J^{(n)}\circ h(x)$. We also have
$\norm{D(1/f_\sigma)(x)}\leq C /f_\sigma(x)$. Hence, the
corresponding terms are also bounded by $C \tilde
L_\sigma^n(|u|)(x)$.

Moreover, $D(e^{-s r^{(n)}\circ h})(x)=-s D(r^{(n)}\circ h)(x) e^{-s
r^{(n)}\circ h(x)}$. The uniform contraction of $h$ and the
boundedness of the derivative of $r\circ \ell$ for $\ell\in \HH$
show that this term is bounded by $C|s| e^{-\sigma r^{(n)}\circ
h(x)}$. Hence, the resulting term is bounded by $C (|t|+1)\tilde
L_\sigma^n (|u|)(x)$.

Finally, $\norm{D(u\circ h)(x)}\leq \expansion^{-n} \norm{Du (hx)}$,
which shows the required bound on the last term.
\end{proof}

From this point on, we will fix once and for all a constant $\Ca>5$
satisfying the conclusion of Lemma \ref{define_Ca}. This lemma
implies that the iterates of $\tilde L_s$ are bounded for the norm
$\norm{\ }_{1,t}$. More precisely, the following holds:
\begin{lem}
\label{acts_boundedly} There exists a constant $C>1$ such that, for
all $s=\sigma+it$ with $\sigma \in [-\sigma_1,\sigma_1]$ and $|t|
\geq 10$, for all $k\in \N$, for all $u \in C^1(\Delta)$,
  \begin{equation}
  \norm{\tilde L_s^k u}_{1,t} \leq C \norm{u}_{C^0} +
  \frac{\expansion^{-k}}{|t|} \norm{ Du}_{C^0}.
  \end{equation}
In particular, $\norm{ \tilde L_s^k u}_{1,t} \leq C \norm{u}_{1,t}$.
\end{lem}
\begin{proof}
The inequality $\norm{\tilde L_s^k u}_{C^0} \leq \norm{u}_{C^0}$ and
Lemma \ref{define_Ca} give
  \begin{equation*}
  \norm{\tilde L_s^k u}_{C^0} +\frac{\norm{D(\tilde L_s^k
  u)}_{C^0}}{|t|}
  \leq \norm{u}_{C^0} + \frac{1}{|t|} \left[ 2\Ca |t| \norm{u}_{C^0}+
  \expansion^{-k} \norm{Du}_{C^0}\right]
  \leq C \norm{u}_{C^0} +
  \frac{\expansion^{-k}}{|t|} \norm{ Du}_{C^0}.
  \qedhere
  \end{equation*}
\end{proof}

To prove Proposition \ref{thm_UNI}, we need to get some
contraction. This is easy to do if the derivative is large compared to
the $C^0$ norm of the function:

\begin{lem}
\label{decays_easily} There exists $N_0 \in \N$ such that any $n\geq
N_0$ satisfies the following property. Let $s=\sigma+it$ with
$\sigma \in [-\sigma_1,\sigma_1]$ and $|t|\geq 10$. Let $v\in
C^1(\Delta)$ satisfy $\sup \norm{Dv} \geq 2\Ca |t| \sup |v|$. Then
  \begin{equation}
  \norm{\tilde L_s^n v}_{1,t} \leq \frac{9}{10} \norm{v}_{1,t}.
  \end{equation}
\end{lem}
\begin{proof}
We have
  \begin{equation}
  \norm{\tilde L_s^n v}_{C^0} \leq \norm{v}_{C^0} \leq
  \frac{1}{2\Ca |t|} \sup \norm{Dv(x)} \leq \frac{1}{2\Ca} \norm{v}_{1,t}.
  \end{equation}
Moreover, for $x\in \Delta$
  \begin{align*}
  \norm{D(\tilde L_s^n v)(x)} &\leq \Ca (1+|t|) \tilde
  L_\sigma^n(|v|)(x)+\expansion^{-n} \tilde L_\sigma^n(\norm{Dv})(x)
  \leq \Ca (1+|t|) \norm{v}_{C^0} + \expansion^{-n} \norm{Dv}_{C^0}
  \\ &
  \leq \left[ \frac{1+|t|}{2} + \expansion^{-n} |t|\right]
  \norm{v}_{1,t}.
  \end{align*}
Hence,
  \begin{equation}
  \norm{\tilde L_s^n v}_{C^0} + \frac{1}{|t|} \norm{D(\tilde L_s^n
  v)}_{C^0}
  \leq \left[ \frac{1}{2\Ca} + \frac{1+|t|}{2|t|} +
  \expansion^{-n}\right] \norm{v}_{1,t}.
  \end{equation}
Since $\Ca \geq 5$ and $|t|\geq 10$, the conclusion of the lemma
holds as soon as $\expansion^{-n}\leq \frac{1}{5}$.
\end{proof}

Hence, to prove Proposition \ref{thm_UNI}, we will mainly have to deal
with functions $v$ satisfying $\sup \norm{Dv} \leq 2 \Ca |t| \sup
|v|$. For technical reasons, it is more
convenient to introduce the following notation.

\begin{definition}
For $t\in \R$, we will say that a pair $(u,v)$ of functions on
$\Delta$ belongs to $\EE_t$ if $u:\Delta\to \R_+$ is $C^1$,
$v:\Delta \to \C$ is $C^1$, $0\leq |v|\leq u$ and
  \begin{equation}
  \forall x\in \Delta,\quad  \max( \norm{Du(x)}, \norm{Dv(x)}) \leq 2\Ca |t|
  u(x).
  \end{equation}
\end{definition}

\begin{lem}
\label{lem:controls_with_chi} There exists $N_1\in \N$ such that any
$n\geq N_1$ satisfies the following property. Let $s=\sigma+it$ with
$\sigma\in [-\sigma_1,\sigma_1]$ and $|t|\geq 10$. Let $(u,v)\in
\EE_t$. Let $\chi \in C^1(\Delta)$ with $\norm{D\chi}\leq |t|$ and
$3/4\leq \chi \leq 1$. Assume that
  \begin{equation}
  \forall x\in \Delta,\quad
  |\tilde L_s^n v(x)|\leq \tilde L_\sigma^n(\chi u)(x).
  \end{equation}
Then $(\tilde L_\sigma^n(\chi u), \tilde L_s^n(v)) \in \EE_t$.
\end{lem}
\begin{proof}
Let $(u,v)\in \EE_t$ with $|t|\geq 10$. Let $n\in \N$. By Lemma
\ref{define_Ca}, for $x\in \Delta$,
  \begin{equation}
  \norm{D(\tilde L_\sigma^n(\chi u))(x)}
  \leq \Ca \tilde L_\sigma^n(\chi u)(x)+ \expansion^{-n} \tilde
  L_\sigma^n( \norm{D(\chi u)})(x).
  \end{equation}
Since $(u,v)\in \EE_t$ and $\norm{D\chi}\leq |t|$,
  \begin{equation}
  \norm{D(\chi u)(x)} \leq |t| u(x)+\norm{Du(x)}\leq (1+2\Ca) |t| u(x)\leq
  \frac{4}{3}(1+2\Ca)|t| (\chi u)(x) .
  \end{equation}
Hence,
  \begin{equation}
  \norm{D(\tilde L_\sigma^n(\chi u))(x)}
  \leq \left[\Ca  + \expansion^{-n}
  \frac{4}{3}(1+2\Ca)|t|\right] \tilde L_\sigma^n(\chi u)(x).
  \end{equation}
If $n$ is large enough, the factor is $\leq 2\Ca |t|$, and we get
$\norm{D(\tilde L_\sigma^n (\chi u))(x)} \leq 2\Ca |t| \tilde
L_\sigma^n(\chi u)(x)$. This is half of what we have to prove.

Concerning $v$, Lemma \ref{define_Ca} gives
  \begin{align*}
  \norm{D( \tilde L_s^n v)(x)}&
  \leq \Ca ( 1+|t|) \tilde L_\sigma^n(|v|)(x) + \expansion^{-n}
  \tilde L_\sigma^n( \norm{Dv})(x)
  \\&
  \leq \Ca (1+|t|) \frac{4}{3}\tilde L_\sigma^n (\chi u)(x)
  +\expansion^{-n} \frac{4}{3} 2\Ca |t| \tilde L_\sigma^n(\chi u)(x).
  \end{align*}
If $n$ is large enough, this quantity is again bounded by $2\Ca |t|
\tilde L_\sigma^n(\chi u)(x)$.
\end{proof}

If $h \in \HH_n$, then $\norm{Dh(x)\cdot y} \leq \expansion^{-n}
\norm{y}$. In particular, since $r$ satisfies Condition (3) of
Definition \ref{def:goodroof}, the first
condition of Proposition \ref{equiv_UNI} gives $n\geq
\max(N_0,N_1)$, two inverse branches $h,k\in \HH_n$ and a continuous
unitary vector field $y_0$ on $\Delta$ such that, for all $x\in
\Delta$,
  \begin{equation}
  \Bigl| D( r^{(n)}\circ h)(x)\cdot y_0(x) -D (r^{(n)}\circ k)(x)\cdot y_0(x) \Bigr|
  \geq
  10 \Ca \max( \norm{ Dh(x)\cdot y_0(x)}, \norm{Dk(x)\cdot y_0(x)}).
  \end{equation}
Smoothing the vector field $y_0$, we get a smooth vector field $y$
with $1\leq \norm{y} \leq 2$ such that, for all $x\in \Delta$,
  \begin{equation}
  \label{eq_UNI}
  \Bigl| D( r^{(n)}\circ h)(x)\cdot y(x) -D (r^{(n)}\circ k)(x)\cdot
  y(x) \Bigr|
  \geq
  9 \Ca \max( \norm{ Dh(x)\cdot y(x)}, \norm{Dk(x)\cdot y(x)}).
  \end{equation}
We fix $n,h,k$ and $y$ as above, until the end of the proof of
Proposition
\ref{thm_UNI}.

\begin{lem}
\label{lem:main_UNI} There exist $\delta>0$ and $\zeta>0$ satisfying
the following property. Let $s=\sigma+it$ with $\sigma\in
[-\sigma_1,\sigma_1]$ and $|t|\geq 10$. Let $(u,v)\in \EE_t$. For
all $x_0\in \Delta$ such that the ball $B(x_0, (\zeta+\delta)/|t|)$
is compactly included in $\Delta$, there exists a point $x_1$ with
$d(x_0,x_1)\leq \zeta/|t|$ such that one of the following
possibilities holds:
\begin{itemize}
\item Either, for all $x\in B(x_1, \delta/|t|)$,
  \begin{multline*}
  \bigl| e^{-s r^{(n)}\circ h(x)} J(hx) (v\cdot f_\sigma)(hx)
  + e^{-s r^{(n)}\circ k(x)} J(kx) (v\cdot f_\sigma)(kx) \bigr|
  \\
  \leq
  \frac{3}{4} e^{-\sigma r^{(n)}\circ h(x)} J(hx) (u\cdot f_\sigma)(hx) +
  e^{-\sigma r^{(n)}\circ k(x)} J(kx) (u\cdot f_\sigma)(kx).
  \end{multline*}
\item Or, for all $x\in B(x_1, \delta/|t|)$,
   \begin{multline*}
  \bigl| e^{-s r^{(n)}\circ h(x)} J(hx) (v\cdot f_\sigma)(hx)
  + e^{-s r^{(n)}\circ k(x)} J(kx) (v\cdot f_\sigma)(kx) \bigr|
  \\ \leq
  e^{-\sigma r^{(n)}\circ h(x)} J(hx) (u\cdot f_\sigma)(hx) +
  \frac{3}{4} e^{-\sigma r^{(n)}\circ k(x)} J(kx) (u\cdot f_\sigma)(kx).
  \end{multline*}
\end{itemize}
\end{lem}
\begin{proof}

Take some constants $\delta>0$ and $\zeta>0$. Let $t\in \R$ with
$|t|\geq 10$. Take $(u,v)\in \EE_t$. Consider $x_0 \in \Delta$ such
that the ball $B(x_0, (\zeta+\delta)/|t|)$ is compactly included in
$\Delta$. If $\delta$ is small enough and $\zeta$ is large enough,
we will find a point $x_1 \in B(x_0,\zeta/|t|)$ for which the
conclusion of the lemma holds.

{\it First case:} Assume first that there exists $x_1 \in
B(x_0,\zeta/|t|)$ such that $|v \circ h(x_1)| \leq u\circ h(x_1)/2$
or $|v \circ k(x_1)| \leq u\circ k(x_1)/2$. We will show that this
point satisfies the required conclusion. The situation being
symmetric, we can assume that $|v \circ h(x_1)| \leq u\circ
h(x_1)/2$.

Since $(u,v)\in \EE_t$, we have $\norm{Du(x)}\leq 2\Ca |t| u(x)$.
This implies $\norm{D(u\circ h)(x)} \leq 2\Ca|t|u\circ h(x)$ since
$h$ is a contraction. We can integrate this inequality along an
almost length-minimizing path between two points $x,x'$:
Gronwall's inequality gives $u(hx') \leq e^{2\Ca |t| d(x,x')}
u(hx)$.

For $x \in B(x_1, \delta/|t|)$, we get
  \begin{equation}
  \norm{ D(v\circ h)(x)} \leq 2\Ca|t| u(hx) \leq 2\Ca |t| e^{2\Ca |t|
  \delta/|t|} u(hx_1).
  \end{equation}
Hence,
  \begin{equation}
  |v(hx) -v(hx_1)| \leq 2\Ca |t| e^{2\Ca \delta} u(hx_1) \delta/|t|.
  \end{equation}
Since $|v(hx_1)|\leq u(hx_1)/2$, we get
  \begin{equation}
  |v(hx)| \leq \left( \frac{1}{2}+ 2\Ca \delta e^{2\Ca \delta}\right)
  u(hx_1)
  \leq
  \left( \frac{1}{2}+ 2\Ca \delta e^{2\Ca \delta}\right) e^{2\Ca
  \delta} u(hx).
  \end{equation}
If $\delta$ is small enough, we get $|v(hx)| \leq \frac{3}{4}u(hx)$
for all $x\in B(x_1, \delta/|t|)$. This concludes the proof.

{\it Second case:} Assume that, for all $x \in B(x_0,\zeta/|t|)$,
holds $|v\circ h(x)| > u\circ h(x)/2$ and $|v\circ k(x)|>
u\circ k(x)/2$.

Let $\phi : [0, \zeta/(2|t|)] \to \Delta$ be the solution of the
equation $\phi'(\tau)= y(\phi(\tau))$ with $\phi(0)=x_0$. Write
$x^\tau=\phi(\tau)$. We will first show that there exists $\tau\leq
\zeta/(8|t|)$ for which $F(x^\tau):= e^{-s r^{(n)}\circ h(x^\tau)}
J\circ h(x^\tau) (v\cdot f_\sigma)(hx^\tau)$ and $G(x^\tau):=e^{-s
r^{(n)}\circ k(x^\tau)} J\circ k(x^\tau) (v\cdot f_\sigma)(kx^\tau)$
have opposite phases. Let $\gamma(\tau)$ be the difference of their
phases.

On the set $h( B(x_0,\zeta/|t|))\cup k( B(x_0, \zeta/|t|))$, the
function $v$ is non vanishing. Hence, it can locally be written as
$v(x)=\rho(x)e^{i \theta(x)}$. Since $Dv(x)=D\rho(x)
e^{i\theta(x)} + i\rho(x) e^{i\theta(x)} D\theta(x)$, the
inequality $\norm{Dv(x)} \leq 2\Ca |t| u(x)$ yields
  \begin{equation}
  \label{eq_der_theta}
  \norm{D\theta(x)} \leq 2\Ca |t| u(x) / \rho(x) \leq 4 \Ca |t|.
  \end{equation}
Since $\gamma(\tau)=- t r^{(n)}(h x^\tau) + \theta(h x^\tau) + t
r^{(n)}(k x^\tau) - \theta(k x^\tau)$, we get
  \begin{multline*}
  \gamma'(\tau) = t \left[
  D(r^{(n)}\circ k)(x^\tau)\cdot y(x^\tau)
  -D(r^{(n)}\circ h) (x^\tau)\cdot y(x^\tau)
  \right]
  \\
  + D \theta( h x^\tau)
  Dh(x^\tau)\cdot y(x^\tau) - D\theta(k x^\tau) Dk(x^\tau)\cdot y(x^\tau).
  \end{multline*}
By \eqref{eq_UNI} and \eqref{eq_der_theta}, we get
  \begin{align*}
  |\gamma'(\tau)| & \geq 9\Ca |t| \max( \norm{ Dh(x^\tau)\cdot y(x^\tau)},
  \norm{Dk(x^\tau)\cdot y(x^\tau)})
  \\ & \hphantom{je tape des trucs longs pour faire beau}
  - 4\Ca |t| \norm{ Dh(x^\tau)\cdot y(x^\tau)} - 4\Ca |t|
  \norm{Dk(x^\tau)\cdot y(x^\tau)}
  \\&
  \geq \Ca |t| \max( \norm{ Dh(x^\tau)\cdot y(x^\tau)},
  \norm{Dk(x^\tau)\cdot y(x^\tau)}).
  \end{align*}
There exists a constant $\gamma_0>0$ such that, for all $x\in
\Delta$ and all $y\in T_x \Delta$ with $1\leq \norm{y}\leq 2$,
$\norm{Dh(x)\cdot y} \geq \gamma_0$ and $\norm{Dk(x)\cdot y} \geq
\gamma_0$. We get finally
  \begin{equation}
  |\gamma'(\tau)| \geq |t| \Ca \gamma_0.
  \end{equation}
If $\zeta=16\pi/(\Ca \gamma_0)$, we obtain $\tau \in
[0,\zeta/(8|t|)]$ for which $F(x^\tau)$ and $G(x^\tau)$ have
opposite phases. Set $x_1=x^\tau \in B(x_0, \zeta/(4|t|))$.

From the definition of $F$ and the inequality $\norm{D(v\circ
h)(x)}\leq  4 \Ca |t| |v(hx)|$ on the ball $B(x_0,\zeta/|t|)$, it
is easy to check the existence of a constant $C$ independent of
$\delta$ such that, for all $x\in B(x_0,\zeta/|t|)$, $\norm{DF(x)}
\leq C |t| |F(x)|$. If $x,x'\in B(x_0, \zeta/(3|t|))$, an almost
length-minimizing path $\gamma$ between $x$ and $x'$ is contained
in $B(x_0, \zeta/|t|)$. Gronwall's inequality along this path
yields $|F(x')| \leq e^{C |t| d(x,x')} |F(x)|$. Moreover, if
$\Gamma_F$ denotes the phase of $F(x)$, we have
$\norm{D\Gamma_F(x)} \leq C |t|$. On the ball $B(x_1,\delta/|t|)$
(which is included in $B(x_0,\zeta/(3|t|))$ as soon as $\delta\leq
\zeta/12$), we get:
  \begin{equation}
  |\Gamma_F(x)-\Gamma_F(x_1)|\leq C \delta \text{ and } e^{-\delta
  C}\leq
  \frac{|F(x)|}{|F(x_1)|} \leq e^{\delta C}.
  \end{equation}
In the same way, if $\Gamma_G$ denotes the phase of $G$, we have for
all $x\in B(x_1,\delta/|t|)$
  \begin{equation}
  |\Gamma_G(x)-\Gamma_G(x_1)|\leq C \delta \text{ and } e^{-\delta
  C}\leq
  \frac{|G(x)|}{|G(x_1)|} \leq e^{\delta C}.
  \end{equation}
Assume for example that $|F(x_1)|\geq |G(x_1)|$ (the other case is
symmetric). If $\delta$ is small enough, we get for all $x\in
B(x_1,\delta/|t|)$
  \begin{equation}
  \label{eq:phase_and_modulus}
  |\Gamma_F(x) - \Gamma_G(x) - \pi| \leq \pi/6 \text{ and }
  |F(x)|\geq |G(x)|/2.
  \end{equation}
We can then use the following elementary lemma:
\begin{lem}
Let $z=re^{i\theta}$ and $z'=r'e^{i\theta'}$ be complex numbers with
$|\theta-\theta'-\pi| \leq \pi/6$ and $r'\leq 2r$. Then $|z+z'| \leq
r+\frac{r'}{2}$.
\end{lem}
\begin{proof}
We can assume that $\theta=0$. Then
  \begin{equation}
  |z+z'|^2 = (r+ r'\cos(\theta'))^2 + (r'\sin(\theta'))^2.
  \end{equation}
Since $\cos(\theta')\leq 0$ and $r'\leq 2r$, we have $r+r'
\cos(\theta') \in [-r,r]$. Moreover, $|\sin(\theta')|\leq 1/2$.
Hence,
  \begin{equation}
  |z+z'|^2 \leq r^2+ {r'}^2/4 \leq (r+r'/2)^2.
  \qedhere
  \end{equation}
\end{proof}
Together with \eqref{eq:phase_and_modulus}, the lemma proves that,
for all $x\in B(x_1,\delta/|t|)$,
  \begin{equation}
  |F(x)+G(x)| \leq |F(x)| + |G(x)|/2.
  \end{equation}
This proves that the second conclusion of Lemma \ref{lem:main_UNI}
holds.
\end{proof}

From this point on, we fix the constants $\zeta$ and $\delta$ given
by Lemma \ref{lem:main_UNI}. Since $\Delta$ is a John domain, there
exist constants $\Cb$ and $\eb$ such that, for all $\epsilon<\eb$,
for all $x\in \Delta$, there exists $x'\in \Delta$ such that
$d(x,x') \leq \Cb \epsilon$ and such that the ball $B(x',\epsilon)$
is compactly contained in $\Delta$. Choose $T_0\geq 10$ such that
$2(\zeta+\delta)/T_0 <\eb$.

\begin{lem}
\label{decays_hard} There exist $\betaa<1$ and $0<\sigma_2<\sigma_1$
satisfying the following property. Let $s=\sigma+it$ with $\sigma\in
[-\sigma_2,\sigma_2]$ and $|t|\geq T_0$. Let $(u,v)\in \EE_t$. Then
there exists $\tilde u : \Delta \to \R$ such that $(\tilde u, \tilde
L_s^n v) \in \EE_t$ and $\int \tilde u^2 \dd\mu \leq \betaa \int u^2
\dd\mu$.
\end{lem}
\begin{proof}
Consider a maximal set of points $x_1,\dots,x_k \in \Delta$ such
that the balls $B(x_i, 2(\zeta+\delta)/|t|)$ are compactly included
in $\Delta$, and two by two disjoint. 
By Lemma \ref{lem:precompact}, this set is finite.
The John domain condition on
$\Delta$ ensures that $\Delta$ is covered by the balls $B(x_i,
\Cc/|t|)$ where $\Cc=(2+\Cb)2(\zeta+\delta)$.

In each ball $B(x_i, (\zeta+\delta)/|t|)$, there exists a ball
$B'_i=B(x'_i, \delta/|t|)$ on which the conclusion of Lemma
\ref{lem:main_UNI} holds for the pair $(u,v)$. We will write
$\type(B'_i)=h$ if the first conclusion of Lemma \ref{lem:main_UNI}
holds, and $\type(B'_i)=k$ otherwise. By Lemma
\ref{lem_exists_bump}, there exists a function $\rho_i$ on $\Delta$
such that $\rho_i=1$ on $B''_i=B(x'_i, \delta/(\Cd |t|))$,
$\rho_i=0$ outside of $B'_i$ and $\norm{\rho_i}_{C^1}\leq \Ce
|t|/\delta$. We define a function $\rho$ on $\Delta$ by
  \begin{equation}
  \rho=\left(\sum_{\type(B'_i)=h} \rho_i\right) \circ T^n
  \end{equation}
on $h(\Delta)$,
  \begin{equation}
  \rho=\left(\sum_{\type(B'_i)=k} \rho_i\right) \circ T^n
  \end{equation}
on $k(\Delta)$, and $\rho=0$ on $\Delta \moins (h(\Delta) \cup
k(\Delta))$. This function satisfies $\norm{\rho}_{C^1} \leq
|t|/\eta_0$ for some constant $\eta_0$ independent of $s,u,v$, and
we can assume $\eta_0<1/4$. Notice that $\eta_0$ depends on
$n,h,k$ and $\delta$, which is not troublesome since these
quantities are fixed once and for all. Define a new function
$\chi= 1-\eta_0 \rho$. It takes its values in $[3/4,1]$, with
$\norm{D\chi} \leq |t|$. Moreover, by construction,
  \begin{equation}
  \label{eq:decreases}
  |\tilde L_s^n v| \leq \tilde L_\sigma^n (\chi u).
  \end{equation}
We set $\tilde u=\tilde L_\sigma^n (\chi u)$. By
\eqref{eq:decreases} and Lemma \ref{lem:controls_with_chi},
$(\tilde u, \tilde L_s^n v) \in \EE_t$.
We have to show that, for some constant $\betaa<1$, $\int \tilde u^2
\dd\mu \leq \betaa \int u^2\dd\mu$ as soon as $\sigma$ is small
enough.

The definition of $\tilde L_\sigma^n$ gives
  \begin{align*}
  \lambda_\sigma^{2n} f^2_\sigma(x) &\tilde u^2(x)
  =\left( \sum_{l\in
  \HH_n} e^{-\sigma r^{(n)}(lx)} J(lx) (\chi\cdot f_\sigma\cdot
  u)(lx)\right)^2
  \\&
  \leq \left( \sum_{l\in \HH_n} J(lx) (f_\sigma \cdot u^2)(lx)
  \right)
  \left( \sum_{l\in \HH_n} e^{-2\sigma r^{(n)}(lx)} J(lx)
  (f_\sigma\cdot \chi^2)(lx)\right)
  \\&
  \leq
  \left( \sup_{\Delta} \frac{f_\sigma}{f_0} \right)
  \left( \sum_{l\in \HH_n} J(lx) (f_0 \cdot u^2)(lx) \right)
  \left( \sup_{\Delta} \frac{f_\sigma}{f_{2\sigma}} \right)
  \left( \sum_{l\in \HH_n} e^{-2\sigma r^{(n)}(lx)} J(lx)
  (f_{2\sigma}\cdot \chi^2)(lx)\right).
  \end{align*}
If $x \in B''_i$ with $\type(B'_i)=h$, we have
  \begin{align*}
  \frac{1}{\lambda_{2\sigma}^n f_{2\sigma}(x)}
  \sum_{l\in \HH_n} e^{-2\sigma r^{(n)}(lx)} J(lx)
  (f_{2\sigma}\cdot \chi^2)(lx)
  &
  = \tilde L_{2\sigma}^n (\chi^2)(x)
  \\&
  = 1 - ( 1-(1-\eta_0)^2) e^{-2\sigma r^{(n)}(hx)} J(hx)
  \frac{f_{2\sigma}(hx)}{ \lambda_{2\sigma}^n f_{2\sigma}(x)}.
  \end{align*}
This is uniformly bounded by a constant $\eta_1<1$. The same
inequality holds if $\type(B'_i)=k$, with $h$ replaced by $k$.
Define a number
  \begin{equation}
  \xi(\sigma)= \left(\sup_{\Delta} \frac{\lambda_{2\sigma}^n f_0(x)
  f_{2\sigma}(x)}{\lambda_\sigma^{2n} f^2_\sigma(x)}\right) \left(
  \sup_{\Delta} \frac{f_\sigma}{f_0} \right) \left( \sup_{\Delta}
  \frac{f_\sigma}{f_{2\sigma}} \right).
  \end{equation}
Let $X=\bigcup B''_i$ and $Y=\Delta \moins X$.  We have proved that
  \begin{equation}
  \label{eq:on_X}
  \forall x \in X,\quad
  \tilde u^2(x) \leq \eta_1 \xi(\sigma) \tilde L_0^n(u^2)(x).
  \end{equation}
If $x\not\in X$, there is no cancellation mechanism, and we simply
have
  \begin{equation}
  \label{eq:on_Y}
  \forall x\in Y,\quad
  \tilde u^2(x) \leq \xi(\sigma) \tilde L_0^n(u^2)(x).
  \end{equation}
The equations \eqref{eq:on_X} and \eqref{eq:on_Y} are not
sufficient by themselves to obtain an inequality $\int \tilde
u^2\dd\mu \leq \betaa \int u^2\dd\mu$, one further argument is
required.

Since $\norm{Du} \leq 2\Ca |t| u$, $\norm{D(u^2)} \leq 4\Ca |t|
u^2$. Hence, $(u^2,u^2)\in \EE_{2t}$. By Lemma
\ref{lem:controls_with_chi}, we obtain $(\tilde L_0^n (u^2), \tilde
L_{2it}^n(u^2))\in \EE_{2t}$. Hence, the function $w=\tilde L_0^n
(u^2)$ satisfies $\norm{Dw}\leq 4\Ca |t| w$. Gronwall's inequality
then implies that, for all points $x,x'\in \Delta$, $w(x')\leq
w(x)e^{4\Ca |t|d(x,x')}$. In particular, there exists a constant $C$
such that, for all points $x,x'$ in a ball $B(x_i, \Cc/|t|)$,
$w(x')\leq C w(x)$. This yields
  \begin{equation}
  \frac{ \int_{B(x_i, \Cc/|t|)} w \dd\mu}{\mu(B(x_i,\Cc/|t|))}
  \leq C \frac{ \int_{B''_i} w \dd\mu}{\mu(B''_i)}.
  \end{equation}
Moreover, $\Leb(B(x_i,\Cc/|t|))/ \Leb(B''_i)$ is uniformly bounded
since $(\Delta,\Leb)$ is a John domain, and the density of $\mu$ is
bounded from above and below. We get another constant $C'$ such that
  \begin{equation}
  \int_{B(x_i,\Cc/|t|)} w \dd\mu \leq C' \int_{B''_i} w \dd\mu.
  \end{equation}
Since the balls $B''_i$ are disjoint, we obtain
  \begin{equation}
  \int_Y w \dd\mu\leq C' \int_X w\dd\mu.
  \end{equation}
Consider finally a large constant $A$ such that $(A+1)\eta_1+C' \leq
A$. With \eqref{eq:on_X} and \eqref{eq:on_Y}, we get
  \begin{multline*}
  (A+1)\int \tilde u^2\dd\mu \leq
  \xi(\sigma)\left[  (A+1)\int_X \eta_1 w\dd\mu + (A+1)
  \int_Y w\dd\mu\right]
  \\
  \leq \xi(\sigma) \left[ (A+1)\eta_1
  \int_X w\dd\mu + A\int_Y w\dd\mu + C' \int_X w\dd\mu \right]
  \leq \xi(\sigma) A \int w \dd\mu.
  \end{multline*}
Since $\int w\dd\mu= \int \tilde L_0^n(u^2)\dd\mu=\int u^2\dd\mu$,
we finally get
  \begin{equation}
  \int \tilde u^2\dd\mu \leq \xi(\sigma)
  \frac{A}{A+1} \int u^2\dd\mu.
  \end{equation}
When $\sigma\to 0$, $\xi(\sigma)$ converges to $1$. Hence, there
exists $\sigma_2>0$ such that $\betaa=\sup_{|\sigma|\leq
\sigma_2}\xi(\sigma)
  \frac{A}{A+1}$ is $<1$.
\end{proof}

Lemmas \ref{decays_easily} and \ref{decays_hard} easily imply
Proposition \ref{thm_UNI}:

\begin{proof}[Proof of Proposition \ref{thm_UNI}]
Is is sufficient to prove that there exist $\beta<1$ and $C>0$ such
that, for all $m\in \N$, for all $s=\sigma+it$ with $\sigma$ small
enough and $|t|\geq T_0$, for all $u\in C^1(\Delta)$,
  \begin{equation}
  \label{contracts_2m}
  \norm{\tilde L_s^{2 m n} u}_{L^2(\mu)} \leq C \beta^m \norm{u}_{1,t}.
  \end{equation}
Indeed, if \eqref{contracts_2m} is proved, consider a general
integer $k$ and write it as $k=2mn+r$ where $0 \leq r \leq 2n-1$.
Then
  \begin{equation}
  \norm{ L_s^k u}_{L^2(\Leb)} \leq C \lambda_\sigma^k \norm{\tilde L_s^k
  u}_{L^2(\mu)}
  \leq C \lambda_\sigma^k \beta^m \norm{ \tilde L_s^r u}_{1,t}
  \leq C \lambda_\sigma^k \beta^m \norm{u}_{1,t},
  \end{equation}
by Lemma \ref{acts_boundedly}. Choosing $\sigma'_0$ small enough so
that $\sup_{ |\sigma|\leq \sigma'_0} \lambda_\sigma \beta^{1/(2n)}
<1$, we obtain the full conclusion of Proposition \ref{thm_UNI}.

Let us prove \eqref{contracts_2m} for $u\in C^1(\Delta)$. Suppose
first that, for all $0\leq p <m$, $\norm{D(\tilde L_s^{pn} u)}_{C^0}
\geq 2\Ca |t| \norm{ \tilde L_s^{pn} u}_{C^0}$. Then Lemma
\ref{decays_easily} gives
  \begin{equation}
  \norm{\tilde L_s^{mn} u}_{1,t} \leq \left(\frac{9}{10}\right)^m
  \norm{u}_{1,t}.
  \end{equation}
Since $\norm{\tilde L_s^{2mn} u}_{L^2(\mu)} \leq  \norm{\tilde
L_s^{2mn} u}_{1,t} \leq C \norm{\tilde L_s^{mn}u}_{1,t}$ by Lemma
\ref{acts_boundedly}, \eqref{contracts_2m} is satisfied.

Otherwise, let $p<m$ be the first time such that $\norm{D(\tilde
L_s^{pn} u)}_{C^0} < 2\Ca |t| \norm{ \tilde L_s^{pn} u}_{C^0}$, and
let $v=\tilde L_s^{pn} u$. Since $(\sup |v|, v) \in \EE_t$, we can
apply Lemma \ref{decays_hard} and obtain a sequence of functions
$u_k$ with $u_0= \sup |v|$, $\int u_k^2 \dd\mu  \leq \beta_0^k \int
u_0^2 \dd\mu$, and $(u_k, \tilde L_s^{kn}v) \in \EE_t$. In
particular,
  \begin{equation}
  \norm{ \tilde L_s^{2mn}u}_{L^2(\mu)}
  = \norm{ \tilde L_s^{(2m-p) n} v}_{L^2(\mu)}
  \leq \norm{u_{2m-p}}_{L^2(\mu)}
  \leq \beta_0^{(2m-p)/2} \sup |v|
  \leq \beta_0^{m/2} \norm{u}_{C^0}.
  \end{equation}
This proves \eqref{contracts_2m} and concludes the proof of
Proposition \ref{thm_UNI}.
\end{proof}

\subsection{A control in the norm $\norm{\cdot}_{1,t}$}

Although it will not be useful in this paper, it is worth
mentioning that Proposition \ref{thm_UNI}, which gives a control
in the $L^2$ norm, easily implies an estimate in the stronger norm
$\norm{\cdot}_{1,t}$. This kind of estimate is especially useful
for the study of zeta functions.

\begin{prop}
There exist $\sigma'_0 \leq \sigma_0$, $T_0>0$, $C>0$ and
$\beta<1$ such that, for all $s=\sigma+it$ with $|\sigma|\leq
\sigma'_0$ and $|t|\geq T_0$, for all $u\in C^1(\Delta)$, for all
$k \in \N$,
  \begin{equation}
  \norm{L_s^k u}_{1,t} \leq C\lambda_\sigma^k
  \min(1, \beta^k |t|) \norm{u}_{1,t}.
  \end{equation}
\end{prop}
\begin{proof}
It is sufficient to prove the existence of $\beta<1$ such that
  \begin{equation}
  \label{eq:3mn}
  \norm{\tilde L_s^{3k} u}_{1,t} \leq C \beta^k |t| \norm{u}_{1,t}
  \end{equation}
if $|\sigma|$ is small enough and $|t|$ is large enough. Indeed,
together with Lemma \ref{acts_boundedly}, it implies the conclusion
of the proposition.

Denote by $\Lip(\Delta)$ the set of Lipschitz functions on $\Delta$,
with its canonical norm
  \begin{equation}
  \norm{w}_{\Lip}=\sup_{x\in \Delta}|w(x)| + \sup_{x\not=x'}
  \frac{|w(x)-w(x')|}{d(x,x')}.
  \end{equation}
We will use the following classical Lasota-Yorke inequality on the
transfer operators $\tilde L_\sigma$, for small enough $|\sigma|$:
there exist $C>0$ and $\beta_1<1$ such that, for all $k\in \N$, for
all $w\in \Lip(\Delta)$,
  \begin{equation}
  \norm{\tilde L_\sigma^k w}_{\Lip} \leq C \beta_1^k \norm{w}_{\Lip}+
  C \norm{w}_{L^1}.
  \end{equation}
Hence,
  \begin{equation}
  \norm{\tilde L_s^{2k} u}_{C^0} \leq \norm{\tilde
  L_\sigma^{k}(|\tilde L_s^{k} u|)}_{C^0}
  \leq C \beta_1^{k} \norm{\tilde L_s^{k} u}_{\Lip}+C
  \norm{\tilde L_s^{k} u}_{L^2}.
  \end{equation}
Moreover, $\norm{\tilde L_s^{k} u}_{\Lip} \leq |t|\norm{\tilde
L_s^{k} u}_{1,t} \leq C |t| \norm{u}_{1,t}$, and $\norm{\tilde
L_s^{k} u}_{L^2} \leq \beta_2^{k} \norm{u}_{1,t}$ for some
$\beta_2<1$, by Proposition \ref{thm_UNI}. Hence, there exists
$\beta_3<1$ such that
  \begin{equation}
  \label{eq:2mn}
  \norm{\tilde L_s^{2k} u}_{C^0} \leq C|t| \beta_3^k \norm{u}_{1,t}.
  \end{equation}
By Lemma \ref{acts_boundedly}, we get
  \begin{equation}
  \norm{\tilde L_s^{3k} u}_{1,t} \leq C \norm{\tilde L_s^{2k}
  u}_{C^0} + \frac{\expansion^{-k}}{|t|} \norm{D(\tilde
  L_s^{2k}u)}_{C^0}.
  \end{equation}
Notice that $\frac{\norm{D(\tilde L_s^{2k}u)}_{C^0}}{|t|}\leq
\norm{\tilde L_s^{2k}u}_{1,t} \leq C \norm{u}_{1,t}$. Together with
\eqref{eq:2mn}, this implies \eqref{eq:3mn} and concludes the proof
of the proposition.
\end{proof}

\subsection{Proof of Theorem \ref{main_thm_expanding}}

\label{sec:proof}

Let $U\in \BB_0$ and $V\in \BB_1$ be such that $\int V\dd\mu_r=0$.
We will prove that there exist $\delta>0$ independent of $U,V$, and
$C>0$ dependent of $U,V$ such that
  \begin{equation}
  \forall t \geq 0, \quad \left|\int U \cdot V\circ T_t \right|
  \leq C e^{-\delta t}.
  \end{equation}
By the closed graph theorem, this will imply Theorem
\ref{main_thm_expanding}.

For $t\geq 0$, let $A_t=\{ (x,a) \in \Delta_r \tq a+t\geq r(x)\}$
and $B_t= \Delta_r \moins A_t$. Then
  \begin{equation}
  \int U \cdot V\circ T_t = \int_{A_t} U\cdot V\circ T_t + \int_{B_t}
  U\cdot V\circ T_t =: \rho(t)+ \bar \rho(t).
  \end{equation}
We have
  \begin{equation}
  \label{eq:barrhot}
  | \bar \rho(t)| \leq C \int_{x\in \Delta} \max(r(x)-t,0)
  \leq C \int_{r(x)\geq t} r(x) \leq C \norm{r}_{L^2} ( \Leb(
  x\tq r(x)>t))^{1/2}.
  \end{equation}
Since $r$ has exponentially small tails, this quantity decays
exponentially. Hence, it is sufficient to prove that $\rho(t)$
decays exponentially to conclude.

Since $\rho(t)$ is bounded, we can define, for $s\in \C$ with $\Re
s>0$,
  \begin{equation}
  \hat \rho(s)=\int_0^\infty e^{-st} \rho(t)\dd t.
  \end{equation}
For $W:\Delta_r \to \R$ and $s\in \C$, set $\hat
W_s(x)=\int_0^{r(x)} W(x,a) e^{-sa}\dd a$ when $x\in \Delta$.

\begin{lem}
\label{formula_hatrho} Let $s\in \C$ with $\Re s>0$. Then
  \begin{equation}
  \hat \rho(s)= \sum_{k=1}^\infty \int_{\Delta} \hat V_s(x)\cdot
  (L_s^k \hat U_{-s})(x) \dd x.
  \end{equation}
\end{lem}
\begin{proof}
We compute
  \begin{align*}
  \hat \rho(s) &=
  \int_{x\in \Delta} \int_{a=0}^{r(x)} \int_{t+a\geq
  r(x)}e^{-st}  U(x,a) V\circ T_t(x,a) \dd t \dd a \dd x
  \\&
  = \sum_{k=1}^\infty \int_{x\in \Delta} \int_{a=0}^{r(x)}
  \int_{b=0}^{r(T^k x)} U(x,a) V(T^k x, b) e^{- s( b+r^{(k)}(x) -a)}
\dd b \dd a \dd x
  \\&
  = \sum_{k=1}^\infty \int_{x\in \Delta}\hat U_{-s}(x) e^{-s r^{(k)}(x)}
  \hat V_s(T^k x) \dd x
  \\&
  = \sum_{k=1}^\infty \int_{x\in \Delta} \hat V_s(x) (L_s^k \hat
  U_{-s})(x) \dd x.
  \qedhere
  \end{align*}
\end{proof}

\begin{lem}
\label{lem:borne_1t} There exists $C>0$ such that, for all
$s=\sigma+it$ with $|\sigma|\leq \sigma_0/4$ and $t\in \R$, the
function $L_s \hat U_{-s}$ is $C^1$ on $\Delta$
and satisfies the inequality
  \begin{equation}
  \norm{ L_s \hat U_{-s}}_{1,t} \leq \frac{C}{\max(1,|t|)}.
  \end{equation}
\end{lem}
\begin{proof}

Let us first prove that there exists $C>0$ such that, whenever
$|\sigma|\leq \sigma_0/4$,
  \begin{equation}
  \label{eq:hkasf}
  \forall x\in\Delta, \quad \left| \hat U_{-s}(x)
  \right| \leq \frac{C}{\max(1,|t|)} e^{(\sigma_0/2) r(x)}.
  \end{equation}
Since $\hat U_{-s}(x)=\int_{a=0}^{r(x)} U(x,a) e^{sa}\dd a$, this is
trivial if $|t|\leq 1$. If $|t|>1$,
an integration by parts gives
  \begin{equation}
  \hat U_{-s}(x)=\int_{a=0}^{r(x)} U(x,a) e^{sa}\dd a = \left[ U(x,a)
  \frac{e^{sa}}{s}\right]_{0}^{r(x)} - \int_{a=0}^{r(x)} \partial_t
  U(x,a) \frac{e^{sa}}{s} \dd a.
  \end{equation}
The boundary terms are bounded by $C e^{(\sigma_0/4) r(x)}/|t|$,
while the remaining term is at most
  \begin{equation}
  C r(x) e^{(\sigma_0/4) r(x)}
  /|t| \leq C' e^{(\sigma_0/2) r(x)} /|t|.
  \end{equation}

This proves \eqref{eq:hkasf}.

We can now compute
  \begin{equation}
  \label{eq:3sigma4}
  |L_s \hat U_{-s}(x)| = \left| \sum_{h\in \HH} e^{-s r\circ h(x)}
  J(hx) \hat U_{-s}(hx) \right|
  \leq \frac{C}{\max(1,|t|)} \sum_{h\in \HH} e^{(\sigma_0/4) r(hx)} J(hx)
  e^{(\sigma_0/2) r(hx)}.
  \end{equation}
This sum is bounded by $\frac{C'}{\max(1,|t|)}$ since $\sigma_0/2+\sigma_0/4
< \sigma_0$.

We have $L_s \hat U_{-s}(x) =\sum_{h\in \HH} e^{-s r\circ h(x)}
J(hx) \int_{a=0}^{r(hx)} U(hx,a) e^{sa} \dd a$. To obtain $D(L_s
\hat U_{-s})(x)$, we can differentiate $e^{-s r\circ h(x)}$, or
$J(hx)$, or $U(hx,a)$ in the integral, or the bound $r(hx)$ of the
integral.

Since $D(e^{-s r\circ h(x)})=-s D(r\circ h)(x) e^{-s r\circ h(x)}$,
and $\norm{D(r\circ h)}$ is uniformly bounded, the corresponding
term is bounded by $C|s|\cdot C/\max(1,|t|)$, by the computation done in
\eqref{eq:3sigma4}. Since $D(J\circ h)(x) \leq C J(hx)$, the
corresponding term is bounded by $C/\max(1,|t|)$. If we differentiate
$U(hx,a)$ in the integral, then the corresponding term is bounded by $C
\sum_{h\in \HH} e^{(\sigma_0/4) r(hx)}J(hx) e^{(\sigma_0/4)
r(hx)} r(hx)$, which is still uniformly bounded. Finally, the last
term satisfies a similar bound.

We have proved that $\norm{ D( L_s \hat U_{-s})}_{C^0} \leq C$ for
some constant $C$. Together with the inequality $\norm{L_s \hat
U_{-s}}_{C^0} \leq C/\max(1,|t|)$, it proves the lemma.
\end{proof}

\begin{lem}
\label{lemme:borneL1} There exists $C>0$ such that, for
$s=\sigma+it$ with $|\sigma|\leq \sigma_0/4$ and $t\in \R$,
  \begin{equation}
  \norm{\hat  V_s}_{L^2} \leq \frac{C}{\max(1,|t|)}.
  \end{equation}
\end{lem}
\begin{proof}
The inequality  \eqref{eq:hkasf} for $\hat
V_s$ is trivial if $|t|\leq 1$, and can be proved by
an integration by parts along the flow direction (using the bounded
variation of $t\mapsto V(x,t)$) if $|t|>1$.
This concludes the proof since $\int_\Delta e^{\sigma_0 r}
<\infty$.
\end{proof}

\begin{cor}
\label{corsigma4} There exists $\sigma_3>0$ (independent of $U,V$)
such that the function $\hat \rho$ admits an analytic extension
$\phi$ to the set $\{ s=\sigma+it \tq |\sigma|\leq \sigma_3, |t|
\geq T_0\}$. This extension satisfies $|\phi(s)|\leq C / t^2$.
\end{cor}
\begin{proof}
For $s=\sigma+it$ with $|\sigma|\leq \sigma_0/4$ and $|t|\geq T_0$,
set $\phi(s)=\sum_{k=1}^\infty \int \hat V_s \cdot L_s^k \hat
U_{-s}$. By Lemma \ref{formula_hatrho}, it coincides with $\hat
\rho$ when $\Re s>0$.

We have to check that the series defining $\phi$ is summable, and
that $\phi$ satisfies the bound $|\phi(s)|\leq C/t^2$. By
Proposition \ref{thm_UNI}, Lemma \ref{lem:borne_1t} and Lemma
\ref{lemme:borneL1}, if $|\sigma|$ is small enough,
  \begin{equation}
  \left| \int \hat V_s \cdot L_s^k \hat U_{-s} \right|
  \leq \norm{\hat V_s}_{L^2} \norm{ L_s^k \hat U_{-s}}_{L^2}
  \leq \norm{\hat V_s}_{L^2}C \beta^{k-1}  \norm{L_s \hat
  U_{-s}}_{1,t}
  \leq \frac{C}{t^2} \beta^k.
  \end{equation}
This last term is summable and its sum is at most
$\frac{C}{(1-\beta)t^2}$.
\end{proof}

\begin{lem}
\label{lemOs} For all $s=it \not=0$, there exists an open disk $O_s$
with center $s$ (independent of $U,V$) such that $\hat \rho$ admits
an analytic extension to $O_s$.
\end{lem}
\begin{proof}
The operator $L_s$ acting on $C^1$ satisfies a Lasota-Yorke
inequality, by Lemma \ref{define_Ca} and the compactness of the unit
ball of $C^1(\Delta)$ in $C^0(\Delta)$. By Hennion's Theorem \cite{H},
its spectral radius 
on $C^1$ is $\leq 1$, and its essential spectral radius is $<1$.

Let us prove that $L_s$ has no eigenvalue of modulus $1$. This is an
easy consequence of the weak-mixing of the flow $T_t$, but we will
rather derive it directly. Assume that there exists a nonzero $C^1$
function $u$ and a complex number $\lambda$ with $|\lambda|=1$ such
that $L_s u = \lambda u$. Then $|u|=|L_s u|\leq L_0 |u|$. Since
$\int |u|=\int L_0 |u|$, we get $|u|=L_0 |u|$. In particular, $|L_s
u|=L_0 |u|$, which means that all the complex numbers $e^{-it r(hx)}
u(hx)$ have the same argument. Take $k\in \N$ such that $k|t|\geq
T_0$. The complex numbers $e^{-itk r(hx)} u^k(hx)$ also have the
same argument. Hence, $|L_{ks}(u^k)|=L_0 |u^k|$. In the same way,
for any $n\in \N$, $|L_{ks}^n( u^k)| = L_0^n |u^k|$. This is a
contradiction, since $L_{ks}^n(u^k)$ tends to $0$ in $L^2$ by
Proposition \ref{thm_UNI}, while $\int L_0^n |u^k|=\int |u^k|$ does
not tend to $0$ when $n\to\infty$.

We have proved that the spectral radius of $L_s$ is $<1$. Hence,
there exists a disk $O_s$ around $s$ and constants $C>0$, $r<1$ such
that, for all $s'\in O_s$ and for all $n\in \N$,
$\norm{L_{s'}^n}_{C^1} \leq C r^n$. Since $L_{s'} \hat U_{-s'}$
is uniformly bounded in $C^1$ by Lemma \ref{lem:borne_1t},
the series $\sum_{k\geq
1} \int_\Delta \hat V_{s'}\cdot L_{s'}^{k-1}( L_{s'} \hat U_{-s'})$ 
is convergent
on $O_s$. By Lemma \ref{formula_hatrho}, it coincides with $\hat
\rho(s')$ for $\Re s'>0$.
\end{proof}

\begin{lem}
\label{lemO0} There exists an open disk $O_0$ with center $0$
(independent of $U,V$) such that $\hat \rho$ admits an analytic
extension to $O_0$.
\end{lem}
\begin{proof}
The transfer operator $L_0$ acting on $C^1$ has an isolated
eigenvalue $1$. For small $s$, $L_s$ is an analytic perturbation of
$L_0$. Hence, it admits an eigenvalue $\lambda_s$ close to $1$.
Denote by $P_s$ the corresponding spectral projection, and $f_s$ the
eigenfunction (normalized so that $\int f_s=1$). On a disk $O_0$
centered in $0$, it is possible to write $L_s= \lambda_s P_s + R_s$
where $P_s$ and $R_s$ commute, and $\norm{R_s^n}_{C^1}\leq C r^n$
for some uniform constants $C>0$ and $r<1$.

The function $s\mapsto \lambda_s$ is analytic in $O_0$, let us
compute its derivative at $0$. Since $\norm{L_s-L_0}_{C^1}=O(s)$ and
$\norm{f_s-f_0}_{C^1}=O(s)$, we have
  \begin{align*}
  \lambda_s &= \int L_s f_s = \int (L_s-L_0)(f_s-f_0)+\int L_0(f_s-f_0)
  +\int L_s f_0
  \\ &
  = O(s^2)+ \int (f_s-f_0) + \int L_0( e^{-sr}f_0)
  = O(s^2) + 0 + \int e^{-sr}\dd\mu
  = 1 -s \int r\dd\mu +O(s^2).
  \end{align*}
Hence, $\lambda'(0)=-\int r\dd\mu \not=0$. Shrinking $O_0$ if
necessary, we can assume that $\lambda_s$ is equal to $1$ only for
$s=0$.

For $s\in O\moins \{0\}$, define a function
  \begin{equation}
  \phi(s)= \frac{1}{1-\lambda_s} \int_\Delta \hat V_s \cdot
  P_s L_s \hat U_{-s} + \sum_{k=0}^\infty \int_\Delta
  \hat V_s \cdot R_s^k L_s \hat U_{-s},
  \end{equation}
where the last series is uniformly converging since
$\norm{R_s^k}_{C^1}\leq C r^k$ and $\norm{L_s \hat
U_{-s}}_{C^1(\Delta)} \leq C$ by Lemma \ref{lem:borne_1t}. 
It coincides with $\hat \rho(s)$ when $\Re
s>0$. When $s\to 0$, the function $\frac{1}{1-\lambda_s}$
has a pole of order exactly one, since $\lambda'(0)\not=0$. Let us
show that $\int \hat V_0 \cdot P_0 L_0 \hat U_0=0$. This will conclude
the proof, since the function $\phi$, being bounded on a
neighborhood of $0$, can then be extended analytically to $0$.

The function $P_0 L_0 \hat U_0$ is proportional to $f_0$. Hence, it is
sufficient to prove $\int \hat V_0 f_0=0$. But
  \begin{equation}
  \int_\Delta \hat V_0(x)  f_0(x)\dLeb(x)= \int_{x\in \Delta}
  \int_{t=0}^{r(x)} V(x,t) \dd t \dd \mu(x) = \int V \dd\mu_r=0.
  \qedhere
  \end{equation}
\end{proof}

We will use the following classical Paley-Wiener theorem:
\begin{thm}
\label{thm:paleywiener} Let $\rho: \R_+ \to \R$ be a bounded
measurable function. For $\Re s>0$, define $\hat
\rho(s)=\int_{x=0}^\infty e^{-sx} \rho(x) \dd x$. Suppose that $\hat
\rho$ can be analytically extended to a function $\phi$ on a strip
$\{ s=\sigma+it \tq |\sigma|< \epsilon, t\in \R\}$ and that
  \begin{equation}
  \int_{t=-\infty}^{\infty} 
  \sup_{|z|< \epsilon} |\phi(z+it)|\dd t <
  \infty.
  \end{equation}
Then there exist a constant $C>0$ and a full measure subset $A
\subset \R_+$ such that, for all $x\in A$, $|\rho(x)| \leq C
e^{-(\epsilon/2) x}$.
\end{thm}
\comm{
\begin{proof}
Take $\sigma>0$. Then $\phi(\sigma+it)=\int_0^\infty
e^{-itx}e^{-\sigma x}\rho(x)\dd x$. By Fourier inversion formula,
there exists a set $A_\sigma\subset \R_+$ of full measure such that,
for $x\in A_\sigma$, $e^{-\sigma
x}\rho(x)=\frac{1}{2i\pi}\int_{t=-\infty}^{\infty} \phi(\sigma+it)
e^{itx}\dd t$.

Let $A=\bigcap A_{1/k}$, it has full measure. Fix $x\in A$. The
function $G: (z,t) \mapsto \phi(z+it) e^{itx}$, defined on
$B(0,\epsilon)\times \R$, is analytic in the first variable and
$|G(z,t)| \leq A(t)$ integrable, by assumption. Hence, $H: z\mapsto
\int_{t=-\infty}^\infty \phi(z+it) e^{itx} \dd t$ is analytic on
$B(0,\epsilon)$ (and bounded by a constant $C$ independent of $x$).
For $z=1/k$, it coincides with $2i\pi
e^{-zx}\rho(x)$. Hence, these two functions coincide on
$B(0,\epsilon)$. In particular, taking $z=-\epsilon/2$, we get
  \begin{equation}
  \rho(x)=e^{-(\epsilon/2)x} H(-\epsilon/2),
  \end{equation}
which is bounded by $C e^{-(\epsilon/2)x}$.
\end{proof}
}
\begin{proof}[Proof of Theorem \ref{main_thm_expanding}]

We can summarize Corollary \ref{corsigma4}, Lemma \ref{lemOs} and
Lemma \ref{lemO0} as follows: there exists $\sigma_4>0$ (independent
of $U,V$) such that $\hat \rho$ admits an analytic extension $\phi$
to the set $\{ s=\sigma+it \tq |\sigma|\leq \sigma_4, t\in \R\}$.
Moreover, there exists $C>0$ such that this extension satisfies
  \begin{equation}
  |\phi(\sigma+it)| \leq C \min\left( 1, \frac{1}{t^2}\right).
  \end{equation}
Together with Theorem \ref{thm:paleywiener}, it implies that
$\rho(t)$ decays exponentially on a subset of $\R_+$ of full
measure. Hence, on a full measure subset of $\R_+$, $\left| \int U
\cdot V \circ T_t\right|\leq C e^{-\delta t}$. Since $t\mapsto \int
U \cdot V\circ T_t$ is continuous by dominated convergence, this
inequality holds in fact everywhere. This concludes the proof of
Theorem \ref{main_thm_expanding}.
\end{proof}

\section{Exponential mixing for hyperbolic semiflows}

\label{par:define_hyperbolic}

In this section, we will prove Theorem \ref{main_thm_hyperbolic},
using Theorem \ref{main_thm_expanding} and an approximation
argument.

\subsection{Estimates on bad returns}

In this paragraph, we will prove the following exponential estimate on
the number of returns to the basis $\Delta$:
\begin{lem}
\label{lem_renouv_temp} Let $\Psi_t(x,a)$ be the number of returns
to $\Delta$ of $(x,a)$ before time $t$, i.e.,
  \begin{equation}
  \Psi_t(x,a)=\sup \{ n \in \N \tq a+t>r^{(n)}(x)\}.
  \end{equation}
For all $\expansion>1$, there exist $C>0$ and $\delta>0$ such that,
for all $t \geq 0$,
  \begin{equation}
  \int_{\Delta_r} \expansion^{-\Psi_t(x,a)}\dLeb \leq C e^{-\delta t}.
  \end{equation}
\end{lem}

\begin{proof}
We have
  \begin{align*}
  \int_{\Delta_r} \expansion^{ - \Psi_t(x,a)} \dLeb &
  = \sum_{n=0}^\infty \expansion^{-n} \Leb\Bigl\{ (x,a) \tq r^{(n)}(x) <
  a+t \leq r^{(n+1)}(x)\Bigr\}
  \\ &
  \leq \sum_{n=0}^\infty \expansion^{-n} \Leb\Bigl\{ (x,a) \tq t\leq
  r^{(n+1)}(x)\Bigr\}.
  \end{align*}
Moreover, for $\sigma>0$,
  \begin{align*}
  \Leb\{ (x,a) \tq t \leq r^{(n+1)}(x)\}&
  =\int_\Delta r(x) 1_{r^{(n+1)}(x) \geq t}
  \leq \left( \int_\Delta r^2\right)^{1/2} 
  \left( \int_\Delta 1_{ r^{(n+1)}(x)
  \geq t}\right)^{1/2}
  \\ &
  \leq C \left( \int_\Delta e^{\sigma r^{(n+1)}(x)}/e^{\sigma t}\right)^{1/2}.
  \end{align*}
If $\sigma$ is small enough,
  \begin{equation}
  \int_\Delta e^{\sigma r^{(n+1)}(x)} = \int_\Delta L_\sigma^{n+1}(1)
  \leq C \lambda_\sigma^{n+1}.
  \end{equation}
Choosing $\sigma$ small enough so that $\expansion^{-1}
\sqrt{\lambda_\sigma}<1$, we obtain $\int_{\Delta_r} \expansion^{ -
\Psi_t(x,a)} \dLeb \leq C e^{-\sigma t/2}$.
\end{proof}

\subsection{Proof of Theorem \ref{main_thm_hyperbolic}}

Let $U,V$ be $C^1$ functions on $\widehat{\Delta}_r$, with $\int U
\dd\nu_r=0$. We will prove that $\int U\cdot V \circ
\widehat{T}_{2t}\dd\nu_r$ decreases exponentially fast in $t$.

Define a function $V_t$ on $\Delta_r$ by $V_t(x,a)=\int_{y\in
\pi^{-1}(x)} V\circ \widehat{T}_t(y,a) \dd\nu_x(y)$. Let $\pi_r :
\widehat{\Delta}_r \to \Delta_r$ be given by
$\pi_r(y,a)=(\pi(x),a)$.

\begin{lem}
\label{lem_delta1} There exist $\delta>0$ (independent of $U,V$) and
$C>0$ such that, for all $t\geq 0$,
  \begin{equation}
  \left| \int_{\widehat{\Delta}_r}
  U\cdot V\circ \widehat{T}_{2t}\dd\nu_r -
  \int_{\widehat{\Delta}_r} U\cdot V_t\circ T_t \circ
  \pi_r \dd\nu_r\right| \leq C e^{-\delta t}.
  \end{equation}
\end{lem}
\begin{proof}
We have
  \begin{multline*}
  \left| \int U\cdot V\circ \widehat{T}_{2t}\dd\nu_r
  - \int U\cdot V_t\circ T_t \circ
  \pi_r \dd\nu_r \right|
  =\left| \int U \cdot (V\circ \widehat{T}_t - V_t\circ \pi_r) \circ
  \widehat{T}_t \dd\nu_r\right|
  \\
  \leq C \int |V\circ \widehat{T}_t - V_t\circ \pi_r| \circ
  \widehat{T}_t \dd\nu_r
  = C \int  |V\circ \widehat{T}_t - V_t\circ \pi_r|\dd\nu_r.
  \end{multline*}
Take $x\in \Delta$. If $\pi(y)=\pi(y')=x$, the contraction
properties of $T$ give $d(\widehat{T}_t (y,a), \widehat{T}_t (y',a))
\leq \expansion^{-\Psi_t(x,a)} d(y,y')$, where $\Psi_t$ is defined in
Lemma \ref{lem_renouv_temp}. Hence, $|V\circ
\widehat{T}_t(y,a)-V\circ \widehat{T}_t(y',a)|\leq C
\expansion^{-\Psi_t(x,a)}$. Averaging over $y'$, we obtain $|V\circ
\widehat{T}_t(y,a) - V_t(x,a)| \leq C \expansion^{-\Psi_t(x,a)}$.
Finally,
  \begin{equation}
  \int_{\widehat{\Delta}_r}  |V\circ \widehat{T}_t - V_t\circ \pi_r|
  \dd\nu_r \leq C
  \int_{\Delta_r}  \expansion^{-\Psi_t(x,a)}\dd\mu_r.
  \end{equation}
This quantity decays exponentially, by Lemma \ref{lem_renouv_temp}
(and since the density of $\mu$ is bounded). 
\end{proof}

\begin{lem}
\label{lem_delta2} There exist  $\delta>0$ (independent of $U,V$)
and $C>0$ such that, for all $t\geq 0$,
  \begin{equation}
  \left| \int_{\widehat{\Delta}_r}
  U\cdot V_t \circ T_t \circ \pi_r \dd\nu_r\right| \leq C e^{-\delta t}.
  \end{equation}
\end{lem}
\begin{proof}
Define a function $\bar U$ on $\Delta_r$ by $\bar U(x,a)=\int_{y\in
\pi^{-1}(x)} U(y,a)\dd\nu_x(y)$. Since $U\in
C^1(\widehat{\Delta}_r)$ and the measures $\nu_x$ satisfy the third
property in the definition of hyperbolic skew-products, the function
$\bar U$ belongs to $\BB_0$. Moreover, $\int_{\Delta_r} \bar U \dd\mu_r =
\int_{\widehat{\Delta}_r} U\dd\nu_r=0$. Hence, Theorem
\ref{main_thm_expanding} (or rather Remark \ref{main_rmk_expanding})
gives
  \begin{equation}
  \left| \int_{\widehat{\Delta}_r} U\cdot V_t \circ T_t \circ \pi_r 
  \dd\nu_r\right|
  = \left| \int_{\Delta_r} \bar U \cdot V_t \circ T_t \dd\mu_r\right|
  \leq C e^{-\delta t} \norm{\bar U}_{\BB_0} \norm{V_t}_{\BB_1}.
  \end{equation}
To conclude the proof, it is thus sufficient to show that
$\norm{V_t}_{\BB_1}$ is uniformly bounded. First of all, since $V$
is bounded, $V_t$ is bounded.

Consider then $x\in \bigcup \Delta^{(l)}$. Take $0<a<r(x)$. If
$T_t(x,a)$ is not of the form $(x',0)$, then $V_t$ is differentiable
along the flow direction at $(x,a)$. Its derivative is given by
  \begin{equation}
  \int_{y\in \pi^{-1}(x)} (\partial_a V)(\widehat{T}_t(y,a)) \dd\nu_x(y),
  \end{equation}
since the flow is an isometry in the flow direction. In particular,
this derivative is bounded by $\norm{V}_{C^1}$.

There is a finite number of points $0<a_1<\dots<a_p<r(x)$ such that
$T_t(x,a_i)$ is of the form $(x',0)$. Indeed, since $r$ is uniformly
bounded from below by a constant $\epsilon_1$, there are at most
$\frac{r(x)}{\epsilon_1}+1$ such points. At each of these points,
$V_t$ has a jump of at most $2\norm{V}_{C^0}$. Finally, the
variation of $a\mapsto V_t(x,a)$ along the interval $(0,r(x))$ is at
most
  \begin{equation}
  \left(\frac{r(x)}{\epsilon_1}+1\right) 2\norm{V}_{C^0} + r(x)
  \norm{V}_{C^1}
  \leq C r(x) \norm{V}_{C^1}.
  \qedhere
  \end{equation}
\end{proof}

Lemmas \ref{lem_delta1} and \ref{lem_delta2} show that, for a
uniform constant $\delta>0$ and for some constant $C>0$ depending on
$U$ and $V$, for all $t\geq 0$,
  \begin{equation}
  \left| \int U\cdot V \circ \widehat{T}_{2t} \dd\nu_r \right| \leq C
  e^{-\delta t}.
  \end{equation}
By the closed graph theorem, the constant $C$ can be chosen of the
form $C' \norm{U}_{C^1} \norm{V}_{C^1}$ for a uniform constant $C'$.
This concludes the proof of Theorem \ref{main_thm_hyperbolic}.

\appendix

\section{A simple distortion estimate}

Here we present an alternative distortion estimate, Theorem \ref
{simpleyoccoz}, which is far from
optimal, but is enough to obtain exponential mixing, while being
based on a much simpler argument.  While much simpler, we
have only noticed it after obtaining the nearly optimal estimate.

For $\AA' \subset \AA$ non-empty,
let $\m_{\AA'}(q)=\min_{\alpha \in \AA'} q_\alpha$, and let
$\m(q)=\m_\AA(q)$.  The other notations are those of \S \ref {distortion}.

\begin{lemma}[Kerckhoff, \cite{K}] \label {T q alpha}

For every $T>0$, $q \in \R^\AA_+$, $\alpha \in \AA$, $\pi \in \RR$, we have
\be
P_q(\gamma \in \Gamma_\alpha(\pi),\, (B_\gamma \cdot q)_\alpha>T q_\alpha
\di \pi)<T^{-1},
\ee
where $\Gamma_\alpha(\pi)$ denotes the set of paths starting at $\pi$ with
no winner equal to $\alpha$.

\end{lemma}

\begin{proof}

Let $\Gamma^{(n)}_\alpha(\pi) \subset \Gamma_\alpha(\pi)$ denote the set of
paths of length at most $n$.  We prove the inequality for
$\Gamma^{(n)}_\alpha(\pi)$ by induction on $n$.  The case $n=0$ is clear. 
The case $n$ follows immediately from the case $n-1$ when none of the rows
of $\pi$ end with $\alpha$.  Assume for instance that the top row of $\pi$
ends with $\alpha$ and the bottom row with $\beta$.  Then every path $\gamma
\in \Gamma_\alpha^{(n)}(\pi)$ starts with the bottom arrow $\gamma_s$
starting at $\pi$.  Let $q'=B_{\gamma_s} \cdot q$.  We have
$q'_\alpha=q_\alpha+q_\beta$ and $P_q(\gamma_s \di \pi)=\frac {q_\alpha}
{q_\alpha'}$.  The inequality follows by the induction hypothesis.
\end{proof}

\begin{thm} \label {simpleyoccoz}

There exists $C>1$ such that for every $q \in \R^\AA_+$, if $\pi \in \RR$
\be
P_q(M(B_\gamma \cdot q)<C \min \{\m(B_\gamma \cdot q),\M(q)\}
\di \pi)>C^{-1}.
\ee

\end{thm}

\begin{proof}

For $1 \leq k \leq d$, let $\m_k(q)=\max \m_{\AA'}(q)$ where the maximum is
taken over all $\AA' \subset \AA$ such that $\#\AA'=k$.  In particular
$\m=\m_d$.  We will show that for $1 \leq k \leq d$ there exists $C>1$ such
that
\be \label {c+2}
P_q(\M(B_\gamma \cdot q)<C \min
\{\m_k(B_\gamma \cdot q),\M(q)\} \di \pi)>C^{-1}
\ee
(the case $k=d$ implying the desired statement).
The proof is by induction on $k$.  For $k=1$ it is obvious.  Assume that it
is proved for some $1 \leq k<d$ with $C=C_0$.
Let $\Gamma$ be the set of minimal paths $\gamma$ starting at $\pi$ with
$\M(B_\gamma \cdot q)<C_0 \min\{\m_k(B_\gamma \cdot q),\M(q)\}$.
Then there exists $\Gamma_1 \subset
\Gamma$ with $P_q(\Gamma_1 \di \pi)>C^{-1}_1$
and $\AA' \subset \AA$ with $\# \AA'=k$
such that if $\gamma \in \Gamma_1$ then $\m_k(B_\gamma \cdot q)=
\m_{\AA'} (B_\gamma \cdot q)$.

For $\gamma_s \in \Gamma_1$, choose a path $\gamma=\gamma_s\gamma_e$ with
minimal length such that $\gamma$ ends at a permutation $\pi_e$ such that
the top or the bottom row of $\pi_e$ (and possibly both)
ends by some element of $\AA
\setminus \AA'$.  Let $\Gamma_2$ be the collection of the $\gamma=\gamma_s
\gamma_e$ thus obtained.
Then $P_q(\Gamma_2 \di \pi)>C^{-1}_2$ and
$\M(B_\gamma \cdot q)<C_2 \M(B_{\gamma_s} \cdot q)$ for
$\gamma=\gamma_s \gamma_e \in \Gamma_2$.

Let $\Gamma_3$ be the set of paths
$\gamma=\gamma_s\gamma_e$ such that $\gamma_s \in \Gamma_2$,
the winner of the last arrow of
$\gamma_e$ belongs to $\AA'$, the winners of the other arrows of $\gamma_e$
belong to $\AA \setminus \AA'$, and we have
$(B_\gamma \cdot q)_\alpha \leq 2d (B_{\gamma_s} \cdot q)_\alpha$ for all
$\alpha \in \AA'$.  By Lemma \ref {T q alpha},
$P_q(\Gamma_3 \di \gamma_s)>\frac {1} {2}$, $\gamma_s \in \Gamma_2$, and
$P_q(\Gamma_3 \di \pi)>(2 C_2)^{-1}$.

Let $\gamma=\gamma_s\gamma_e \in \Gamma_3$, $\gamma_s \in \Gamma_2$.
If $\M(B_\gamma \cdot
q)>2d\M(B_{\gamma_s} \cdot q)$, we take $\gamma_1$ with $\gamma_s \leq
\gamma_1 \leq \gamma$, of minimal length such that $\M(B_{\gamma_1} \cdot
q)>2d \M(B_{\gamma_s} \cdot q)$; there exists $\alpha \in \AA \setminus
\AA'$ such that $\M(B_{\gamma_1} \cdot q)=(B_{\gamma_1} \cdot
q)_\alpha \leq 4d\M(B_{\gamma_s} \cdot q)$.
Moreover we have $\m_{\AA'}(B_{\gamma_1} \cdot q)>(C_0 C_2 4d)^{-1}
M(B_{\gamma_1} \cdot q)$ in this case.
If $\M(B_\gamma \cdot q) \leq 2d \M(B_{\gamma_s} \cdot q)$, the loser
$\alpha$ of the last arrow of $\gamma$ belongs to $\AA \setminus \AA'$ and
satisfies $(B_\gamma \cdot q)_\alpha \geq (C_0 C_2 2d)^{-1} \M(B_\gamma
\cdot q)$ which allows again to conclude: in any case there exists
$\gamma_1$ with $\gamma_s \leq \gamma_1 \leq \gamma_e$ and $\AA'_1$ with
$\#\AA'_1=k+1$ such that $\M(B_{\gamma_1} \cdot q) \leq 4 d C_0 C_2
\min \{\m_{\AA'_1}(B_{\gamma_1} \cdot q),\M(q)\}$.  Since the set $\Gamma_4$
of all $\gamma_1$ thus obtained satisfies $P_q(\Gamma_4 \di \pi) \geq
P_q(\Gamma_3 \di \pi)>(2 C_2)^{-1}$, (\ref {c+2}) holds with $k+1$
instead of $k$.
\end{proof}

\section{Spectral gap} \label {sg}

This section is concerned with the natural action of $\SL(2,\R)$ on a
connected component of a stratum $\CC^{(1)}$.  Though we have not used it
elsewhere in this paper, this action is very important in several works on
the Teichm\"uller flow, see for instance the work on Lyapunov exponents of
\cite {F}.

We recall that the mere {\it existence} of this
action has already important implications: for instance the action of
non-compact one-parameter subgroups (which are conjugate either to the
{\it Teichm\"uller flow} or the {\it horocycle flow}) is automatically
mixing with respect to any ergodic invariant measure
for the $\SL(2,\R)$ action.  Thus, ergodicity of the Teichm\"uller flow
(\cite {M}, \cite {V1}) with respect to the absolutely continuous invariant
measure on $\CC^{(1)}$ implies mixing (which can be obtained also directly
\cite {V2}).

Here we will show how our analysis of the Teichm\"uller flow can be used to
show that the $\SL(2,\R)$ action has a {\it spectral gap}.  To put this
concept in context, we recall some more general definitions.

\begin{definition}

Let $G$ be a (locally compact $\sigma$-compact)
group.  A (strongly continuous) unitary
representation of $G$ is said to have {\it almost invariant vectors}
if for every
$\epsilon>0$ and for every compact subset $K \subset G$, there exists a unit
vector $v$ such that $\|g \cdot v-v\|<\epsilon$ for all $g \in K$.

A unitary action which does not have almost invariant vectors is said to be
{\it isolated from the trivial representation}.

If $G$ is a semi-simple Lie group (such as $\SL(2,\R)$),
a representation which is isolated from
the trivial representation is also said to have a {\it spectral gap}.

\end{definition}

Given a probability preserving
action of $\SL(2,\R)$, it thus makes sense to ask whether
the corresponding unitary representation on $L^2_0$ (the space of
zero-average $L^2$ functions) has a spectral gap.  Ergodicity of the action
is of course a necessary condition, being equivalent to the inexistence of
invariant unit vectors.
It may happen for a group $G$ that any unitary representation which has
almost invariant vectors has indeed an invariant unit vector: this is one of
the equivalent definitions of Kazhdan's property (T), and has several
consequences.  As it is well known, $\SL(2,\R)$ does not have property (T),
so the spectral gap is indeed a non-automatic property in this case.

The spectral gap for the $\SL(2,\R)$ action on $\CC^{(1)}$ can be also seen
more geometrically as a statement about the foliated
Laplacian on $\CC^{(1)}/\SO(2,\R)$,\footnote
{The space $\CC^{(1)}/\SO(2,\R)$ is foliated by quotients of
$\SL(2,\R)/\SO(2,\R)$, which is a model for $2$-dimensional hyperbolic
space.  In particular there is a natural leafwise metric of constant
curvature $-1$, which allows us to define the foliated Laplacian, whose
spectrum is contained in $[0,\infty)$.} or of the Casimir
operator: the spectrum (for the action on $L^2_0$) does not contain $0$.

\comm{
This is an important non-automatic property of an $\SL(2,\R)$ action.
There are several equivalent ways to
define what it means for the $\SL(2,\R)$ action to have a spectral gap:
\begin{enumerate}
\item $0$ is an isolated eigenvalue of the foliated Laplacian on
$\CC^{(1)}/\SO(2,\R)$.\footnote
{The space $\CC^{(1)}/\SO(2,\R)$ is foliated by quotients of
$\SL(2,\R)/\SO(2,\R)$, which is a model for $2$-dimensional hyperbolic
space.  In particular there is a natural leafwise metric of constant
curvature $-1$, which allows us to define the foliated Laplacian, whose
spectrum is contained in $[0,\infty)$.  Ergodicity of the action implies
that $0$ is an eigenvalue of multiplicity $1$.}
\item The trivial representation is an isolated irreducible factor of the
unitary representation of $\SL(2,\R)$ on $L^2(\nu_{\CC^{(1)}})$.\footnote
{The unitary representations of $\SL(2,\R)$ have been classified,
c.f.\ 
\cite {Rt}.  The space of all irreducible representations 
has a natural topology.
There is a one-dimensional
{\it trivial representation}, by the identity,
all the other representations being infinite-dimensional. 
Ergodicity of the action means that the identity representation has
multiplicity $1$.}
\end{enumerate}
}

The connection between the spectral gap for the $\SL(2,\R)$ action and
rates of mixing for non-compact one-parameter subgroups was used most
notably by Ratner \cite {Rt}.  In her work, estimates on the
rates of mixing are deduced from the spectral gap.  That one could also go
the other way around seems to be also understood (the argument is much
easier than for the direction used by Ratner).  It is possible however
that this is the first time that it has been useful to
consider this connection in the other direction.

The existence of a spectral gap has several ramifications.  It is even
interesting to just ``go back'' to rates of mixing using the work of Ratner.
It implies {\it polynomial decay of correlations} for
the horocycle flow.  It even gives back extra information regarding the
Teichm\"uller flow: it implies that exponential mixing holds for observables
which are only H\"older along the $\SO(2,\R)$ orbits (this notion of
regularity is made precise in \cite {Rt}).  Further applications include
exponential estimates for the Ball Averaging Problem, see \cite {MNS}.

The initial line of the arguments given here (reduction to a ``reverse
Ratner estimate'') was explained to us by
Nalini Anantharaman, Sasha Bufetov and Giovanni Forni.  The proof of
the ``reverse Ratner estimate'' was explained to us by Giovanni Forni.

\begin{prop} \label {ergodic}

Let us consider an ergodic action of $\SL(2,\R)$ by measure-preserving
automorphisms of a probability space.  Let $\rho$ be the corresponding
representation on the space $H$ of $L^2$ zero average functions.
Assume that there exist $\delta\in (0,1)$ and a dense
subset of the subspace of $\SO(2,\R)$-invariant functions $H' \subset H$
consisting of functions $\phi$ for which the correlations
$\langle \phi,\rho(g_t) \cdot \phi \rangle$, $g_t=
\left (\begin{matrix} e^t & 0 \\ 0 & e^{-t} \end{matrix} \right )$,
decay like $O(e^{-\delta t})$.  Then $\rho$ is isolated from the trivial
representation.
\end{prop}

\begin{proof}

Let us decompose $\rho$ into irreducible representations.  Thus
$H=\int H_\xi \dd\mu(\xi)$, and there are irreducible actions $\rho_\xi$
of $\SL(2,\R)$ on each $H_\xi$ which integrate to $\rho$.

Bargmann's classification (see \cite {Rt} and the references therein)
shows that all non-trivial
irreducible representations fall into one of three series of
representations: the principal, the complementary and the discrete series.
Thus we have the corresponding decomposition $\mu=\mu_p+\mu_c+\mu_d$.  We
recall some basic facts that follow from this classification:
\begin{enumerate}
\item If $\rho_\xi$ is in the complementary series, then there exists
$s=s(\xi) \in (0,1)$, such that $\rho_\xi$ is isomorphic to the
following representation $\rho_s$: the Hilbert space is 
  \begin{equation}
  \HH_s=\left\{ f: \R \to \C \tq \norm{f}^2=\int_{\R\times \R} \frac{f(x)
  \overline{f(y)} }{ |x-y|^{1-s}}\dd x \dd y<\infty \right\},
  \end{equation}
and the action is given by
  \begin{equation}
  \rho_s \left(\begin{array}{cc} a & b\\ c& d \end{array}\right)f (x)=
  \frac{1}{(cx+d)^{1+s}} f\left( \frac{ax+b}{cx+d} \right).
  \end{equation}
\item The (integrated) representation $\rho$ is isolated from the trivial
representation if and only if
there exists $\epsilon>0$ such that $s(\xi)<1-\epsilon$
for $\mu_c$-almost every $\xi$.
\item The space of $\SO(2,\R)$ invariant vectors $H'_\xi \subset H_\xi$
is one-dimensional (in the
case of the principal and complementary series), or zero-dimensional (in the
case of the discrete series).
\end{enumerate}

\comm{
Bargmann's classification
(see \cite {Rt} and the references therein)
shows that for each $\xi$ there exists $s=s(\xi) \in \C$ such that
the action on $H_\xi$ is isomorphic to some specific irreducible action.
Here the possible values of $s$ are either purely imaginary (the principal
series), $s \in (0,1)$ (the complementary series), or $s$
an integer (the discrete series).  To say that the trivial
representation is isolated is to say that, for some $\epsilon>0$,
$\mu$ gives zero measure to the set $\{\xi \tq s(\xi)\in
(1-\epsilon,1)\}$.
It is a statement thus concerned only
with the representations in the complementary series.
}

Let $H' \subset H$ be the set of $\SO(2,\R)$ invariant functions.  Then
$H'=\int H'_\xi \dd\mu(\xi)$.  The point of the proof is the
following lemma:
\begin{lem}
\label{positive_correlations}
If $\rho_\xi$ is in the complementary series and $\phi_\xi \in H'_\xi$ is
a non-zero vector, then
$\langle \phi_\xi,\rho_\xi(g_t)
\cdot \phi_\xi \rangle$ is {\em positive} and
  \begin{equation}
  \label{eq:decay_exp}
  \lim_{t \to \infty} \frac {1} {t} 
  \ln \langle \phi_\xi,\rho_\xi(g_t) \cdot \phi_\xi \rangle = -1+s(\xi).
  \end{equation}
\end{lem}

Let us show how to conclude the proof using the lemma.
Suppose by contradiction that $\rho$ is not isolated from
the trivial representation. There exists a function
$\phi=\int \phi_\xi \dd\mu(\xi)\in H'$ whose correlations decay
like $O(e^{-\delta t})$ and such that
  \begin{equation}
  \label{charge_tout}
  \mu_c\{ \xi \tq \phi_\xi \not=0 \text{ and } s(\xi)\in (1-\delta/2,1)
  \}>0.
  \end{equation}
Write $\phi=\phi_p+\phi_c$ where $\phi_p$ is the part of $\phi$
corresponding to representations in the principal series, and
$\phi_c$ corresponds to the complementary series (as discussed above,
$\phi_\xi=0$ for $\mu_d$-almost every $\xi$). Then
  \begin{equation}
  \langle \phi,\rho(g_t) \cdot \phi \rangle
  = \langle \phi_p,\rho(g_t) \cdot \phi_p \rangle
  +
  \langle \phi_c,\rho(g_t) \cdot \phi_c \rangle.
  \end{equation}
By the results of Ratner \cite{Rt}, the correlations of $\phi_p$
decay at least as $t e^{-t}$. Moreover, by \eqref{eq:decay_exp}, positivity,
and \eqref{charge_tout}, the second term is larger than $C e^{-\delta
t/2}$ for large $t$.
This contradicts the speed of decay of correlations
of $\phi$.
\comm{
On ne peut pas toujours prendre une fonction $\phi$ telle que
$\phi_\xi \not=0$ pour presque tout $\xi$ : par exemple, dans
$L^2([0,1])$, l'ensemble des fonctions dont le support est de mesure
$<1$ est dense...
}
\end{proof}

\begin{proof}[Proof of Lemma \ref{positive_correlations}]
\comm{
Let $s=s(\xi)$. The representation $\rho_\xi$ is (isomorphic to) the
following one: the Hilbert space is 
  \begin{equation}
  \left\{ f: \R \to \C \tq \norm{f}^2=\int_{\R\times \R} \frac{f(x)
  \overline{f(y)} }{ |x-y|^{1-s}}\dd x \dd y<\infty \right\},
  \end{equation}
and the action of $\SL(2,\R)$ is given by
  \begin{equation}
  \left(\begin{array}{cc} a & b\\ c& d \end{array}\right)f (x)=
  \frac{1}{(cx+d)^{1+s}} f\left( \frac{ax+b}{cx+d} \right).
  \end{equation}
}
A function $f \in \HH_s$ is invariant under the $\SO(2,\R)$ action if and
only if it is smooth and 
satisfies the differential equation $(1+x^2)f'(x)+(1+s)xf (x)=0$,
i.e., $f(x)=\frac{c}{(1+x^2)^{(1+s)/2}}$.

For such a function $f$, the correlations are given by
  \begin{align}
  \label{eq:ratner}
  \langle f,\rho_s(g_t)\cdot f \rangle&
  = |c|^2 e^{t(1+s)} \int_{\R\times \R} \frac{\dd x\dd y}{ (1+x^2)^{(1+s)/2}
  (1+e^{4t}y^2)^{(1+s)/2} |x-y|^{1-s}}
  \\
\nonumber&
  =|c|^2 e^{t(-1+s)} \int_{\R \times \R} \frac{\dd x\dd y}{ (1+x^2)^{(1+s)/2}
  (1+y^2)^{(1+s)/2} |x-e^{-2t}y|^{1-s}}.
  \end{align}
This shows that the correlations are positive and that
  \begin{equation}
  \label{eq:ratner2}
  \liminf e^{t(1-s)}\langle f,\rho_s(g_t) \cdot f \rangle
  \geq |c|^2 \int_{\R \times \R} \frac{\dd x\dd y}{ (1+x^2)^{(1+s)/2}
  (1+y^2)^{(1+s)/2} |x|^{1-s}}>0.
  \end{equation}
Moreover, Ratner has proved in \cite[Theorem 1]{Rt} 
the upper bound $\limsup e^{t(1-s)}\langle f,\rho_s(g_t)
\cdot f \rangle <\infty$ (the convergence of the last integral in
\eqref{eq:ratner} to the integral in \eqref{eq:ratner2} can also be
verified directly).
This concludes the proof of the lemma.
\end{proof}

Since our Main Theorem implies exponential decay of correlations for
compactly supported smooth functions, it implies that the hypothesis of
Proposition \ref {ergodic}
is satisfied.  Corollary \ref {isolated} follows.


\begin{thebibliography}{MMY}



\bibitem[AF]{AF} Avila, A.; Forni, G.
Weak mixing for interval exchange transformations and translation flows.
Preprint (www.arXiv.org).  To appear in Annals of Math.

\bibitem[AV]{AV}
Avila, A.; Viana, M.
Simplicity of Lyapunov spectra: proof of the Zorich conjecture.
Preprint (www.arXiv.org).

\bibitem[Aar]{aaronson:book}
Aaronson, J. {\em An introduction to infinite ergodic theory},
volume~50 of {\em
  Mathematical Surveys and Monographs}.
American Mathematical Society, 1997.



\bibitem[At]{At}
Athreya, J.
Quantitative recurrence and large deviations for the Teichm\"{u}ller flow.
Preprint (www.arXiv.org).

\bibitem[BV]{BV}
Baladi, V.; Vall\'{e}e, B.
Exponential decay of correlations for
surface semi-flows without finite Markov partitions.  Proc. Amer. Math. Soc.
133  (2005),  no. 3, 865--874.

\bibitem[Bu]{Bu}
Bufetov, A.
Decay of correlations for the Rauzy-Veech-Zorich induction map on the space
of interval exchange transformations.  ESI Preprint.

\bibitem[Do]{Do}
Dolgopyat, D.
On decay of correlations in Anosov flows.  Ann. of Math.
(2)  147  (1998),  no. 2, 357--390.

\bibitem[EM]{EM} Eskin, A.; Masur, H. Asymptotic formulas on flat
surfaces.  Ergodic Theory Dynam. Systems  21  (2001),  no. 2, 443--478.

\bibitem[Fo]{F}
Forni, G.
Deviation of ergodic averages for area-preserving flows
on surfaces of higher genus.
Ann. of Math. (2) 155 (2002), no. 1, 1--103.

\bibitem[He]{H}
Hennion, H.
Sur un th{\'e}or{\`e}me spectral et son application aux noyaux
  lipschitziens.
Proc. Amer. Math. Soc. 118 (1993), 627--634.

\bibitem[K]{K}
 Kerckhoff, S. P. 
Simplicial systems for interval exchange maps and measured foliations. Ergodic Theory Dynam. Systems 5 (1985), no. 2, 257--271.


\bibitem[KZ]{KZ}
Kontsevich, M.; Zorich, A.
Connected components of the moduli spaces of Abelian differentials with
prescribed singularities. Invent. Math. 153 (2003), no. 3, 631--678.


\bibitem[MNS]{MNS}
Margulis, G. A.; Nevo, A.; Stein, E. M. Analogs of Wiener's ergodic theorems
for semisimple Lie groups. II.  Duke Math. J.  103  (2000),  no. 2,
233--259.


\bibitem[MMY]{MMY}
Marmi, S.; Moussa, P.; Yoccoz, J.-C. The cohomological equation
for Roth type interval exchange transformations. J. Amer.
Math. Soc. 18 (2005), 823-872.

\bibitem[Ma]{M}
Masur, H.
Interval exchange transformations and measured
foliations. Ann. of Math. (2) 115 (1982), no. 1, 169--200.

\bibitem[Rt]{Rt} Ratner, Marina
The rate of mixing for geodesic and
horocycle flows.  Ergodic Theory Dynam. Systems  7  (1987),  no. 2,
267--288.

\bibitem[R]{R}
Rauzy, G.
Echanges d'intervalles et transformations induites.
Acta Arith. 34, (1979), no. 4, 315--328.

\bibitem[Ve1]{V1}
Veech, W.
Gauss measures for transformations on the space of interval exchange maps.
Ann. of Math. (2) 115 (1982), no. 1, 201--242.

\bibitem[Ve2]{V2}
Veech, W.
The Teichm\"{u}ller geodesic flow.
Ann. of Math. (2) 124 (1986), no. 3, 441--530.

\bibitem[Z]{Z}
Zorich, A.
Finite Gauss measure on the space of interval exchange transformations.
Lyapunov exponents.
Ann. Inst. Fourier (Grenoble) 46 (1996), no. 2, 325--370.

\end{thebibliography}
\end{document}